\documentclass[12pt]{amsart}
\usepackage{amssymb}
\usepackage{amsmath}
\usepackage{amsopn}
\newtheorem{theorem}{Theorem}%[section]      %Usage: {\th Statement}    %Usage: {\th Statement}
\newtheorem{proposition}{Proposition}[section]           %  ETC ...
\newtheorem{corollary}[proposition]{Corollary}%[section]
\newtheorem{lemma}[proposition]{Lemma}%[section]

\newcommand{\CC}{{\mathbb C}}
\newcommand{\NN}{{\mathbb N}}
\newcommand{\PP}{{\mathbb P}}        %%% "OPEN" LETTERS %%%%

\newcommand{\ZZ}{{\mathbb Z}}

\def\rf#1{(\ref{#1})}
\def\ot{\otimes}
\def\a{\alpha}    \def\la{\lambda}
\def\b{\beta} 
\def\D{\Delta}
\def\g{{\mathfrak g}} \def\n{{\mathfrak n}}
 \def\H {{\mathfrak h}'}
\def\H {{\mathfrak h}}
\def\m{{\mathfrak m}}
 \def\bb {{\mathfrak b}}
\def\Dh{D}
\def\ov{\overline}
\def\End{\text{End}\,}
\def\Zp{Z^{\n_{}}}
\def\Zm{Z_{\n_-}}
\def\Sp{S^{\n_{}}}
\def\Sm{S_{\n_-}}
\def\P{{\rm \: \! P\:\!}}
\def\B{{\rm \: \! B\:\!}}
\def\BB{{\rm \: \! C\:\!}}

\def\Znn{{}_{\n_-}\A_\n}
\def\Ann{{}_{\n_-}{\A}_{\n}}

\def\Cnn{{}_{\n}{\A}_{\n_-}}
\def\Annl{{}_{\n_-}{\A}_{\n,\lambda}}
\def\Blnn{{}_{\lambda,\n_-}{\widetilde{\A}}_{\n}}
\def\Bnnl{{}_{\n_-}{\widetilde{\A}}_{\n,\lambda}}
\def\Alnn{{}_{\lambda,\n_-}{\A}_{\n}}

\def\nml{{\lambda}}

\def\PP{{\rm \: \! p\:\!}}
\def\PPP{\bar{\rm \: \! p\:\!}}

\def\zprim{{z}'}
\def\tzprim{{\tilde z}'}

\def\shift{\tau}
\def\T{T}
\def\N{\text{Nr}}
\def\A{{\mathcal A}}

\def\ad{\text{ ad\,}}

\def\e{\Hat{e}}

\def\Hh{\Hat{h}}

\def\q{{\mathrm q}}
\def\nq{\tilde{\q}}
\def\oq{\q}

\def\qq{\check{\q}}
\def\qqq{\breve{\q}}

\def\W{{\mathcal V}}
\def\w{ v}
\def\ow{\overline{w}}

\def\nnq{\n}
\def\ch{h}
\def\TT{\text{T}}
\def\AAA{{}}

\begin{document}
\bigskip
\hfill 
\begin{flushright}
CPT-P11-2006\\
ITEP-TH-13/06
\end{flushright}

\bigskip
\begin{center}
{\Large\bf Mickelsson algebras and Zhelobenko operators}

\bigskip
{\bf  S. Khoroshkin$^{\circ}$ \ \ and \, O. Ogievetsky$^{\star}$\footnote{On leave of absence
from P.N. Lebedev Physical Institute, Theoretical Department, Leninsky prospekt 53, 119991
Moscow, Russia}
%%$^{\star}$\footnote{E-mail: khor@itep.ru}
}\bigskip\\
$^\circ${\it Institute of Theoretical and Experimental Physics, 117259
Moscow, Russia}\smallskip\\
$^\star${\it Centre de Physique Th\'eorique\footnote{Unit\'e Mixte de Recherche
(UMR 6207) du CNRS et des Universit\'es Aix--Marseille I,
Aix--Marseille II et du Sud Toulon -- Var; laboratoire affili\'e \`a la FRUMAM (FR 2291)},
Luminy, 13288 Marseille, France}
\end{center}%\\
\bigskip
%\hfill{\today}
%\bigskip
\begin{center}
\textsc{Abstract}
\end{center}

{\footnotesize{
We construct a family of automorphisms of Mickelsson algebra, satisfying
braid group relations. The construction uses 'Zhelobenko cocycle' and
includes the dynamical Weyl group action as a particular case.
}}
\setcounter{tocdepth}{1}
{\small
\tableofcontents}
% $^\star$ {\it Institute of Theoretical \& Experimental Physics,
%117259 Moscow, Russia}

\section{Introduction}
Mickelsson algebras were introduced in \cite{M} for the study of
Harish-Chandra modules of  reductive groups.
The Mickelsson algebra, related to a real reductive group, acts in the
space of highest weight vectors of its maximal compact subgroup, and each
irreducible Harish-Chandra module of the initial reductive group
 is uniquely characterized by this action.

A similar construction takes place for any associative algebra $\A$, which
contains a universal enveloping algebra $U(\g)$ (or its $q$-analog)
of a contragredient Lie algebra
$\g$ with a fixed Gauss decomposition $\g=\n_-+\H+\n$.
Namely, we define the Mickelsson algebra $S^\n(\A)$ as the quotient
 of the normalizer $N(\A\n)$ of the ideal $\A\n$ over the ideal $\A\n$.
 In this case for any representation $V$ of $\A$ the  Mickelsson
algebra $S^\n(\A)$ acts in the space $V^\n$ of $\n$-invariant vectors.
This construction performs
a reduction of a representation of $\A$ over the action of $U(\g)$
 and can be viewed as a counterpart of hamiltonian reduction.
 It was applied for various problems of representation theory,
 see the survey \cite{T2} and references therein.

The structure of Mickelsson algebra %$S^\n(\A)$
 simplifies after the
localization over a certain multiplicative subset of $U(\H)$, where $\H$
is the Cartan subalgebra of $\g$. The corresponding algebra $Z^\n(\A)$
is generated by a finite-dimensional space of generators, which obey
quadratic-linear relations. These generators can be defined with the help
of an extremal projector of Asherova-Smirnov-Tolstoy \cite{AST}.
An application of the extremal projector to the study of the Mickelsson algebras
$Z^\n(\A)$ was developed by Zhelobenko \cite{Zh}.
Besides, Zhelobenko developed so called 'dual methods', where he gave a
construction of another generators of the Mickelsson algebra by means of
a family of special operators, which form a  cocycle on the Weyl group
\cite{Zh2}.

Later Mickelsson algebras regenerate in the theory of dynamical quantum
groups. Their basic ingredients, intertwining operators between Verma modules and their tensor
products with finite-dimensional representations actually form
special Mickelsson algebras. Matrix coefficients of
these intertwining operators are very useful in  quantum integrable
models \cite{ES}. The powerful instrument for the study of the algebra of intertwining
operators is a recurrence relation on the structure constants of the
algebra, known as ABRR equation, see \cite{ABRR}.

Tarasov and Varchenko \cite{TV} found the symmetries of the algebra of
intertwining operators,
which originate from  morphisms of Verma modules.
They satisfy braid group relations and transform the weights
by means of a shifted Weyl group action. These symmetries got the name
 'dynamical Weyl group'.  The theory of dynamical Weyl
groups was generalized to the quantum groups setup in \cite{EV}.

The form of  operators of the dynamical Weyl group is very close
to the factorized expressions for the extremal projector and for
 Zhelobenko cocycles \cite{Zh2}. However,
the precise statements and the origin of such a relation are not clear.
One of our goals  is to clarify this relation.

In this paper we establish a family of symmetries in a wide class of
Mickelsson algebras. They form a representation of the related braid group
by automorphisms of the Mickelsson algebra $\Zp(\A)$ and
transform Cartan elements by means of the shifted Weyl group action.
Each generating automorphism is a product of the Zhelobenko 'cocycle' map $q_{\a_i}$ and
an automorphism $T_i$ of the algebra $\A$, extending the action of the Weyl
group in $U(\g)$ (or Lusztig automorphism of $U_q(\g)$).
 The main new point of our approach is the homomorphism property of
Zhelobenko maps, that was not noticed before. Unfortunately, the proof of
this fact is not short and requires calculations with the extremal
projector.

The construction of automorphisms of Mickelsson algebra $\Zp(\A)$ is quite
general. In particular, it covers examples of Mickelsson algebras, related
to reductions $\g'\supset\g$ of one reductive Lie algebra to another,
and of the smash product $U(\g)\ltimes S(V)$ of $U(\g)$ and the symmetric
algebra of a $U(\g)$-module $V$, where it becomes the dynamical Weyl group
action after a specialization of Cartan elements of $\g$. It can be applied
for the construction of finite-dimensional representations of Yangians and
 of quantum affine algebras, see \cite{KhN,KhN2}.

The paper is organized as follows. In Section \ref{section2} we collect
a necessary information about the extremal projector and required extensions
of $U(\g)$.

 In Section \ref{section3} for a fixed contragredient
Lie algebra $\g$ of finite growth
we introduce a class of associative algebras $\A$, which
we call $\g$-admissible. They contain $U(\g)$ as a
 subalgebra and the adjoint action of $\g$ in $\A$ has special
 properties. In particular, $\A$ is equivalent to a tensor product of
$U(\g)$ and some subspace $\W\in\A$ as a $\g$-module with respect to
the adjoint action.
 Mickelsson algebras, related to admissible
algebras, have distinguished properties. The crucial one is the existence
of two distinguished spaces of generators $z_\w$ and $z'_\w$, $\w\in\W$.

Section \ref{section4} is an exposition of the 'Zhelobenko cocycle' \cite{Zh2}.
We present it with complete proofs in order to eliminate unnecessary
restrictions, assumed in \cite{Zh2}.
The story starts with the map $\q_\a$, which relates  universal Verma
modules, attached to different  maximal nilpotent subalgebras of $\g$.
The product of such operators over the system of positive roots maps
vectors $\w\in\W$ to the generators $z'_\w$ of Mickelsson algebras.
This invariant description proves the cocycle conditions for the maps
$\q_\a$.

Section \ref{section5} describes homomorphic properties of Zhelobenko
 maps. We prove first that the Zhelobenko
map $\q_\a$ establishes an isomorphism of a
double coset algebra and the Mickelsson algebra. This implies that
the compositions $\check{\q}_i$ of Zhelobenko maps with extensions of Weyl group
automorphisms are automorphisms of the Mickelsson algebras, satisfying braid
group relations.

In Sections \ref{section6} and \ref{section7}
 we calculate the images of generators of Mickelsson algebras  and of
 standard modules over them with respect to $\check{\q}_i$ and
 show that the dynamical Weyl group is a particular case of our
construction.
Section \ref{sectionzm} is devoted to the Mickelsson algebra $\Zm(\A)$,
related to $\n_-$-coinvariants of $\A$-modules.

Section \ref{section8} is a sketch of extensions of the
constructions to quantum groups $U_\nu(\g)$. The new important detail here
 is that compositions of Zhelobenko maps $\q_\a$ with Lusztig automorphisms
coincide with the compositions of $\q_\a$ with the adjoint action of Lusztig
automorphisms, see Proposition \ref{propnu1} and Proposition \ref{propnu4}.
This allows to prove both homomorphism properties and braid group
relations. We conclude  with remarks about the range of assumptions
on $\g$-admissible algebras, used in the paper.

\setcounter{equation}{0}
\section{Extremal projector}\label{section2}
In this section we review Zhelobenko's approach to extremal projector
of Asherova-Smirnov-Tolstoy \cite{AST}. The exposition follows \cite{Zh}
in  main details.

\subsection{Taylor extension of $U(\g)$}\label{section1.1}

%\addtocontents{toc}{\contentsline {section}{\numberline  \ref{section1.1} Taylor extension of $U(\g)$ }{\pageref{section1.1}}}

Let $\g$ be a contragredient Lie algebra of finite growth
 with symmetrizable Cartan matrix $a_{i,j}$, $i,j=1,...r$.
Let
\begin{equation}\label{1}
\g=\n_-+\H+\n_{}\,
\end{equation}
be its Gauss decomposition,
where $\H$ is a Cartan subalgebra, $\n_{}=\n_+\subset \bb_{}$ and
$\n_-\subset\bb_-$ are nilradicals
of two opposite Borel subalgebras $\bb_{}=\bb_+$ and $\bb_-$.
We use the notation $\Pi$ for the system of simple positive roots;
$\Delta_\pm$ and $\Delta=\Delta_+\coprod\Delta_-$
for the systems of positive, negative and all roots; $\Delta_\pm^{re}$ and
$\Delta^{re}=\Delta_+^{re}\coprod\Delta_-^{re}$
for the  systems of positive, negative and all real  roots.
Let $(,)$ be the scalar product in $\H^*$, such that
$(\a_i,\a_j)=d_ia_{i,j}=d_ja_{j,i}$, for $\a_i,\a_j\in\Pi$ and $d_i\in\NN$.

Denote by $Q\subset \H^*$ the root lattice, $Q=\ZZ\cdot\Delta$,
and put $Q_\pm=\ZZ_{\geq 0}\cdot\Delta_\pm$. For any
$\mu\in Q^\pm$ we denote by $U(\n_{})_\mu$ and  $U(\n_{-})_\mu$   the subspace
of elements $x$ of $U(\n_{})$ and $U(\n_{-})$, such that
$[h,x]= \langle \mu,h\rangle x$.
We accept the normalization of Chevalley generators $e_{\a_i}\in\n_{}$,
$e_{-\a_i}=f_{\a_i}\in\n_-$, and of coroots $\ch_{\a_i}=\a_i{}^\vee\in\H$,
where $\a_i\in\Pi$, such that
\begin{align*}
[e_{\a_i},e_{-\a_j}]&=\delta_{i,j}\ch_{\a_i}\,,&&\\
[\ch_{\a_i},e_{\pm\a_j}]&=\pm a_{i,j}e_{\pm\a_j}\,,&&\\
\text{ad}^{1-a_{i,j}}_{e_{\pm\a_{i}}}e_{\pm\a_j}&=0&\text{if}
\quad i\not= j\,& .
\end{align*}
{}For any $\gamma\in\Delta$  we define a coroot $h_\gamma\in\H$  by the rule
$$(\a,\a)h_\a+(\b,\b)h_\b=
(\a+\b,\a+\b)h_{\a+\b}\qquad\text{if}\quad\a,\b,\a+\b\in\Delta_+.$$

Let $W$ be the Weyl group of $\g$. For any $w\in W$ we denote by
$\T_w :U(\g)\to U(\g)$ a lift of the map $w:\H\to\H$ to the
 automorphism of the algebra $U(\g)$, satisfying braid group relations
$T_{ww'}=T_wT_{w'}$ if $l(ww')=l(w)+l(w')$, where $l(w)$ is the length of
$w$. For instance, we may choose $T_w$,  as in  \cite{L}.
 We accept a shortened
notation $\T_i$ for automorphisms $\T_{s_{\a_i}}$, where $\a_i\in\Pi$.

Denote by $\Dh$ the localization of the free commutative algebra $U(\H)$
with respect to the multiplicative set of denominators, generated by
$$ \{\ch_\a+k|\a\in\Delta\,,k\in\ZZ\}\,.$$
To any $\mu\in\H^*$ we associate an automorphism $\shift_\nu$
of the  algebra $\Dh$. It is uniquely defined by the conditions
\begin{equation}\label{shift}
\shift_\mu (h)=h+\langle h,\mu\rangle\qquad\text{for any}\quad h\in\H\,.
\end{equation}
Denote by $U'(\g)$ the extension of $U(\g)$ by means of $\Dh$:
$$U'(\g)=U(\g)\ot_{U(\H)}\Dh\approx\Dh\ot_{U(\H)}U(\g)\, .$$
Note that $U'(\g)$ is a $\Dh$-bimodule and any automorphism $\T_w$ admits
a canonical extension to an automorphism of $U'(\g)$, which we denote by
the same symbol.

Choose a normal ordering (see \cite{T} for the definition) $\gamma_1\prec\gamma_2\prec\ldots\prec\gamma_n$
of the system $\Delta_+$ of positive
roots of $\g$ ($n=|\Delta_+|$ may be infinite).
Let $e_{\pm\a}$, and $h_\a$, where $\a\in\Delta$ be Cartan-Weyl generators, constructed by recursive
procedure, attached to this order. We assume that they are normalized in
such a way that
\begin{equation}\label{P4}
[e_{\a},e_{-\a}]=\ch_{\a}\,, \qquad
[\ch_{\a},e_{\pm\beta}]=\pm\langle \ch_\a,\b\rangle
e_{\pm\beta}\, .
\end{equation}
%where $\langle h_\a,\b\rangle=\frac{2(\a,\beta)}{(\a,\a)}$.

For any $\ov{k}\in
\ZZ_{\geq 0}^n$,
$\ov{k}=(k_1,...,k_n)$ with $\sum k_i<\infty$ denote by $e^+_{\ov{k}}$ the
monomial $e^+_{\ov{k}}=e_{\gamma_1}^{k_1}\cdots
 e_{\gamma_n}^{k_n}\in U(\bb_+)$
and by $e^-_{\ov{k}}$ the monomial $e_{-\gamma_1}^{k_1}\cdots e_{-\gamma_n}^{k_n}\in
 U(\bb_-)$.
{}For every $\mu\in Q$ denote by
$\left(F_{\g,\n}\right)_\mu$ the vector space of series
\begin{equation}\label{P4a}
x_\mu=\sum\limits_{\ov{k},\ov{r}\in\ZZ_{\geq 0}^n}e^-_{\ov{k}}
x_{\ov{k},\ov{r}}
e^+_{\ov{r}}\,,\qquad x_{\ov{k},\ov{r}}\in \Dh\,
\end{equation}
of the total weight $\mu$. Set
$$F_{\g,\n}=\oplus_{\mu\in Q}\left(F_{\g,\n}\right)_\mu\,.$$

\begin{proposition} \label{prop1} {\rm (See \cite{Zh}, Section 3.2.3)}
The space $F_{\g,\n}$ is an associative algebra with respect to the
multiplication
of formal series. Its definition does not depend on a choice of the
normal ordering $\prec$.
\end{proposition}

Clearly, $F_{\g,\n}$ contains $U'(\g)$ as a subalgebra.
We call $F_{\g,\n}$ a {\it Taylor extension} of $U'(\g)$, related to
the decomposition \rf{1}.

A choice of the normal ordering is a technical tool for a description of the
algebra. It is used for a construction of particular bases in weight
components of the algebras $U(\n_\pm)$. Instead, one can fix, for any
$\nu\in Q_+$, a basis
$e^+_{\nu,j}$ of the finite-dimensional space $U(\n_{})_\nu$ and, for any
$\nu\in Q_-$, a basis
$e^-_{\nu,j}$ of the finite-dimensional space $U(\n_-)_\nu$. Then the
space $F_{\mu,\g}$ consists of formal series
$$ x_\mu=\sum\limits_{\nu\in Q_+,\ \nu'\in Q_-,\, j,j'}
e^-_{\nu',j'}x_{\nu',j',\nu,j}e^+_{\nu,j}\,,\qquad
x_{\nu',j',\nu,j}\in \Dh\,$$
of the total weight $\mu\in Q$.
\smallskip

\subsection{Universal Verma module and extremal projector}
\label{section1.2}
Set
$$M_{\n}(\g)=U'(\g)/U'(\g)\n_{}\, .$$
The space $M_{\n}(\g)$ is a left $U'(g)$-module and a $\Dh$-bimodule. It is called the
{\it universal Verma module}. Since $M_{\n}(\g)$ is a $U(\H)$-bimodule, we
have an adjoint action of $U(\H)$ in $M_{\n}(\g)$, defined as $ad_h(m)=[h,m]$
for any $h\in \H$ and $m\in M_{\n}(\g)$. We have the weight decomposition
of $M_{\n}(\g)$ with respect to the adjoint action of $U(\H)$:
$$M_{\n}(\g)=\oplus_{\mu\in Q_-}\left(M_{\n}(\g)\right)_\mu\, .$$
Denote by $E'(\g)$ the algebra of
all endomorphisms of $M_{\n}(\g)$, which commute with the right action of
$U(\H)$. We have a linear map $\xi: U'(\g)\to E'(\g)$, induced by the
multiplication in $U'(\g)$ and establishing in
$M_{\n}(\g)$ the structure of the left $U'(\g)$-module.
For any $\mu\in Q$, define
$$E_\mu(\g)=\{a\in E'(\g)\,|\,[\xi(h),a]=\langle \mu,h\rangle
a\quad\text{for any } h\in\H\}\, ,$$
and set
$$E(\g)=\oplus_{\mu\in Q}E_\mu(\g)\, .$$

\begin{proposition} \label{prop2} {\rm (See \cite{Zh}, Sections 3.2.4.-3.2.5)}
The map $\xi$ induces an isomorphism of algebras:
\begin{equation}
\label{P1}
\xi: F_{\g,\n}\to E(\g)\, .
\end{equation}
\end{proposition}

The proof of Proposition \ref{prop2} is strongly based on the
nondegeneracy of Shapovalov form.

Recall \cite{Sh} that the Shapovalov form $A(x,y): U(\g)\ot U(\g)\to U(\H)$
is defined by the relation $A(x,y)=\beta(x^t\cdot y)$, where
$x\mapsto x^t$ is the Chevalley antiinvolution
($e_{\pm\a}^t=e_{\mp\a},\ h^t=h,\ (xy)^t=y^tx^t$)
 in $U(\g)$ and
$\beta: U(\g)\to U(\H)$ is the projection with respect to the decomposition
$$U(\g)=U(\H)\oplus \left( \n_-U(\g)+U(\g)\n_{}\right)\,.$$
The Shapovalov form vanishes on the left ideal $U(\g)\n_{}$ and is
nondegenerate on $U(\n_-)$. It admits an extension to a nondegenerate form
on $M_{\n}(\g)$ with values in $\Dh$.

For the proof of the isomorphism \rf{P1} we choose for any
$\nu\in Q^-$ a basis
$e_{\nu,j}^-$ of  $U(\n_-)_\nu$, which is
 orthogonal with respect to the Shapovalov form, and take
 $e_{-\nu,j}^+=(e_{-\nu,j}^-)^t$.

\begin{proposition} \label{prop3} {\rm (See \cite{Zh}, Section 3.2.8)}
There exists a unique element $\P_\n\in F_{\g,\n}$, satisfying equations
\begin{equation}\label{P2}
e_\a \P_\n= \P_\n e_{-\a}=0\qquad\text{for all }\quad\a\in\Delta_+
\end{equation}
with zero term $\beta(\P_\n)=1$. It is a self-adjoint projector of zero weight,
$\P_\n^2=\P_\n$, $\P_\n^t=\P_\n$.
\end{proposition}

The definition and the construction of the projector $\P_\n$ depends on a
choice of the nilpotent subalgebra $\n\subset\g$. When $\n$ coincides with a
fixed nilpotent subalgebra $\n$, entering the decomposition \rf{1}, we omit
the label $\n$ and denote the projector simply by $\P$.

Note  an important property of $\P$, which follows from Proposition
 \ref{prop3}:
\begin{equation}\label{P2a}
\P-1\in F_{\g,\n}\n_{}\bigcap \n_-F_{\g,\n}\ .
\end{equation}
By means of Proposition \ref{prop2}, the element $\P$ is described as an element of $E(\g)$, which projects
the universal Verma module $M_{\n}(\g)$ to the subspace
$\left(M_{\n}(\g)\right)^{\n_{}}=M_{\n}(\g)_0$
of $\n_{}$-invariants along 
$\n_-M_{\n}(\g)=\oplus_{\gamma<0}\left(M_{\n}(\g)\right)_\gamma$. The element $\P$ is called the
{\it extremal projector}. It was discovered in \cite{AST}.
\medskip

There are three distinct cases of the use of the extremal projector.
\begin{itemize}
\item[1.] Let $V$ be a $F_{\g,\n}$-module. Then $\P$ projects $V$ on the subspace
$V^{\n_{}}$ of $\n_{}$-invariants along the subspace $\n_-V$.

\item[2.] Let $V$ be a module over $U'(\g)$, locally finite with respect to
$\n_{}$. Then  $\P$ projects $V$ on the subspace
$V^{\n_{}}$ of $\n_{}$-invariants along the subspace $\n_-V$.

\item[3.] Let $V$ be a module over $U(\g)$, locally finite with respect to
$\bb_{}$. Assume that $\mu\in\H^*$ satisfies the conditions:
\begin{equation}\label{P3}
\langle \mu+\rho, \ch_\a\rangle\not=-1,-2,...\qquad \text{for any
}\a\in\Delta_+\, ,
\end{equation}
where $\rho=\frac{1}{2}\sum_{\a\in\Delta_+}\a$.

 Denote by $V_\mu$ the generalized weight subspace of $V$ of the weight
$\mu$. Then $\P$ projects $V_\mu$ on $V_\mu\cap V^{\n_{}}$ along 
$V_\mu\cap \n_-V$.
\end{itemize}

In the following, for any left $U(\g)$-module $M$, in which the action
of the projector $\P$ is defined, we denote the corresponding
element of $\End M$ by $\PP$, and for any right $U(\g)$-module $N$,
in which the action
of the projector $\P$ is defined, we denote the corresponding
element of $\End N$ by $\PPP$.

The operator $\PP$ satisfies the relations
\begin{equation}\label{?1}
\PP(e_{-\gamma}m)=e_{\gamma}\PP(m)=0\qquad \text{for any}\quad \gamma\in
\Delta_+,m\in  M,\qquad \PP^2=\PP\,.
\end{equation}
The operator $\PPP$ satisfies the relations
\begin{equation}\label{?5}
\PPP(ne_{\gamma})=\PPP(n)e_{-\gamma}=0\qquad \text{for any}
\quad \gamma\in \Delta_+,n\in N\qquad \PPP^2=\PPP\,.
\end{equation}

\subsection{Multiplicative formula for extremal projector}
The extremal projector $\P$ for simple Lie algebras was discovered and
investigated by Asherova, Smirnov and Tolstoy \cite{AST}. They presented
a multiplicative expression for $\P$, which was later generalized to
affine Lie superalgebras  and their $q$-analogs. We reproduce here
the formula of \cite{AST}.

{}For any $\a\in\Delta_+$ and $\lambda\in\H^*$
let $f_{\a,n}[\lambda]$ and $g_{\a,n}[\lambda]$ be the following elements
of $\Dh$:
\begin{equation}
\label{P5}
f_{\a,n}[\lambda]=\prod\limits_{j=1}^n(h_\a+\langle
h_\a,\lambda\rangle+j)^{-1},\qquad
g_{\a,n}[\lambda]=\prod\limits_{j=1}^n(-h_\a+\langle
h_\a,\lambda\rangle+j)^{-1}.
\end{equation}
Define $\P_{\a}[\lambda]\in F_{\g,\n}$ and $\P_{-\a}[\lambda]\in F_{\g,\n_-}$
by relations
\begin{equation}
\label{P6}
\P_\a[\lambda]=\sum\limits_{n=0}^\infty
\frac{(-1)^n}{n!}f_{\a,n}[\lambda]e_{-\a}^ne_\a^n,\qquad
\P_{-\a}[\lambda]=\sum\limits_{n=0}^\infty
\frac{(-1)^n}{n!}g_{\a,n}[\lambda]e_{\a}^ne_{-\a}^n
\end{equation}
Set
\begin{equation}\label{P6a}
\P_\a=\P_\a[ \rho],\qquad \P_{-\a}=\P_{-\a}[ \rho]\ .
\end{equation}

\begin{proposition}\label{prop4} {\rm (See \cite{AST}.)} Let
$\gamma_1\prec...\prec \gamma_n$ be a normal ordering of $\Delta_+$. Then
the extremal projector $\P$ is equal to the product
\begin{equation}
\label{P7}
\P=\prod\limits_{\gamma\in\Delta_+}^\prec \P_\gamma\,,
\end{equation}
where the order in the product coincides with the chosen normal order
$\prec$.
\end{proposition}

Analogously, the product
$\prod\limits_{\gamma\in\Delta_+}^\prec \P_{-\gamma}$ is equal in
$F_{\g,\n-}$ to the projector $\P_{\n_-}$.

For a generalization of \rf{P7} to arbitrary contragredient Kac-Moody Lie
algebras of finite growth and their $q$-analogs, see \cite{T,KhT}.

We define also  elements $\P[\lambda]$ and  ${\P}_{-}[\lambda]$
for any $\lambda\in\H^*$ by the
relations
\begin{equation}\label{P8}
\P[\lambda]=\prod\limits_{\gamma\in\Delta_+}^\prec \P_\gamma[\lambda],\qquad
{\P}_{-}[\lambda]=
\prod\limits_{\gamma\in\Delta_+}^\prec \P_{-\gamma}[\lambda]
\end{equation}
It is known (see \cite{Zh}, and also  Section \ref{section2.8} of the present paper),
 that $\P[\lambda]$
and ${\P}_{-}[\lambda]$ do not depend on a choice of the normal order. In this notation,
$\P=\P[\rho]$ and $\P_{\n_-}={\P}_{-}[\rho]$.
%\normalize

\setcounter{equation}{0}
\section{Mickelsson algebras}\label{section3}

\subsection{$\g$-admissible algebras}\label{section2.1}

Let $\A$ be an associative algebra, which contains $U(\g)$ as a subalgebra.
Then $\A$ has a natural structure of a $U(\g)$-bimodule. Since $U(\g)$ is a Hopf
algebra, the bimodule structure produces an adjoint action of $U(\g)$ in
$\A$. We have
\begin{equation}\label{ad}
\ad_g(x)=\sum_ig'_ixS(g''_i)\, ,
\end{equation}
 where the coproduct
$\Delta(\g)$ of the element
$g\in U(\g)$ has a form $\Delta(g)=\sum_i g'_i\ot g''_i$ and $S(g)$ is the
antipode of $g$.
 The adjoint action $\ad_g$ of an element $g\in\g$ is the
commutator, %\begin{equation}\label{ad1}
$\ad_g(a)=ga-ag$\,.
%\end{equation}

%We also use the opposite adjoint action
%$\tilde{\ad}$ of $U^{op}(\g)$ in $\A$, defined by the relation
%\begin{equation}\label{ad2}
%\tilde{\ad}_g(a)=\sum_iS(g'_i)xg''_i\,.
%\end{equation}
% The adjoint action $\tilde{\ad}_g$ of an element $g\in\g$ looks as the
%commutator
%%\begin{equation}\label{ad3}
%$\tilde{\ad}_g(a)=ag-ga\,$.
%%\end{equation}

In the sequel we use the following notation for the adjoint action:
$$\Hat{g}(a)\equiv \ad_g(a)\, ,%%\qquad \tilde{g}(a)\equiv \tilde{\ad}_g(a),
\qquad g\in U(\g), a\in\A\,.$$

  We call $\A$ a
{\it $\g$-admissible } algebra, if
\begin{itemize}
\item[(a)] there is a subspace $\W\subset \A$, invariant with respect
to the adjoint action of $U(\g)$, such that the multiplication $m$ in $\A$
induces the following isomorphisms of vector spaces
 \begin{equation*}
{\rm (a1)}\qquad m: U(\g)\ot \W\to\A\,,\qquad {\rm (a2)}\qquad
 m: \W \ot U(\g)\to\A\, ;
\end{equation*}
\item[(b)] the adjont action of real root vectors $e_\gamma\in U(\g)$
 in $\W$ is locally nilpotent.
The adjoint action of the Cartan subalgebra $\H$ in $\W$ is semisimple.
% Here $\gamma\in\Delta^{re}$.
\end{itemize}
Sometimes we call $\W$ an {\it $ad$-invariant generating subspace} of
the $\g$-admissible algebra $\A$.

The  condition (a) says, in particular, that  $\A$ is a free left
 $U(\H)$-module  and a free right $U(\H)$-module.

Since the adjoint action of real root vectors in $U(\g)$ is locally nilpotent,
the conditions (a) and (b) imply that the adjoint action of real root
vectors  $e_\gamma$ is
locally nilpotent in $\A$, that is, for any $a\in \A$ and
$\gamma\in\Delta^{re}$
 vectors $\ad^n_{e_\gamma}(a)$ are zero for sufficiently big $n$.
 Thus the restriction
of the adjoint action of $U(\g)$ to any ${\mathfrak sl}_2$-subalgebra, generated by
real root vectors $e_{\pm \gamma}$, where $\gamma\in\Delta_+^{re}$,
 is locally finite in $\A$. Since
  the adjoint action of the Cartan subalgebra $U(\H)$ is semisimple,
$\A$ admits the weight decomposition with respect to the adjoint action
of $U(\H)$.

There are two main classes of $\g$-admissible algebras.

\begin{itemize}
\item[1]. Let $\g$ be a reductive finite-dimensional Lie algebra and
$\g_1$ a contragredient Lie algebra of finite growth, which contains $\g$.
Then the adjoint action of $\g$ in $\g_1$ is locally finite and $\A=U(\g_1)$ is
a $\g$-admissible algebra.

\item[2]. Let $\W$ be a  $U(\g)$-module algebra with a locally nilpotent action
of real root vectors. This means that
$\W$ is an associative algebra, equipped with a structure of a
$U(\g)$-module, such that the action of real root vectors
is locally finite. These two structures are related: the action of Lie
generators $g\in \g$ satisfy the Leibniz rule:
\begin{equation}\label{b31}
g(\w_1\w_2)=g(\w_1)\w_2+\w_1g(\w_2)\, .
\end{equation}
Denote by  $U(\g)\ltimes \W$ the smash product of $U(\g)$ and
$\W$. It is an associative algebra, generated by elements $g\in U(\g)$
and $\w\in \W$, satisfying the relation
\begin{equation}
\label{R1}
g\w-\w g=g(\w)\, ,\qquad g\in\g,\qquad \w\in \W\, ,
\end{equation}
and, more generally,
\begin{equation}
\sum_i g'_i \w S(g''_i)=g(\w)\, ,\qquad g\in U(\g),\qquad \w\in \W\, ,
\label{R2}
\end{equation}
where $S$ is the antipode in $U(\g)$, $\Delta(g)=\sum_i g'_i\ot g''_i$
is the comultiplication in $U(\g)$.
The smash product $U(\g)\ltimes \W$ is a $\g$-admissible algebra.
\item[2a]. Let $V$ be a  $U(\g)$-module and $\End^0V$ be the
algebra of the endomorphisms of $V$, finite with respect to the adjoint action
of real root vectors of $U(\g)$. Then the tensor product $U(\g)\ot \End^0V$ is a $\g$-admissible
algebra. This construction is a particular case of the previous one: the tensor
product $U(\g)\ot \End^0V$ is a smash product of $\CC\ot \End^0V$ and of
diagonally imbedded $U(g)$, generated by the elements $g\ot 1+1\ot g$,
$g\in\g$.
\end{itemize}

\subsection{Mickelsson algebras. Definitions}\label{section2.2}

Let $\A$ be an associative algebra, which contains $U(\g)$.
Let $\N(\A\n_{})$ be the normalizer of the left ideal
$\A\n_{}$ :
$$x\in \N(\A\n_{})\equiv \n_{}x\subset \A\n_{} .$$
 Denote by $S^{\n_{}}(\A)$  the quotient space
$$S^{\n_{}}(\A)=\N(\A\n_{})/\A\n_{}
\,.$$

%%%%%workspace%%%%%%%%%%%%%%%%
\begin{proposition}
\label{prop5}${}$
\begin{itemize}
\item[(i)] The space $S^{\n_{}}(\A)$
 is an algebra with respect to the multiplication in $\A$;
\item[(ii)] Let $M$ be an $\A$-module. Then the space $M^{\n_{}}$ of
$\n_{}$-invariant vectors in $M$ is a $S^{\n_{}}(\A)$-module.
\end{itemize}
\end{proposition}
\smallskip

The algebra $S^{\n_{}}(\A)$ is called
{\it Mickelsson algebra} \cite{M}. Since $\H$ normalizes $\n_\pm$, we have
the inclusion
$U(\H)\subset S^{\n_{}}(\A)$.
 Denote by $\A'$ the localization
$$\A'=\Dh\ot_{U(\H)}\A\, .$$
By the condition (a) of a $\g$-admissible algebra,
 we have a canonical imbedding of $\A$ into $\A'$ and thus an adjoint
action of $U(\g)$ in $\A'$, compatible with the adjoint action of $U(\g)$ in
$\A$.

 Define
{\it Mickelsson algebra} $\Zp(\A)$  as the quotient
$$\Zp(\A)=\N(\A'\n_{})/\A'\n_{}\,,$$
where $\N(\A'\n_{})$ is the normalizer of the left ideal $\A'\n_{}$ of $\A'$.
The algebra  $\Zp(\A)$  is a localization of the algebras
 $\Sp(\A)$:
$$\Zp(\A)=\Dh\ot_{U(\H)}\Sp(\A).$$
 We can change the order of taking quotients and subspaces in the definition
of Mickelsson algebra. Then the Mickelsson algebra $Z^{\n_{}}(\A)$ is
defined as a subspace of $\n$-invariants in a left $U(\g)$-module
$M_\n(\A')=\A'/\A'\n$  :
\begin{equation*}%\label{R16}
\Zp(\A)=\left(M_{\n}(\A')\right)^{\n_{}}
=\{m\in M_{\n}(\A')\,|\,\n_{}m=0\}\, .
\end{equation*}

The algebra $\Zp$ acts in the space $M^{\n_{}}$ of $\n_{}$-invariants of
any $\A'$-module $M$.

\subsection{Double coset algebra}\label{sectiondcs}
Suppose that a $g$-admissible algebra $\A$  satisfies the additional
 {\it local highest weight} condition:
\medskip

\begin{itemize}
\item[(HW)]\ For any $\w\in\W$, the adjoint action of elements
$x\in U(\n)_\mu$ on $\w$ is nontrivial, $\Hat{x}(\w)\not=0$,
 only for a finite number of
$\mu\in\H^*$.
\end{itemize}
\smallskip

With this assumption  the
quotient $M_\n(\A')=\A'/\A'n$ has a structure of a left
$F_{\g,\n}$-module, extending the
 action of $\A'$ by left multiplication. In particular, the extremal
projector $\P$ acts in the left $F_{\g,\n}$-module $M_\n(\A')$.

The properties of the extremal projectors imply the relation
\begin{equation}\label{?2}
\Zp(\A)= \text{Im}\,\PP\subset M_\n(\A')\,,
\end{equation}
where $\PP\in\End M_\n(\A')$ is the action of $\P$ in $M_\n(\A)$, see 
Section \ref{section1.2}.

Denote by $\Ann$ the double coset space
\begin{equation}\label{Znn}
\Ann=\n_-\A'\backslash \A'/\A'\n\equiv\A'/\left(\n_-\A'+\A'\n\right)\,.
\end{equation}
Equip $\Ann$ with a binary operation $\circ:\Ann\ot\Ann\to\Ann$:
\begin{equation}
\label{?3}
a\circ b=a\P b\ \stackrel{def}{=}a\ \cdot\PP(b)\ .
\end{equation}
The rule \rf{?3} means the following. For a class $\bar{x}$ in $\Ann$,
we take its representative $x\in \A'$. For a class $\bar{y}$ in $\Ann$, we
take its representative $y\in M_\n(\A')$. Consider an element $x\PP(y)$ in
$M_\n(\A')=\A'/\A'\n$.
Then  its class modulo $\n_-M_\n(\A')$ defines an element
$\bar{x}\circ\bar{y}$ of $\Ann$. It does
not depend on a choice of representatives.

We call the double coset space $\Ann$, equipped with the operation \rf{?3}
the {\it double coset algebra} $\Znn$.

Define  linear maps $\phi^+:\Zp(\A)\to\Znn$ and $\psi^+:\Znn\to\Zp(\A)$
by the rules
\begin{equation}\label{?4}
\phi^+(x)=x\mod\n_-M_\n(\A'),\qquad \psi^+(y)=\PP(y),\qquad x\in\Zp(\A),y\in\Znn\ .
\end{equation}
Let us explain the formula  $\psi^+(y)=\PP(y)$. For a class
$y$ in $\Znn=\A'/\n_-\A'+\A'\n$, we choose its
representative $\bar{y}\in M_\n(\A')=\A'/\A'\n$ and take $\PP(\bar{y})$.
The result does not depend on a choice of a representative and is denoted
by $\psi^+(y)$.

\begin{proposition}\label{propgener1}
Assume that a $\g$-admissible algebra $\A$ satisfies the condition (HW). Then
\begin{itemize}
\item[(i)]
The operation \rf{?3} equips $\Znn$ with a structure of an associative
algebra.
\item[(ii)]
The linear maps $\phi^+$ and $\psi^+$ are inverse to each other
and establish an isomorphism of algebras $\Zp(\A)$ and $\Znn$.
\end{itemize}
\end{proposition}

{\it Proof}. Let $x\in \Zp(\A')$. Then
$\PP(x)=x\mod \A'\n$
due to \rf{P2a}. Thus $\psi^+\cdot\phi^+=Id_{\Zp(\A)}$. On the other hand,
due to the same property of $\P$, for any $y\in\A'/\A'\n$,
$$\PP(y)=y \mod \n_-\A'+\A'\n\ .$$
Thus $\phi^+\cdot\psi^+=Id_{\Znn}$. So the maps $\phi^+$ and $\psi^+$ are
inverse to each other.

Let now $x\in\A'$ and $y\in\A'$ be representatives of classes $\tilde{x}$ and
$\tilde{y}$ in $\Znn$, $\bar{x}=x\mod\A'\n$ and $\bar{y}= y\mod\A'\n$ be
their images in $M_\n(\A')$. We have $\psi^+(\tilde{x})=\PP(\bar{x})$,
$\psi^+(\tilde{y})=\PP(\bar{y})$ and $\psi^+(\tilde{x}\circ \tilde{y})=
\PP(x\cdot \PP(\bar{y}))$.

On the other hand, the multiplication rule $m$ in $\Zp(\A)$ can be written
as follows. Let ${z},{u}\in\Zp(\A')$. Let ${z}'\in\N(\A'\n)$ be a
representative of a class $z\in \A'/\A'\n$.
 Then $m(z\ot u)={z}'\cdot{ u}$ as an element
of $M_\n(\A')$. By \rf{?2}, $m(z,u)=\PP({z}'\cdot {u})$.
Thus we have in $\Zp(\A')$
$$m(\PP(x)\ot\PP(y))= \PP({\PP(x)}'\cdot \PP(y))\,,$$
where ${\PP(x)}'$ is a representative of a class $\PP(x)$ in $\A'$.
By the property \rf{P2a},
$${\PP(x)}'=x+x'+x''\,,$$
where $x'\in \n_-\A'$ and $x''\in\A'\n$. Thus
$$\PP({\PP(x)'}\cdot \PP(y))=\PP(x\cdot\PP(y))$$
due to \rf{P2a} and \rf{?1}. Thus $\psi^+$ is a homomorphism, which proves
simultaneously (i) and (ii). \hfill{$\square$}

\subsection{Generators of Mickelsson algebras}
\label{section2.6}

Let $\A$ be a $\g$-admissible algebra with an ad-invariant generating
subspace $\W$, satisfying the highest weight condition (HW). By the condition
(a) of a $\g$-admissible algebra (see  Section \ref{section2.1}) and the PBW theorem
for the algebra $U(\g)$, any element of
$\A'$ can be uniquely presented in the following form
$$x=\sum\nolimits_if_id_ie_i\w_i,\qquad \text{where}\quad f_i\in U(\n_-),\, d_i\in\Dh,\, e_i\in U(\n),
\, \w_i\in\W\,.$$
Due to the highest weight condition (HW), we can move all $e_i$ to the right
and get a presentation
$$x=\sum\nolimits_if'_id'_i\w'_ie'_i,\qquad \text{where}\quad f'_i\in U(\n_-),\,
 d'_i\in\Dh,\, e'_i\in U(\n),
\, \w'_i\in\W\, .$$
In the double coset space this presentation gives
$$x=\sum\nolimits_id'_i\w'_i,\quad \mod \n_-\A'+\A'\n\,,
\qquad \text{where}\quad
 d'_i\in\Dh,\,
 \w'_i\in\W\,.$$

\begin{proposition}
\label{propr3}
Let $\A$ be a $\g$-admissible algebra satisfying the highest weight condition (HW). Then
\begin{itemize}
\item[(i)] Each element of the double coset algebra $\Znn$ can be uniquely
presented in a form $x=\sum_id_i\w_i$, where $d_i\in\Dh$, $\w_i\in\W$, so
that $\Znn$ is a free left (and right) $\Dh$ module, isomorphic to
$\Dh\ot\W$ ($\W\ot\D$).
\item[(ii)]
{}For each $\w\in \W$ there exists a unique element $z_\w\in \Zp(\A)$
 of the form
\begin{equation}
\label{r7}z_\w= \w+\sum_{i=1}d_i f_i \w_i\ ,\qquad
f_{i}\in \n_-U(\n_-),\ d_{i}\in \Dh, \ \w_i\in \W
\end{equation}
such that the algebra $\Zp(\A)$ is a free left (and right) $\Dh$-module,
generated by the elements $z_\w$. The element $z_\w$ is equal to $\PP(\w)$.
\end{itemize}
\end{proposition}

{\it Proof}. The part (i) is already proved. Applying Proposition
\ref{propgener1}, we see that any element of $\Zp(\A)$ can be presented in a
form $\sum_i d_i\PP(w_i)$, where $d_i\in\Dh$, $\w_i\in\W$. The element
$z_{\w}=\PP(\w)$ has a form \rf{r7} due to the definition of the operator $\PP$
and is uniquely characterized by this presentation. \hfill{$\square$}

The Mickelsson algebras have distinguished generators of another type.
Their existence is imposed by the following proposition.

Let $\A$ be an arbitrary $\g$-admissible algebra with an ad-invariant generating
subspace $\W$.

\begin{proposition}${}$
\label{theorem00}
{}For each $\w\in \W$ there exists
 at most one element $z'_\w\in\Zp(\A)$ of the form
\begin{align}
\label{r8}\zprim_\w&= \w+\sum\nolimits_{i}d_i\w_i   f_i\ ,
&&f_{i}\in \n_-U(\n_-),\ d_{i}\in \Dh, \ \w_i\in \W\ .
\end{align}
\end{proposition}
\medskip

{\it Proof}. Consider the case of the algebra $\Zp(\A)$.
 If $m\in M_\n(\A')=\A'/\A'\n$ is a highest weight
vector, that is $e_\a m=0$ for any $\a\in\Delta_+$, then $dm$ is also a
highest weight vector for any $d\in \Dh$. Thus (i) is equivalent
to the statement
that for any $\gamma\in\H^*$ there is no highest weight vector
of the form
\begin{equation}\label{rm1}
x=\sum\nolimits_i  d_i \w_i f_i,
 \qquad f_{i}\in \n_-U(\n_-),\ d_{i}\in \Dh, \ \w_i\in \W,
\end{equation}
where all the terms have the weight $\gamma$ with respect to the adjoint
representation of $\H$. In other words, we should prove that the conditions
 $[e_{\a_i},x]=0$
in $M_{\n}(\A)$ imply $x=0$ in $M_{\n}(\A)$ if all $f_{j}\in \n_-U(\n_-)$.

By the condition (a2) of a $\g$-admissible algebra and the PBW theorem for $U(\g)$
the elements $\w_i f_i$ form a basis of $M_\n(\A')$ over $\Dh$ if
$ \w_i$ form a basis of $\W$ and $f_i$ form a basis of $U(\n_-)$.
Consider the terms of the right hand side of \rf{rm1}
 with $\w_j$ having minimal weights with respect
to other weights which occur in \rf{rm1}. Then the expression
$[e_{\a_i},x]$ contain terms $\w_j[e_{\a_i},f_j]$ which are nonzero
for some $\a_i$ if $f_{j}\in \n_-U(\n_-)$. This is because all the highest
weight vectors of $M_\n(\g)$ have zero weight. Thus $x$ cannot be a highest
weight vector.
\hfill{$\square$}
%Since any element of $M_\n(\A')$ can be presented in a form
%$\sum\nolimits_i d_i \w_i f_i$, where
%$ f_{i}\in U(\n_-),\ d_{i}\in \Dh, \ \w_i\in \W$, the statement (i) implies
%the statement (ii) for the algebra $\Zp(\A)$.

\medskip

%Assume now that both conditions (HW) and (LW) are satisfied. This takes
%place, for instance, 
We now specify to a case when $\A$ is an admissible algebra over
a finite-dimensional reductive Lie algebra.% $\g$.

\begin{theorem} \label{proposition3a} Let $\g$ be a
finite-dimensional reductive Lie algebra and
$\A$ a $\g$-admissible algebra with generating subspace $\W$.
Then for any $\w\in \W$ there exists
a unique element $\zprim_\w\in \Zp(\A)$ \rf{r8}.
The algebra $\Zp(\A)$ is generated by elements
$\zprim_\w$ as a free left (and right)
 $\Dh$-module.
\end{theorem}

Proof of  Theorem \ref{proposition3a} will be given in the next Section.

\subsection{ Relations between two sets of generators}
\label{section2.8}

Extend the notation of canonical generators of Mickelsson algebras
to the elements of $\Dh\ot_\CC\W$ and $\W\ot_\CC\Dh$. We set for any $d\in\Dh$ and
$\w\in\W$
\begin{equation}\begin{split}\label{r10}
z_{d\ot\w}&=d\cdot z_\w,\quad z'_{d\ot\w}=d\cdot z'_\w,\quad
z_{\w\ot d}= z_\w\cdot d,\quad z'_{\w\ot d}=z'_\w\cdot d
\quad\text{in}\quad
Z^\n(\A)\,.
\end{split}
\end{equation}

{}Fix a positive real root $\a$.
 We define now certain operators in a vector space
$\W\ot\Dh$.  Accept the notation $A^{(1)}$ for the operator $A\ot 1$ in
a vector space $\W\ot \Dh$ and $A^{(2)}$ for the operator $1\ot A$.

Let $\a$ be a real  root, For any $\mu\in\H^*$ and $n\geq 0$ define
operators $\bar{f}_{n,\a}^{(1)}[\mu]$,
$\bar{g}_{n,\a}^{(1)}[\mu] $, $\B_{\a}^{(1)}[\mu]$ and
$\BB^{(1)}_{-\a}[\lambda]$
  $\in$ $ \End(\W\ot \Dh)$ by the relations
\begin{equation}\begin{split}\label{r11}
\bar{f}_{n,\a}^{(1)}[\mu]&=\prod_{k=1}^{n}
\left(\Hat{h}_\a^{(1)}+h^{(2)}_\a+\langle h_\a,\mu\rangle+k
\right)^{-1},\\
\bar{g}_{n,\a}^{(1)}[\mu]&=\prod_{k=1}^{n}
\left(-\Hat{h}_\a^{(1)}-h^{(2)}_\a+\langle h_\a,\mu\rangle+k
\right)^{-1},
\end{split}
\end{equation}
\begin{equation}%\label{r12}
\begin{split}
\B_{\a}^{(1)}[\mu]&=\sum_{n=0}^\infty
\frac{(-1)^n}{n!}\bar{f}_{n,\a}^{(1)}[\mu]\big(\e_{-\a}^{(1)}\big)^n
\big(\e_\a^{(1)}\big)^n \!, \\
\BB_{-\a}^{(1)}[\mu]&=\sum_{n=0}^\infty
\frac{(-1)^n}{n!}\bar{g}_{n,\a}^{(1)}[\mu]\big(\e_{\a}^{(1)}\big)^n
\big(\e_{-\a}^{(1)}\big)^n \!.
\end{split}
\end{equation}
Here $\Hat{h}_\a^{(1)}=\ad_{h_\a}^{(1)}$ is the adjoint action of $h_\a$ in
$\W$, $h^{(2)}_\a$ is the operator of multiplication by $h_\a$ in $\Dh$.

 Define also
operators $\bar{f}_{n,\a}^{(2)}[\mu]$,
$\bar{g}_{n,\a}^{(2)}[\mu] $, $\BB_{\a}^{(2)}[\mu]$ and
$\B^{(2)}_{-\a}[\lambda]$
  $\in$ $ \End(\Dh\ot)$ by the following relations
\begin{equation}\begin{split}\label{r11a}
\bar{f}_{n,\a}^{(2)}[\mu]&=\prod_{k=1}^{n}
\left(\Hat{h}_\a^{(2)}-h^{(1)}_\a+\langle h_\a,\mu\rangle+k
\right)^{-1},\\
\bar{g}_{n,\a}^{(2)}[\mu]&=\prod_{k=1}^{n}
\left(-\Hat{h}_\a^{(2)}+h^{(1)}_\a+\langle h_\a,\mu\rangle+k
\right)^{-1},
\end{split}
\end{equation}
\begin{equation}\label{r12}
\begin{split}
\BB_{\a}^{(2)}[\mu]&=\sum_{n=0}^\infty
\frac{(-1)^n}{n!}\bar{f}_{n,\a}^{(2)}[\mu]\big(\e_{-\a}^{(2)}\big)^n
\big(\e_\a^{(2)}\big)^n \!, \\
\B_{-\a}^{(2)}[\mu]&=\sum_{n=0}^\infty
\frac{(-1)^n}{n!}\bar{g}_{n,\a}^{(2)}[\mu]\big(\e_{\a}^{(2)}\big)^n
\big(\e_{-\a}^{(2)}\big)^n \!.
\end{split}
\end{equation}
Here $\Hat{h}_\a^{(2)}=\ad_{h_\a}^{(2)}$ is the adjoint action of $h_\a$ in
$\W$, $h^{(1)}_\a$ is the operator of multiplication by $h_\a$ in $\Dh$.

Let $\g$ be a reductive finite-dimensional Lie algebra. Let $\prec$ be a
normal ordering of the system $\Delta_+$ of positive roots. Set
\begin{equation*}\begin{split}%\label{r13}
\B^{(1)}[\lambda]=\prod\limits_{\gamma\in\Delta_+}^\prec
\B_{\gamma}^{(1)}[\lambda],\qquad
\BB^{(1)}_-[\lambda]=\prod\limits_{\gamma\in\Delta_+}^\prec
\BB_{-\gamma}^{(1)}[\lambda]\, ,\\
\BB^{(2)}[\lambda]=\prod\limits_{\gamma\in\Delta_+}^\prec
\BB_{\gamma}^{(2)}[\lambda],\qquad
\B^{(2)}_-[\lambda]=\prod\limits_{\gamma\in\Delta_+}^\prec
\B_{-\gamma}^{(2)}[\lambda].
\end{split}
\end{equation*}

\begin{theorem}\label{theorem0}%${}$
Let $\A$ be an admissible algebra
 with a generating subspace $\W$ over a finite-dimensional reductive
Lie algebra $\g$. Then
 for any $\w\in \W$ we have the following
equality in $Z^{\n}(\A)$
\begin{equation}\label{r14}
z_\w=z'_{\B^{(1)}[\rho](\w\ot 1)}\,.
\end{equation}
\end{theorem}
In particular, operators $\B^{(1)}[\rho](\w):\W\ot\Dh\to\W\ot\Dh$
 do not depend on a choice of
the normal order $\prec$.

{\it Proof}. Consider the left $F_{\g,\n}$-module $M_\n(\A)=\A'/\A'\n$.
The multiplication in $\A'$ induces an isomorphism of
vector spaces $M_\n(\A)$ and  $\W\ot M_\n(\g)$, where
$M_\n(\g)=U'(\g)/U'(\g)\n$:
\begin{equation}\label{r15}m:\W\ot M_\n(\g)\to M_\n(\A)\,. \end{equation}
With this identification the tensor product
$\W\ot M_\n(\g)$ becomes a $F_{\g,\n}$-module. As a $U(\g)$-module it
coincides with the tensor product of $\W$, equipped with a structure of
the adjoint representation of $U(\g)$, and of the left $U(\g)$-module
$M_\n(\g)=U'(\g)/U'(\g)\n$.
This follows from the Leibniz rule
 $$g(\w\cdot x)=(g\w-\w g)\cdot x+\w\cdot (gx)=
\Hat{g}(\w)\cdot x+\w\cdot gx\,,$$
for any $g\in\g$, $\w\in\W$ and $x\in U'(\g)/U'(\g)\n$.

The elements of $\Dh$ act by the following rule:
 for any $d\in\Dh$, $\w\in\W$ of the weight $\mu_\w$
and $x\in M_\n(\A)$ we have $d\cdot (\w\ot x)=\w\ot \tau_{\mu_\w}(d)x$.
Due to local finiteness of the adjoint action in $\W$
these prescriptions define correctly the action of $F_{\g,\n}$ in
 $\W\ot M_\n(\g)$.

Under the identification \rf{r15} we have
\begin{equation*}%\label{r16}
z_\w=\PP(\w\ot 1)\,.
\end{equation*}
In order to express $z_\w$ via $z'_\w$, we should write it in a form
$\w'_i\ot d_i +\text{lower order terms}$, where 'lower order terms'
contain vectors, whose second tensor component lies in
 $\Dh\cdot\n_-U(\n_-)$.
Write $\P$ as a series over ordered monomials in
$e_{\gamma_i}$, and $e_{-\gamma_i}$, where $\gamma_i\in\Delta_+$, with coefficients
being rational functions of $h_{\gamma_i}$ such that in
any monomial all $e_{-\gamma_i}$ stand before all $e_{\gamma_j}$, in accordance to
the rules of $F_{\g,\n}$,
$$\P=\P(h_{\gamma_i},e_{-\gamma_i},e_{\gamma_i})\ .$$
 By a coproduct rule,
in the action of $\P$ in  $\W\ot M_\n(\g)$ we substitute, instead of
$e_{\pm \gamma_i}$, the sum $\e_{\pm \gamma_i}^{(1)}+e_{\pm \gamma_i}^{(2)}$,
and, instead of
$h_{ \gamma_i}$, the sum $\Hat{h}_{ \gamma_i}^{(1)}+h_{ \gamma_i}^{(2)}$,
each term acting in the corresponding tensor component. The action of
$e_{ \gamma_i}^{(2)}$ on $\w\ot 1$ vanishes, the action of
$e_{- \gamma_i}^{(2)}$ gives 'lower order terms'.
So the term we are looking for is equal to
$$\PP\left(\left(\Hat{h}^{(1)}_{\gamma_i}+ {h}^{(2)}_{\gamma_i}\right),
\e_{-\gamma_i}^{(1)},\e_{\gamma_i}^{(1)}\right)(\w\ot 1)
\,.$$
This is precisely $\B^{(1)}[\rho](\w\ot 1)$. \hfill{$\square$}

The operators $\B^{(1)}[\rho]$, $\B_-^{(2)}[\rho]$,
 $\BB^{(2)}[-\rho]$ and $\BB_-^{(1)}[-\rho]$,
  are closely related to the
operators $\P[\lambda]$ and ${\P}_-[\lambda]$, see \rf{P8}.
Namely, denote by $\rho^{(1)}$ and $\rho^{(2)}$ the expressions
$$\rho^{(1)}=\frac{1}{2}\sum_{\gamma\in\Delta_+} h_\gamma^{}\ot
\gamma,\qquad
\rho^{(2)}=\frac{1}{2}\sum_{\gamma\in\Delta_+}\gamma\ot h_\gamma^{}\,.$$
Then
\begin{equation*}%\label{r16a}
\begin{split}
\B^{(1)}[\rho]=\Hat{\PP}^{(1)}[\left(\rho+\rho^{(2)}\right)],\qquad
\BB_-^{(1)}[-\rho]=\Hat{{\PP}}_-^{(1)}[-\left(\rho+\rho^{(2)}\right)],\\
\B_-^{(2)}[\rho]=\Hat{{\PP}}_-^{(2)}[\left(\rho+\rho^{(1)}\right)],\qquad
\BB^{(2)}[-\rho]=\Hat{\PP}^{(2)}[-\left(\rho+\rho^{(1)}\right)].
\end{split}
\end{equation*}

The changes of coordinates, described in  Theorem \ref{theorem0} and 
Theorem \ref{theorem0a} below can be inverted by means of the relation
\begin{equation}\label{r20}
\P_\a[\lambda]P_{-\a}[-\lambda]=\frac{h_\a+\langle h_\a,\lambda\rangle}
{\langle h_\a,\lambda\rangle}\ .
\end{equation}
This relation makes sense for a generic $\lambda$
in any finite-dimensional representation of
${\mathbf sl}_2$, see  \cite{TV}, Theorem 10.

So we have 
\begin{equation}\label{R21}
\begin{split}
\B^{(1)}[\rho]^{-1}&=\BB_-^{(1)}[-\rho]\prod_{\a\in\Delta_+}
\frac{h_\a^{(2)}+\langle h_\a,\rho\rangle}
{\Hat{h}_\a^{(1)}+{h}_\a^{(2)}+\langle h_\a,\rho\rangle},\\
\B_-^{(2)}[\rho]^{-1}&=\BB^{(2)}[-\rho]\prod_{\a\in\Delta_+}
\frac{h_\a^{(1)}+\langle h_\a,\rho\rangle}
{{h}_\a^{(1)}-\Hat{h}_\a^{(2)}+\langle h_\a,\rho\rangle}\ .
\end{split}
\end{equation}

{\it Proof of  Theorem \ref{proposition3a}}. The inversion relations
\rf{R21} imply that
\begin{equation}
\label{inv2}
z'_\w=z_{\gamma_1{\BB}_-^{(1)}[-\rho](1\ot\w)},
%\tilde{z}'_\w=
%\tilde{z}_{\gamma_2{\BB}^{(2)}[-\rho](\w\ot 1)}\,,
\end{equation}
where $\gamma_1=\prod_{\a\in\Delta_+}
\frac{h_\a^{(2)}+\langle h_\a,\rho\rangle}
{\Hat{h}_\a^{(1)}+{h}_\a^{(2)}+\langle h_\a,\rho\rangle}$.
%and $\gamma_2=\prod_{\a\in\Delta_+}
%\frac{h_\a^{(1)}+\langle h_\a,\rho\rangle}
%{{h}_\a^{(1)}-\Hat{h}_\a^{(2)}+\langle h_\a,\rho\rangle}$.
Thus we have correctly defined elements $z'_\w$,
which proves the first statement of the Theorem. Other statements  follow
from the corresponding statements of Proposition \ref{propr3}.
 \hfill{$\square$}
\medskip

{\bf Remark}.
 If we replace the ring $\Dh$ by the field of fractions
$\tilde{D}={\rm Frac}(U(\H))$, the statement of  Theorem \ref{proposition3a}
 does not require the precise inversion relation \rf{inv2}. We just note that
all the elements $\w_i$, entering the right hand side of \rf{r7},
belong to a finite-dimensional ad-invariant subspace $V\subset\W$, generated
 by $\w$. We move then all $f_i$ to the right and get a sum of elements
$z'_{\w_k}$ with coefficients in $\Dh$. If we allow coefficients
in $\tilde{D}$, such a transformation defines an injective operator in
 $\tilde{D}\ot V$, which is a finite-dimensional vector space over
 $\tilde{D}$. Thus this operator  is invertible.
 An advantage of the precise formula  \rf{inv2} is that it shows that
the coefficients of the inverse matrix belong to $\Dh$.

\setcounter{equation}{0}
\section{Zhelobenko maps}\label{section4}
\subsection{Maps $\q_{\a}^{(k)}$ and $\q_\a$}
\label{section4.1}
Let $\a$ be a  real root of $\g$. For any $x\in\A$ and  $k\geq 0$
denote by $\q_{\a}^{(k)}(x)$ the following element of
$\A'/\A'e_\a$:
\begin{equation}
\label{q1}
\oq_{\a}^{(k)}(x)=\sum\limits_{n=k}^\infty
\frac{(-1)^n}{(n-k)!}\e_\a^{n-k}(x)\cdot e_{-\a}^n\cdot{g_{n,\a}}\
 \quad \mod\ \A'e_\a,
\end{equation}
where
\begin{equation}
\label{q2}
g_{n,\a}=\left(h_\a(h_\a-1)\cdots (h_\a-n+1)\right)^{-1}\ .
\end{equation}
The assignment \rf{q1} has the following properties \cite{Zh}:
\begin{equation}
\begin{split}
%\begin{align}
\label{q4}
(i)\qquad
\oq_{\a}^{(k)}(xe_{-\a})&=0\ ,%\notag
\\
(ii)\qquad
[h,\oq_\a^{(k)}(x)]&= \oq_\a^{(k)}([h,x])\ ,\quad\qquad\qquad h\in\H\, ,
\\
(iii)\qquad\quad
\oq_\a^{(k)}(xh)&= \oq_\a^{(k)}(x)(h+\langle h,
\a\rangle)\ ,\qquad h\in\H\, ,\\
(iv)\qquad\quad
 e_\a \oq_\a^{(k)}(x)&= \oq_\a^{(k)}(e_\a x)=-k\oq_{\a}^{(k-1)}(x)\ .%\notag
\end{split}\end{equation}%\end{align}
Note that the quotient $\A'/\A'e_{\a}$ admits left and right actions of
$\Dh$, so the commutator $[h,\oq_\a^{(k)}(x)]$ is well defined in (ii).
The property (iii) in the notation \rf{shift} can be written as
\begin{equation}
\label{q4a}
\oq_\a^{(k)}(xh)=\oq_\a^{(k)}(x)\shift_\a(h),\qquad h\in U(\H)\, .
\end{equation}
We extend the assignment \rf{q1} to the linear map
$\q_{\a}^{(k)}:\A'\to\A'/\A'e_\a$ by means of the relation $(iii)$
and denote by $\q_\a$ the linear map $\q_{\a}^{(0)}:\A'\to\A'/\A'e_\a$,
\begin{equation}
\label{q5}
\q_{\a}(x)=\sum\limits_{n=0}^\infty
\frac{(-1)^n}{n!}\e_\a^{n}(x)\cdot e_{-\a}^n\cdot{g_{n,\a}}\ \quad \mod\
\A'e_\a\ .
\end{equation}
The relations \rf{q4}, (i), (iv) show that the image of $\q_{\a}$ belongs to the
 normalizer of $\A'e_\a$ in $\A'$ and the ideals $e_\a\A'$ and $\A'e_{-\a}$ are
 in the kernel of $\q_\a$.

Consider the one-dimensional vector space $\CC e_\a$ as an abelian
 Lie algebra $\n_\a=\CC e_\a$ and one-dimensional vector space $\CC e_{-\a}$
 as an abelian Lie algebra $\n_{-\a}=\CC e_{-\a}$.
{}Following the notation of  Section \ref{section2.2},
denote by $Z^{\n_\a}(\A')=\N(\A'e_\a)/\A'
e_\a$ the Mickelsson algebra, related to the reduction to the corresponding
$sl_2$-subalgebra.

\begin{proposition}
\label{propqa}
The map $\q_\a$ defines an isomorphism of vector spaces
${}_{\n_{\a}}\A_{\n_{-\a}}$ $\equiv$
$e_{\a}\A'\backslash \A'/\A'e_{-\a}$ and  $ Z^{\n_\a}(\A')$,
such that for any $x\in {}_{\n_{\a}}\A_{\n_{-\a}} $, and
$d\in\Dh$
\begin{align}\label{q30}
[d,\q_\a(x)]&=\q_\a([d,x]),&\q_\a(xd)&=\q_\a(x)\shift_\a(d)
\,.
%\\ \label{q27a}
%\q_\a(w)&=z'_w,& \w&\in \W\,.
\end{align}
\end{proposition}
\smallskip

{\it Proof.} First of all note that the properties \rf{q4}, (i) and (iv),
say that the map $\q_\a$ vanishes on $\n_\a\A'+\A'\n_{-\a}$ and thus defines a map
of ${}_{\n_{\a}}\A_{\n_{-\a}}$ to  $\A'/\A'\n_{\a}$. Its image
belongs to $\Zp(\A)$ by  \rf{q4}, (iv).

 Let $\g_\a$ be the $sl_2$ subalgebra of $\g$, generated by $e_\a$,
$e_{-\a}$ and $h_\a$.

Since $\a$ is a real root, the adjoint action of $\g_\a$ in $\g$ is locally
finite and semisimple. So there is a decomposition $\g=\g_\a+\mathfrak{p}$,
invariant with respect to the adjoint action of $\g_\a$. Poincare-Birkhoff-Witt
theorem implies that the multiplication in $U(\g)$ defines an isomorphism
of tensor products $U(\g_\a)\ot S(\mathfrak{p})$ and
$ S(\mathfrak{p})\ot U(\g_\a)$ with $U(\g)$. Here $S(\mathfrak{p})$ is
regarded as a subspace of $U(\g)$, which consists of symmetric
noncommutative polynomials on $\mathfrak{p}$. The space $S(\mathfrak{p})$
is invariant with respect
to the adjoint action of $\g_\a$. Thus $U(\g)$ is $\g_\a$-admissible.
Since the adjoint representation of $U(\g_\a)$ in
$U(\g)$ is locally finite, and $\A$ is a $\g$-admissible algebra, it is
$\g_\a$-admissible as well. Let $\W_\a\subset\A$ be the subspace of $\A$,
invariant with respect to the adjoint action of $\g_a$ and the
multiplication in $\A$ induces an isomorphism of vector spaces
$U(\g_\a)\ot\W_\a$ and $\A$.

The double coset space $ {}_{\n_{\a}}\A_{\n_{-\a}}$ is a free
$\Dh$-module, generated by the vector space $\W_\a$ (see the formulas in Section
\ref{section2.6} with a replacement of $\n$ by $\n_{-\a}$). On the other hand,
the Mickelsson algebra $ Z^{\n_\a}(\A')\,$ is also a free $\Dh$-module,
generated, in the notation of the Remark in  Section \ref{section2.6},
 by the vectors $z'_{\n_\a,\w}$, $\w\in\W_\a$, see Theorem
\ref{proposition3a}. By the structure of the map $\q_\a$, see \rf{q1}, we
have
$$
\q_\a(\w)=\w+\sum \w_i f_i d_i,
$$
where
$ \w_i\in\W_\a$, $f_i$ is a polynomial on $e_{-\a}$ without  a constant term,
$d_i\in \Dh$. In other words,
\begin{equation}\label{q27a}
\q_\a(\w)=z'_{\n_\a,\w}\,.
\end{equation}
 The relations \rf{q27a} and \rf{q4a} prove the proposition. \hfill{$\square$}

Let us restrict the map $\q_\a$ to the normalizer $\N(\A'\n_{-\a})$. Due to
\rf{q4}, (i), this restriction defines a map
$\q_\a|_{\N(\A'\n_{-\a})}: Z^{\n_{-\a}}(\A)\to Z^{\n_\a}(\A)$.
We have also a map in other direction,
$\q_{-\a}|_{\N(\A'\n_{\a})}: Z^{\n_{\a}}(\A)\to Z^{-\n_\a}(\A)$.

\begin{proposition}\label{propinverse}
We have equalities
\begin{equation}\label{inverse}
\begin{split}
\q_{-\a}\q_\a(x)&=(h_\a+1)^{}x(h_\a+1)^{-1}\qquad \text{for any}\quad x\in
Z^{\n_{-\a}}(\A)\,,\\
\q_{\a}\q_{-\a}(y)&=(h_\a+1)^{-1}y(h_\a+1)^{}\qquad \text{for any}\quad y\in
Z^{\n_\a}(\A)\,.
\end{split}
\end{equation}
\end{proposition}

In particular, the restriction of the map $\q_\a$ to the normalizer
$\N(\A'e_{-\a})$ defines
an isomorphism of the vector spaces $Z^{\n_{-\a}}(\A)$ and $Z^{\n_{\a}}(\A)$.
The inverse is given by the formula
$$y\mapsto (h_\a+1)^{-1}\q_{-\a}(y)(h_\a+1)^{},\qquad y\in Z^{\n_{\a}}(\A)\,.$$
{\it Proof}.  Take $y\in Z^{\n_\a}(\A')$. We have
\begin{equation*}\begin{split}
\q_{\a}\q_{-\a}(y)&=\q_{\a}\left(\sum_{n\geq 0}\frac{(-1)^n}{n!}
\e_{-\a}^n(y)e_{\a}^n\cdot
g_{n,-\a}\right)\\
&=\q_{\a}\left(\sum_{n\geq 0}\frac{1}{n!}
\e_{-\a}^n(y)e_{\a}^n\right)\cdot
\left((h_\a+2)(h_\a+3)\cdots (h_\a+n+1)\right)^{-1}\\
&=\q_{\a}\left(\sum_{n\geq 0}\frac{(-1)^n}{n!}
\e_{\a}^n\e_{-\a}^n(y)\right)\cdot
\left((h_\a+2)(h_\a+3)\cdots (h_\a+n+1)\right)^{-1}\\
&=\sum_{m,n\geq 0}\frac{(-1)^{n+m}}{n!m!}
\e_{\a}^m\e_{\a}^n\e_{-\a}^n(y)\cdot e_{-\a}^m
\left((h_\a+2)\cdots (h_\a+n+1)\right)^{-1}\cdot g_{m,\a}\ .
\end{split}
\end{equation*}
Since  $\q_{\a}\q_{-\a}(y)\in Z^{\n_\a}(\A)$, we can replace it by
$\PP_\a\left(\q_{\a}\q_{-\a}(y)\right)$,
where $\P_\a$ is the projection operator $\P$, related
to the $sl_2$-algebra $\g_\a$, and  $\PP_\a$ is the action of $\P_\a$ in
 $\A'/\A'\n_\a$.
  Since $\PP_\a(e_{-\a}z)=0$ for
any $z\in\A'/\A'\n_\a$, we have
\begin{equation*}
\begin{split}
&\q_{\a}\q_{-\a}(y)=\PP_\a\left(\sum_{m,n\geq 0}\frac{(-1)^{n}}{n!m!}
\e_{-\a}^m\e_{\a}^m\e_{\a}^n\e_{-\a}^n(y)\left((h_\a+2)\cdots
 (h_\a+n+1)\right)^{-1}\cdot g_{m,\a}\right)\\
=&\PP_\a\left(\sum_{m,n\geq 0}\frac{(-1)^{n}}{n!m!}
\prod_{k=1}^m(\Hh_\a-h_\a+k-1)^{-1}
\prod_{k=1}^n(-\Hh_\a+h_\a+k+1)^{-1}
\e_{-\a}^m\e_{\a}^m\e_{\a}^n\e_{-\a}^n(y)\right).
\end{split}
\end{equation*}
The expression inside brackets can be interpreted as
$$\hat{\PP}_\a^{(2)}[-h_\a^{(1)}\ot\rho-\rho]\cdot
\hat{\PP}_{-\a}^{(2)}[h_\a^{(1)}\ot\rho+\rho](1\ot y)$$
in $\Dh\ot Z^{\n_\a}(\A)$. Due to \rf{r20}, it is equal to
$\frac{\Hh_\a-h_\a^{(1)}-1}{-h_\a^{(1)}-1}(1\ot y)=
\frac{-\Hh_\a+h_\a^{(1)}+1}{h_\a^{(1)}+1}(1\ot y)\,.$
It means that
\begin{equation*}\begin{split}
\q_{\a}\q_{-\a}(y)=&\PP_\a\left(\frac{-\Hh_\a+h_\a+1}{h_\a+1}y\right)=
\PP_\a\left((h_\a+1)^{-1}y(h_\a+1)\right)\\
=&(h_\a+1)^{-1}\PP_\a(y)(h_\a+1)=(h_\a+1)^{-1}y(h_\a+1)\,.
\end{split}\end{equation*}
In the last line we used again the relation $y=\PP_\a(y) \mod \A'\n_\a$
for any  $y\in Z^{\n_\a}(\A)$.
The second relation is proved in an analogous manner.
\hfill{$\square$}

\subsection{Maps $\q_{\a,{\m}}^{(k)}$ and $\q_{\a,{\m}}$}

Let $\a$ be a  real  root of $\g$, $e_\a$ the corresponding root vector
with respect to the decomposition \rf{1}. Let $\m$  be a maximal nilpotent
subalgebra of $\g$, conjugated to $\n$ by means of an element of
the Weyl group $W$,
 such that $e_\a$ is a simple positive  root vector of
$\m$. Set $\m_-=\m^t$, where $x\mapsto x^t$ is the Chevalley antiinvolution
in $U(\g)$, see  Section \ref{section1.2}.
 For any  $x\in\A$ and $k\geq 0$ denote by $\oq_{\a,\m}^{(k)}(x)$ the
following element of $\A'/\A'\m_{}$:
\begin{equation}
\label{q1a}
\oq_{\a,\m}^{(k)}(x)=\sum\limits_{n=k}^\infty
\frac{(-1)^n}{(n-k)!}\e_\a^{n-k}(x)\cdot e_{-\a}^n\cdot{g_{n,\a}}\
 \quad \mod\ \A'\m_{},
\end{equation}
where $g_{n,\a}$ is given by the relation \rf{q2}. In other words,
\begin{equation}\label{q26}
\oq_{\a,\m}^{(k)}(x)=\pi_{m,\a}\cdot \oq_{\a}^{(k)}(x)\, ,
\end{equation}
where $\pi_{m,\a}:\A/\A e_\a\to \A/\A \m$ is the natural factorization map.
Due to \rf{q26}, the assignment \rf{q1a} satisfies the relations \rf{q4}
and admits an extension to the map
$\q_{\a,\m}^{(k)}:\A'\to \A'/\A'\m$, satisfying the relation
\begin{equation}\label{q8}
[d,\q_{\a,\m}^{(k)}(x)]=\q_{\a,\m}^{(k)}([d,x]),\qquad
\q_{\a,\m}^{(k)}(xd)=\q_{\a,\m}^{(k)}(x)\shift_\a(d)\, ,
\end{equation}
 for any $x\in\A'$
and $d\in\Dh$. We denote by $\q_{\a,\m}^{}$ the map
$\q_{\a,\m}^{(0)}:\A'\to \A'/\A'\m$.

{}For any element $w\in W$, where $W$ is the Weyl group of $\g$, and a maximal nilpotent
subalgebra $\m\subset\g$, we denote by
 $\m_{}^w\subset \g$ the nilpotent subalgebra
$\m_{}^w=\T_w(\m_{})$, see  Section \ref{section1.1}.

\begin{lemma}\label{lemma1} For any $k\geq 0$,
\begin{equation*}%\label{q6}
\q_{\a,\m}^{(k)}(\A'\m_{}^{s_\a})=0\ .
\end{equation*}
\end{lemma}

{\it Proof}. Denote by $\m(\a)\subset\m$ the parabolic subalgebra of $\m$,
generated as a vector space by all root vectors of $\m$ except $e_\a$.
 Due to \rf{q4}, (i), it is
sufficient to prove that $\q_{\a,k}(\A\m(\a))=0$. The basic theory of
root systems for simple Lie algebras says, that
$\Hat{e}_{\pm\a}(\m(\a))\subset \m(\a)$. Thus for any $n\geq 0$ we have
$\Hat{e}^n(\A\m(\a))\subset \A\m(\a)$ and
$\m(\a) e_{-\a}^ng_{n,\a}^{-1}\subset \A'\m(\a)$, which imply the statement of
the Lemma. \hfill{$\square$}

Due to Lemma \ref{lemma1}, the maps $\q_{\a,\m}^{(k)}$ induce linear maps
of $\A'/\A'\m_{}^{s_\a}$ to $\A'/\A'\m_{}$. We denote them by
the same symbol:
$$\q_{\a,\m}^{(k)}: \A'/\A'\m_{}^{s_\a} \to \A'/\A'\m_{}\,.$$

%%%%%%%%%%%%

\begin{proposition} ${}$
\label{prophom1}
\begin{itemize}\item[(i)]
The map $\q_{\a,\m}$ transforms the normalizer $\N(\A'\m^{s_\a})$ to
 the Mickelsson
algebra $Z^m(\A)$ $=$ $\N(\A'\m)/\A'\m$.
\item[(ii)] The restriction of the map $\q_{\a,\m}$ to the normalizer
$\N(\A'\m^{s_\a})$ defines an isomorphism of vector spaces
$Z^{m^{s_\a}}(\A)$ and $Z^m(\A)$, satisfying \rf{q8} with $k=0$.
\end{itemize}
\end{proposition}

{\it Proof}. We prove first the statement (i). Let $\Delta_+(\m)$ be
a system of positive roots, related to the decomposition
$\g=\m_-+\H+\m$. By assumption, $\a$ is a simple root of $\Delta_+(\m)$.

Let $x$ be an element of the normalizer of the ideal $\A'\m^{s_\a}$,
$x\in \N(\A'\m^{s_\a})$. It means that $e_{-\a}x\in \A'\m^{s_\a}$  and
$e_\gamma x\in \A'\m^{s_\a}$
for any root $\gamma\in\Delta_+(\m)$, $\gamma\not=\a$.

Since $e_\a \q_{\a,\m}=0$ by \rf{q4}, (iv), we have to prove that
\begin{equation}\label{hom8}
e_\gamma \q_{\a,\m}(x)=0\qquad \text{for any}\qquad \gamma\in \Delta_+(\m)
\setminus\a\,.
\end{equation}
{}Fix a positive root $\gamma\in \Delta_+(\m)\setminus\a$.
Let $\gamma_0,...,\gamma_m$ be an '$\a$-string' of roots, starting with
$\gamma$,
that is $\gamma_0=\gamma$, and $\gamma_{k+1}=\gamma_k+\a$. Since $\a$ is simple, all
roots $\gamma_k$ are positive, $\gamma_k\in \Delta_+(\m)\setminus\a$ and
 for any $k= 0,\ldots,m$ we have
\begin{equation}\label{q22}
\Hat{e}_{\a}^k(e_{\gamma})=a_ke_{\gamma_k},\quad k=0,\ldots,m,\
 a_k\in\CC\,,\quad\qquad \Hat{e}_{\a}^k(e_{\gamma})=0,\quad k>m\, .
\end{equation}
{}For any $y\in\A$ we have by \rf{q22}
$$
\q_{\a}(e_{\gamma}y)=\sum_{k=0}^ma_k\sum_{n=k}^\infty\frac{(-1)^n}{(n-k)!}
e_{\gamma_k}\Hat{e}_{\a}(y)e_{-\a}^n\cdot g_{n,\a}\, ,
$$
therefore, 
\begin{equation}\label{q24}
e_{\gamma}\q_{\a}(y)=\q_{\a}(e_{\gamma}y)-
\sum_{k=1}^ma_k e_{\gamma_k}\q_{\a}^{(k)}(y)\, .
\end{equation}
Iterations of \rf{q24} and a factorization by $\A'\m$ give the relation
\begin{equation}
\label{hom9}
e_\gamma\q_{\a,\m}(x)=\q_{\a,\m}
(e_{\gamma}x)+
\sum_{k=1}^mb_k \q_{\a,\m}^{(k)}(e_{\gamma_k}
x)\,,\qquad b_k\in\CC .
\end{equation}
The right hand side of \rf{hom9} is zero by assumption.
This proves \rf{hom8} and the statement (i) of the Proposition.

Let us prove (ii). The root $-\a$ is a simple positive root for the algebra
$\m^{s_\a}$. Thus by (i) the map $\q_{-\a,\m^{s_\a}}$ maps the normalizer
$\N(\A'\m)$ and the Mickelsson algebra $Z^{\m}(\A)$
to the Mickelsson algebra $Z^{\m^{s_\a}}(\A)$.
This implies that the map
$$\q'_{-\a,\m^{s_\a}}={\rm Ad}_{h_\a+1}^{-1}\cdot\q_{-\a,\m^{s_\a}}\,,$$
which sends $x\in \N(\A'\m)$ to $(h_\a+1)^{-1}\cdot
\q_{-\a,\m^{s_\a}}(x)\cdot(h_\a+1)$, also maps
 $Z^{\m}(\A)$ to the Mickelsson algebra $Z^{\m^{s_\a}}(\A)$. Proposition
\ref{propinverse} says that $\q'_{-\a,\m^{s_\a}}$ is inverse to the
map of $Z^{\m^{s_\a}}(\A)$ to $Z^{\m}(\A)$, induced by the restriction of
 $\q_{\a,\m^{}}$ to $\N(\A'\m^{s_\a})$. This proves the statement (ii).
 \hfill{$\square$}

\subsection{Maps $\q_{\ow,\m}$}\label{section3.3}

Let $\m$ be a maximal nilpotent subalgebra of $\g$, conjugated to $\n$ by
an automorphism $\T_{w'}$, where $w'\in W$, $\m=\T_{w'}(\n)$,
 and $\m_-=\T_{w'}(\n_-)$  the opposite maximal
nilpotent subalgebra. Assume that $w\in W$ satisfies the condition
\begin{equation}\label{q27}
\dim \T_w(\m)\cap \m_-=l(w)\, ,
\end{equation}
where $l(w)$ is the length of $w$ in $W$. Since $\m=\T_{w'}(\n)$,
the condition \rf{q27}
is equivalent to the  relation $l(w'w)=l(w')+l(w)$.
Let $\overline{w}$ be a reduced decomposition of the element $w\in W$:
\begin{equation}\label{q19}
%\overline{w}:\
 \ow=\{w=s_{\a_{i_1}}s_{\a_{i_2}}\cdots s_{\a_{i_n}}\}\,.
\end{equation}
 Let $n=l(w)$. Define a sequence $w_1,w_2$, $\ldots,$ $w_{n}, w_{n}=w$
 of  elements of the Weyl group $W$, a sequence
 $\gamma_1,\gamma_2,...,\gamma_n$
 of positive roots,
 and a sequence $\m_0,...,\m_n$ of maximal
nilpotent subalgebras of $\g$  by the following
inductive rule :
\begin{align}\label{q19a}
 w_0&=1, \quad &w_{k+1}&=w_k\cdot s_{\a_{i_{k+1}}}\,,\\
%\gamma_1&=w'(\a_{i_1}),&
\label{q19b}\gamma_{k}&=w'w_{k-1}^{}(\a_{i_{k}})\, ,&
    %%\text{\rm{and}}&
\m_{k}&=\T_{w'w_{k-1}}(\n)\, .
\end{align}
The relations \rf{q27}, \rf{q19a}--\rf{q19b} imply that for any
$k=1,\ldots,n$ the root vector $e_{\gamma_k}$ is a simple positive root
vector for the algebra $\m_k$, and the composition
\begin{equation}
\label{q10}
\q_{\overline{w},\m}=\q_{\gamma_1,\m_1}\cdot \q_{\gamma_2,\m_2}\cdot\ldots\cdot
\q_{\gamma_n,\m_{n}}\ 
\end{equation}
is a well defined map
 $\q_{\overline{w},\m}: \A'/\A'\m_{}^{w}\to\A'/\A'\m_{}$.
The index $\overline{w}$ reminds that this map is related by the construction
 to a reduced decomposition \rf{q19}.

 Denote by $\Delta_w(\m)$ the set of
roots $\gamma_1,\ldots,\gamma_n$, defined in \rf{q19b}. Alternatively,
$\Delta_w(\m)$ consists of all roots $\gamma\in\Delta(\m)=w'(\Delta_+)$,
such that $w^{-1}(\gamma)\in\Delta(\m_-)=w'(\Delta_-)$.

\begin{lemma} {\rm (See \cite{Zh}, Section 5.2.4)}
%\noindent
\label{propq1}
Let $x\in\A'$. Then, in the notation \rf{q19a}-\rf{q10}, we have
\begin{align*}
(i)&&\q_{\overline{w},\m}(xd)&=\q_{\overline{w},\m}(x)\cdot
\shift_{w'(\rho)-w'w(\rho)}(d)&
\text{for any}\quad d&\in\Dh,\\
(ii)&&e_\a\cdot\q_{\overline{w},\m}(x)&=0 &\text{for any}\quad\a&\in\Delta_w(\m),\\
(iii)&&\q_{\overline{w},\m}(e_\a x)&=0 &\text{for any}\quad\a&\in\Delta_w(\m),\\
(iv)&&e_{w'w(\a)}\q_{\overline{w},\m}(x)&=\q_{\overline{w},\m}(e_{w'w(\a)}x)
&\text{if}\quad l(w'ws_\a)&>l(w'w)\, ,\\
(v)&&e_{w'(\a)}\q_{\overline{w},\m^{s_\a}}(x)&=\q_{\overline{w},\m^{s_\a}}
(e_{w'(\a)}x)
&\text{if}\quad l( w's_\a w)&>l(w'w)\, .
\end{align*}
\end{lemma}

{\it Proof}. The property (i) is a direct consequence of \rf{q8}, (ii).
Indeed, the relation \rf{q8}, (ii), for $k=0$  implies that
for any $x\in\A'$ and $d\in\Dh$ we have
 $$\q_{\overline{w},\m}(xd)=\q_{\overline{w},\m}(x)\shift_{\gamma_1+\ldots+
\gamma_n}(d)=
\q_{\overline{w},\m}(x)\cdot \shift_{w'(\rho)-w'w(\rho)}(d).$$
Suppose that the relations (ii)-(iv) take place for all $\m=w'(\n)$, and  all
reduced decompositions of any
element $w$ of the Weyl group of the length less than $n$, such that $l(w'w)=
l(w')+l(w)$. Let \rf{q19}
be a reduced decomposition of the element $w$ of the length $n$.

 In the notation of \rf{q19a}--\rf{q10} we have the equality
 $\,\q_{\overline{w},\m}=\q_{\overline{w}_{n-1},\m}\q_{\gamma_{n},\m_n}$,
with $\,l(w_{n-1})=n-1$. Thus, by the
induction assumption, for any $x\in\A'$,
$e_{\gamma_i}\q_{\overline{w}_{n-1},\m}(x)=0$ for $i=1,...,n-1$ and
 $e_{\gamma_{n}}\q_{\overline{w}_{n-1},\m}(x)=
\q_{\overline{w}_{n-1},\m}e_{\gamma_{n}}(x)$. This implies  equalities
$e_{\gamma_i}\q_{\overline{w}_{},\m}(x)=0$ for $i=1,...,n-1$ and
$$e_{\gamma_{n}}\q_{\overline{w}_{},\m}(x)=\q_{\overline{w}_{n-1},\m}(x)
e_{\gamma_{n}}
\q_{\gamma_{n},\m_n}(x)\, .$$
The last line is zero due to \rf{q8}, (iii). This proves the induction step
for the statement (ii).

On the other hand, present $w$ as a product $w=s_{\gamma_1}w''$, where
$w''$ is an element of the length
$n-1$ with a given reduced decomposition
 $\ow''=\{w''=s_{\a_{i_2}}\cdot\ldots\cdot s_{\a_{i_n}}\}$.
 We have a decomposition
$\q_{\overline{w},\m}=
\q_{\gamma_1,\m}\cdot{\q}_{\overline{w}'',\m_2}$, where
\begin{equation}\label{q20}
{\q}_{\overline{w}'',\m_2}=\q_{\gamma_2,\m_2}\q_{\gamma_{3},\m_3}
\cdot\ldots\cdot\q_{\gamma_{{n}},\m_n}\ .
\end{equation}
The induction assumptions say that ${\q}_{\overline{w}'',\m_2}(e_{\gamma_i}x)=0$
 for $i=2,...,n$ by (ii) and
$${\q}_{\overline{w}'',\m_2}(e_{\gamma_1}x)=
e_{\gamma_1}{\q}_{\overline{w}'',\m_2}(x)$$
 by (v). Thus ${\q}_{\overline{w},\m}(e_{\gamma_i}x)=0$
 for $i=2,...,n$ and
$${\q}_{\overline{w},\m}(e_{\gamma_1}x)=
\q_{\gamma_1,\m}{\q}_{\overline{w}'',\m_2}(e_{\gamma_1}x)=
\q_{\gamma_1,\m}(e_{\gamma_1}{\q}_{\overline{w}'',\m_2}(x))=0$$
by \rf{q4},(iv). This proves the induction step for (iii).

Let us prove the induction step for (iv).
Take a simple root $\a$, such that $l(w'ws_\a)>l(w'w)$. Let
$\gamma=w'w(\a)$. Then
$\gamma$ is a positive root and the sequence
\begin{equation}\label{q21}
\gamma_1,...,\gamma_n,\gamma
\end{equation}
is convex, that is, the sum of any two elements of the sequence lies
between them, if the sum is a root.
Let $\mu_0,\ldots,\mu_m$ be the finite "$\gamma_1$-sequence' of positive roots, starting
with $\gamma$, that is,
$\mu_0=\gamma$, $\mu_{k+1}=\mu_k+\gamma_1$. Then each $\mu_k$ belongs to the
set \rf{q21}, $\mu_k=\gamma_{i_k}$, where $i_k\in\{2,...,n\}$ if $k>0$,
 and $\e_\a^k(e_\gamma)=a_ke_{\mu_k}$ with $a_k\in\CC$, and
$\e_\a^k(e_\gamma)$=0 for $k>m$.
This implies, see the proof of Proposition \rf{prophom1}, that
\begin{equation}\label{q25}
e_{\gamma}\q_{\gamma_1}(y)=\q_{\gamma_1}(e_{\gamma}y)+
\sum_{k=1}^mb_k \q_{\gamma_1}^{(k)}(e_{\mu_k}y),\qquad b_k\in\CC\ .
\end{equation}
The relation \rf{q25} implies the equality
\begin{equation*}
%\label{q25a}
e_\gamma\q_{\overline{w},\m}(x)=\q_{\gamma_1,\m}
(e_{\gamma}{\q}_{\overline{w}'',\m_2}(x))+
\sum_{k=1}^mb_k \q_{\gamma_1,\m}^{(k)}(e_{\mu_k}{\q}_{\overline{w}'',\m_2}
(x))\ ,
\end{equation*}
where ${\q}_{\overline{w}'',\m_2}$ is defined in \rf{q20}.
The induction assumption says that
$e_{\gamma}{\q}_{\overline{w}'',\m_2}(x)=
{\q}_{\overline{w}'',\m_2}(e_{\gamma} x)$ and
$e_{\mu_k}{\q}_{\overline{w}'',\m_2}=0$ for any $k\geq 1$. This
implies the induction step for (iv). The statement (v) is proved in an
analogous manner.

\subsection{Map $\q_{w_0}$.}

In this Section we assume that $\g$ is a finite-dimensional reductive
Lie algebra and $\A$ is a $\g$-admissible algebra. In this case the Weyl
group $W$ is finite and the adjoint action of $\g$ in $\A$ is locally
 finite.

Set $\m=\n$ and for any $w\in W$ and a reduced decomposition
$\ow=\{w=s_{\a_{i_1}}\cdots s_{\a_{i_k}}\}$ denote the map $\q_{\ow,\n}$ as
$\q_{\ow}$:
$$\q_{\ow}\equiv\q_{\ow,\n}\,.$$
Let  $\ow_0=\{w_0=s_{\a_{i_1}}\cdots s_{\a_{i_n}}\}$ be a reduced decomposition
of the longest element $w_0\in W$. The map
$\q_{\ow_0}$ has the following properties:
 the right ideal $\n\A'$ is in the kernel of
$\q_{\ow_0}$ by Lemma \ref{propq1}, (ii), and
the image of $\q_{\ow_0}$ is in the normalizer of the left ideal
 $\A'\n$   by Lemma \ref{propq1}, (iii). It means that it induces
the map of the double coset space $\Cnn=\n\A'\backslash \A'/\A'\n_-$
to the Mickelsson
algebra $Z^\n(\A)$:
\begin{equation*}%\label{q28}
\q_{\ow_0}: \Cnn%\n\A'\backslash \A'/\A'\n_-
\to Z^\n(\A)\,.
\end{equation*}

 Let $\W\subset\A$ be an $\ad$-invariant generating  subspace of $\A$.
Let $z'_w$, where $w\in \W$, be the canonical generators of the $\Dh$-module
 $Z^\n(\A)$, see  Theorem \ref{proposition3a}.

\begin{proposition}
\label{propqw}
The mapping $\q_{\ow_0}$ of vector spaces defines an isomorphism of the double coset space
 $ \Cnn$ and  $ Z^{\n}(\A)\,,$
such that for any $x\in  \Ann$,
$d\in\Dh$ and $\w\in\W$ we have
\begin{align}\label{q30w}
[d,\q_{\ow_0}(x)]&=\q_{\ow_0}([d,x]), \qquad
\q_{\ow_0}(xd)=\q_{\ow_0}(x)\shift_{\rho-w_0(\rho)}(d)\,,
\\ \label{q27w}
\q_{\ow_0}(\w)&=z'_\w\,.
\end{align}
\end{proposition}

{\it Proof}. The arguments are the same as in the proof of Proposition
\ref{propqa}.

\noindent
The double coset space $\Cnn$ is a free
$\Dh$-module, generated by the vector space $\W$ . On the other hand,
the Mickelsson algebra $ Z^{\n}(\A)\,$ is also a free $\Dh$-module,
generated by the vectors $z'_\w$, $\w\in\W_\a$, see Theorem
\ref{proposition3a}. By the structure of the map $\q_{\ow_0}$,
$$
\q_{\ow_0}(\w)=w+\sum \w_i g_i ,
$$
where $g_i\in U'(\g)$, and\
$ \w_i\in\W$ have the weight strictly bigger then the weight of $\w$,
that is, $\mu(\w-\w_i)\in Q_+$, $\mu(\w-\w_i)\not=0$, where $\mu(x)\in\H^*$
denotes the weight of $x$. We can present further
any $g_i$ as a sum $g_i=\sum_j f_{i,j}d_{i,j}e_{i,j}$, where $f_{i,j}\in
U(\n_-)$, $d_{i,j}\in\Dh$ and $e_{i,j}\in U(\n)$. The terms, where
$e_{i,j}\not= 1$, vanish by definition in $ Z^{\n}(\A)\,$, so we have
the equality \rf{q27w}.\hfill{$\square$}

Proposition \ref{propqw} implies that
 the map $\q_{\ow_0}$ does not depend on a reduced
decomposition of $w_0$. We denote it further by $\q_{w_0}$,
$$\q_{w_0}\equiv \q_{\ow_0}\,.$$

\begin{corollary}\label{cor1}
The restriction of the map
$\q_{w_0}$ to the normalizer $\N(\A'\n_-)$
 defines an isomorphism of vector spaces
$ Z^{\n_-}(\A)$ and  $ Z^{\n}(\A)\,,$
such that
\begin{align}\notag
[d,\q_{\ow_0}(x)]&=\q_{\ow_0}([d,x]), \qquad
\q_{w_0}(xd)
=\q_{w_0}(x)\shift_{\rho-w_0(\rho)}(d),&
d&\in \Dh\,,\\ \label{q27ww}
\q_{w_0}(z_{\n_-,\w})&=z'_{\n,\w},& \w&\in \W\,.
\end{align}
\end{corollary}
Here ${z}_{\n_-,\w}$ are  generators of the Mickelsson algebra
$Z^{\n_-}(\A)$ of the 'first type', see Proposition \ref{propr3},
 ${z}'_{\n_,\w}\equiv {z'}_\w$ are generators
of the Mickelsson algebra
$Z^{\n}(\A)$ of the 'second type', see Proposition \ref{theorem00}.

{\it Proof}. We have $\n_-=\n^{w_0}$ and all the statements
of the Corollary follow by induction from 
Proposition \ref{prophom1}. We should prove only the equality \rf{q27ww}.
By the definition \rf{r7}, the element $z_{\n_-,\w}\in Z^{\n_-}(\A)$ has a
form
$$z_{\n_-,\w}=\w+ \sum_{i=1}d_i e_i \w_i\ ,
e_{i}\in \n U(\n),\ d_{i}\in \Dh, \ \w_i\in \W\,,$$
that is, $z_{\n_-,\w}=\w \mod \n\A'$. By the properties of the map
$\q_{w_0}$, we have
$\q_{w_0}(z_{\n_-,\w})=\q_{w_0}(\w)=z'_\w\,.$
\hfill{$\square$}

\subsection{Cocycle properties} In this Section
$\g$ is an arbitrary contragredient Lie algebra of finite growth, $\A$ is
a $g$-admissible algebra, $\m$  a maximal nilpotent subalgebra of $\g$,
 conjugated to $\n$ by an element of the Weyl group $W$, and $w$ an
 element of $W$, satisfying the condition \rf{q27}.

\begin{proposition}\label{prop3.5}  The maps $\q_{\ow,\m}$
do not depend on a choice of the reduced decompositions of $w\in W$.
\end{proposition}

{\it Proof}.  First we take an element $w\in W$, equal to the
  longest element  of a reductive subalgebra $\g'\subset\g$ of rank two. Then the
algebra $\g'=\n'_-+\H+\n'$ is generated by the elements $e_{\pm\gamma_i}$
and $h\in\H$,
where $\gamma_1,\ldots,\gamma_n$ are the members of the sequence \rf{q19b},
such that all $e_{\gamma_i}\in\m$ by the condition \rf{q27}.
By Proposition \ref{propqw}, the map
$$\q_{\ow,\n'}:\A'\to\A'/\A'\n'$$
does not depend on a reduced decomposition of $w$. The map
$$\q_{\ow,\m}:\A'\to\A'/\A'\m$$ is the composition of $\q_{\ow,\n'}$ and the
natural projection of $\A'/\A'\n'$ to $\A'/\A'\m$ and thus also does not
depend on a choice of the reduced decomposition of $w$.

This result implies the statement of the proposition for general $w$, since any
two reduced decompositions are related by a sequence of flips of reduced
 decompositions of longest elements of rank two subalgebras.
\hfill{$\square$}
\bigskip

 With the use of Proposition \ref{prop3.5}, we simplify further the
notation for the maps $\q_{\ow,\m}$ and $\q_{\ow}$ and write them as
$\q_{w,\m}$ and $\q_w$:
$$\q_{w,\m}\equiv\q_{\ow,\m},\qquad \q_w\equiv\q_{\ow}\equiv
 \q_{\ow,\n}\, .$$
Proposition \ref{prop3.5} can be formulated as the condition
\begin{equation}\label{coc}
\q_{w'w,\m}=\q_{w',\m}\q_{w,\m^{w'}}, \qquad\text{if}\quad
l(w'w)=l(w')+l(w)\,.
\end{equation}
This statement is known
 as 'Zhelobenko cocycle condition'.

\setcounter{equation}{0}
\section{Homomorphism properties of Zhelobenko maps}
\label{section5}

\subsection{Homomorphism property of maps $\q_\a$}

Let $\A$ be an admissible algebra over a contragredient Lie algebra $\g$
of finite growth. Let $\a$ be a real root of $\g$, $\g_\a$ the
$\mathfrak{sl}_2$-subalgebra of $\g$, generated by $e_{\pm\a}$ and $h_\a$,
$\n_{\pm\a}=\CC e_{\pm\a}$.
Since the algebra $\A$ is $\g$-admissible and the adjoint action of $\g_\a$ in $\g$
is locally finite, $\A$ is $\g_\a$-admissible as well, see Section
\ref{section4.1}. In this setting we
proved in  Section \ref{sectiondcs} that the double coset space
${}_{\n_\a}{\A}_{\n_{-\a}}=\n_{\a}\A'\backslash \A'/ \A'\n_{-\a}$
can be equipped with a structure of an associative algebra with the
multiplication rule
\begin{equation*}
a\circ b= a\PP_{{-\a}}(b)=\PPP_{{-\a}}(a)b,\,
\end{equation*}
in the notation of  Section \ref{section1.2} and Section
\ref{sectiondcs}.

The following theorem is  the basic point for applications of the Zhelobenko
operators in the representation theory of Mickelsson algebras.

\begin{theorem}
\label{theorem1}
The map
$\q_\a:\, {}_{\n_\a}{\A}_{\n_{-\a}} \to Z^{\n_\a}(\A)$
is a homomorphism of algebras.
\end{theorem}

 Theorem \ref{theorem1} and  Proposition \ref{propqa} imply that $\q_\a$
establishes an {\it isomorphism} of the double coset algebra
$ {}_{\n_\a}{\A}_{\n_{-\a}}$ and Mickelsson algebra $Z^{\n_\a}(\A)$.
On the other hand,  Theorem \ref{theorem1} implies the equality
\begin{equation}\label{hom18}
\q_{\a}(xy)=\q_{\a}(x)\q_{\a}(y)\,,\qquad
\text{for any}\quad x\in\A', \ y\in \N(\A'\n_{-\a})\,.
\end{equation}
Indeed, if $y\in \N(\A'\n_{-\a})$ and $\bar{y}$ is the class
 of $y$ in $\A'/\A'\n_{-\a}$ then $\bar{y}=\PP_{-\a}(\bar{ y})$ and the
first equality in \rf{hom8} is a corollary of the first statement of Theorem
\ref{theorem1}.
\medskip

{\it Proof}.
 Let $\P_{-\a}$ be the extremal projector, related to the decomposition
$\g_\a=\n_\a+\CC h_\a+\n_{-\a}$ of the
algebra $\g_\a$.
It is given  by the relations \rf{P6} and \rf{P6a}. The corresponding operators
$\PP_{-\a}:\A'/\n_{-\a} \A'\to \A'/\n_{-\a} \A'$ and
 $\PPP_{-\a}:  \A'\n_\a\backslash\A'\to  \A'\n_\a\backslash\A' $
can be written as
\begin{equation}
\label{hom0a}
\PP_{-\a}(x)=\sum_{m\geq0}\frac{1}{m!}\frac{1}{(h_\a-2)\cdots(h_\a-m-1)}
{e}_{\a}\e_{-\a}^m(x) \quad \mod \A'e_{-\a}\,,
\end{equation}
\begin{equation}\label{hom0b}
\begin{split}
\PPP_{-\a}(x)=&\sum_{m\geq0}\frac{(-1)^m}{m!}
\e_{\a}^m(x)e_{-\a}^m
\frac{1}{(h_\a-2)\cdots(h_\a-m-1)}\ \quad \mod e_\a\A'\ .
\end{split}
\end{equation}

We should establish an equality
\begin{equation}\label{hom4}
\q_\a\left(\PPP_{-\a}(x)y\right)=\q_\a(x)\q_\a(y)
\quad\text{for any}\ x\in\n_\a\A'\backslash\A',\quad y\in\A'\,.
\end{equation}
Suppose  $y\in\A$ is a weight vector of the weight $\mu_y\in\H^*$, such that
$[h_\a,y]=\mu_yy$ for any $h\in\H$. We have
%\begin{equation}
%\begin{split}

%\begin{align}\notag
\begin{equation}\begin{array}{l}
\q_\a\left(\PPP_{-\a}(x)y\right)=
\q_\a\left(\sum_{m\geq0}{\frac{(-1)^m}{m!}}
\e_{\a}^m(x)e_{-\a}^m
{\frac{1}{(h_\a-2)\cdots(h_\a-m-1)}}y\right)\\[1em] %\notag
\ \ =\q_\a\left(\sum_{m\geq0}{\frac{(-1)^m}{m!}}
\e_{\a}^m(x)e_{-\a}^m y
{\frac{1}{(h_\a-2+\mu_y)\cdots(h_\a-m-1+\mu_y)}}\right)\\[1em] %\notag
\ \ =\q_\a\left(\sum_{m\geq0}{\frac{(-1)^m}{m!}}
\e_{\a}^m(x)e_{-\a}^m y\right)
{\frac{1}{(h_\a+\mu_y)\cdots(h_\a-m+1+\mu_y)}}\\[1em] %\notag
\ \ =\sum_{n,m\geq0}{\frac{(-1)^{n+m}}{n!m!}}\e_\a^n\left(
\e_{\a}^m(x)e_{-\a}^m y\right)e_{-\a}^n{\frac{1}{g_{n,\a}}}
{\frac{1}{(h_\a+\mu_y)\cdots(h_\a-m+1+\mu_y)}} %\notag
\\[1em]
%\end{align}
%\begin{align}\notag
\ \ =\sum_{n,m\geq0}{\frac{(-1)^{n+m}}{n!m!}}\ \\[1em] %\notag
\ \ \ \ \cdot\sum_{k=0}^n
\left(\begin{array}{c}n\\ k\end{array}\right)
\e_{\a}^k\left(\e_{\a}^m(x)e_{-\a}^m\right)
\e_{\a}^{n-k}(y)e_{-\a}^n{\frac{1}{g_{n,\a}}}
{\frac{1}{(h_\a+\mu_y)\cdots(h_\a-m+1+\mu_y)}}\\[1em] %\label{hom5}
\ \ =\sum_{k,m=0}^\infty{\frac{(-1)^{m}}{m!k!}}
\e_{\a}^k\left(\e_{\a}^m(x)e_{-\a}^m\right)
{\frac{1}{(h_\a+2k)\cdots(h_\a+2k-m+1)}}\q_\a^{(k)}(y)\,.
%\end{split}
\end{array}\label{hom5}\end{equation}
%\end{align}
The sum of the terms with $k=0$ in \rf{hom5} is $\q_\a(x)\cdot q_\a(y)$.
So it is sufficient to prove the following identity in $\A'/\A'e_\a$:
\begin{equation}\label{hom6}
\sum_{k=1}^\infty\sum_{m=0}^\infty\frac{(-1)^{m}}{m!k!}
\e_{\a}^k\left(\e_{\a}^m(x)e_{-\a}^m\right)
\frac{1}{(h_\a+2k)\cdots(h_\a+2k-m+1)}\q_\a^{(k)}(y)=0
\end{equation}
for any $x,y\in\A$.

We have for any $k\geq 1$:
\begin{align}\notag
&\e^k_\a\left(\e_\a^m(x)e_{-\a}^m\right)\\
=&\e^{k-1}_\a \notag
\left(\e_\a^{m+1}(x)e_{-\a}^m+m\e_\a^m(x)
e_{-\a}^{m-1}(h_\a-m+1)\right)\\ \label{hom7}
=&\e^{k-1}_\a
\left(\e_\a^{m+1}(x)e_{-\a}^m\right)+m\e^{k-1}_\a
\left(\e_\a^{m}(x)e_{-\a}^{m-1}\right)(h_\a-m+1)\\ \notag
&-2m(k-1)
\e^{k-2}_\a
\left(\e_\a^{m}(x)e_{-\a}^{m-1}\right)e_\a\, .
\end{align}
Denote the left hand side of \rf{hom6} by $S$. Substitute \rf{hom7} into $S$. We get
\begin{align*}
S=&\sum_{k=1}^\infty\sum_{m=0}^\infty\frac{(-1)^{m}}{m!k!}\ \\
&\cdot\left(\e_{\a}^{k-1}\left(\e_{\a}^{m+1}(x)e_{-\a}^m\right)
\frac{1}{(h_\a+2k)\cdots(h_\a+2k-m+1)}\q_\a^{(k)}(y)\right.\\
&+m\e_{\a}^{k-1}\left(\e_{\a}^{m}(x)e_{-\a}^{m-1}\right)
\frac{(h_\a-m+1)}{(h_\a+2k)\cdots(h_\a+2k-m+1)}\q_\a^{(k)}(y)\\
&\left.-2m(k-1)\e_{\a}^{k-2}\left(\e_{\a}^{m}(x)e_{-\a}^{m-1}\right)e_\a
\frac{1}{(h_\a+2k)\cdots(h_\a+2k-m+1)}\q_\a^{(k)}(y)\right)\ .
\end{align*}
Substitute $(h_\a-m+1)=(h_\a+2k-m+1)-2k$ into the second sum  and use the
relation \rf{q4}, (iv) in the third sum. We get
\begin{align*}
S=&\sum_{k=1}^\infty\sum_{m=0}^\infty\frac{(-1)^{m}}{m!k!}\ \\
&\cdot\left(\e_{\a}^{k-1}\left(\e_{\a}^{m+1}(x)e_{-\a}^m\right)
\frac{1}{(h_\a+2k)\cdots(h_\a+2k-m+1)}\q_\a^{(k)}(y)\right.\\
&+m\e_{\a}^{k-1}\left(\e_{\a}^{m}(x)e_{-\a}^{m-1}\right)
\frac{1}{(h_\a+2k)\cdots(h_\a+2k-m+2)}\q_\a^{(k)}(y)\\
&-2km\e_{\a}^{k-1}\left(\e_{\a}^{m}(x)e_{-\a}^{m-1}\right)
\frac{1}{(h_\a+2k)\cdots(h_\a+2k-m+1)}\q_\a^{(k)}(y)\\
&+\left.2mk(k-1)\e_{\a}^{k-2}\left(\e_{\a}^{m}(x)e_{-\a}^{m-1}\right)
\frac{1}{(h_\a\!+2k-2)\cdots(h_\a+2k-m-1)}\q_\a^{(k-1)}(y)\right)\, .
\end{align*}
Change indices of summation: $m$ to $m+1$ in the second and third
sums;  $m$ to $m+1$ and  $k$ to $k+1$ in the last sum. Then
\begin{align*}
S=&\sum_{k=1}^\infty\sum_{m=0}^\infty\frac{(-1)^{m}}{m!k!}\ \\
&\cdot\left(\e_{\a}^{k-1}\left(\e_{\a}^{m+1}(x)e_{-\a}^m\right)
\frac{1}{(h_\a+2k)\cdots(h_\a+2k-m+1)}\q_\a^{(k)}(y)\right.\\
&-\e_{\a}^{k-1}\left(\e_{\a}^{m+1}(x)e_{-\a}^{m}\right)
\frac{1}{(h_\a+2k)\cdots(h_\a+2k-m+1)}\q_\a^{(k)}(y)\\
&+2k\e_{\a}^{k-1}\left(\e_{\a}^{m}(x)e_{-\a}^{m-1}\right)
\frac{1}{(h_\a+2k)\cdots(h_\a+2k-m)}\q_\a^{(k)}(y)\\
&-\left.2k\e_{\a}^{k-2}\left(\e_{\a}^{m}(x)e_{-\a}^{m-1}\right)
\frac{1}{(h_\a\!+2k)\cdots(h_\a+2k-m)}\q_\a^{(k)}(y)\right)=0\,.
\end{align*}
The theorem is proved. \hfill{$\square$}

\subsection{Properties of maps $\q_{\a,\m}$ and $\q_{w,\m}$}

Let $\a$ be a real root of the Lie algebra $\g$, $e_\a$ the corresponding root vector
with respect to the decomposition \rf{1} and $\m=\n^w$, where $w\in W$,
a maximal nilpotent
subalgebra of $\g$, such that $e_\a$ is a simple positive root vector
 of $\m$.

\begin{theorem}\label{theorem2} For any $x\in\A',\ y\in \N(\A'\m^{s_\a})$
\begin{align*}%\label{hom10}
\q_{\a,\m}(xy)=\q_{\a,\m}(x)\q_{\a,\m}(y)\,.
\end{align*}
\end{theorem}

{\it Proof}.  Since $y\in\N(\A'\m^{s_\a})$, its image
$\bar{y}$ in $M_{\m^{s_\a}}(\A')=\A'/\A'\m^{s_\a}$ satisfies the relations
$e_\gamma \bar{y}=0$ for all $\gamma\in\Delta_+(\m^{s_\a})$. In particular,
we have the equality
\begin{equation}\label{hom12}
e_{-\a}\bar{y}=0\,\qquad \text{in}\quad \A'/\A'\m^{s_\a}\,.
\end{equation}
Let $\g_\a\subset \g$ be the $\mathfrak{sl}_2$-subalgebra, generated by
$e_{\pm\a}$ and $h_\a$. Since $\A$ is $\g$-admissible, it is
$\g_\a$-admissible as well and $M_{\m^{s_\a}}(\A')$ is locally nilpotent with
respect to  $e_{-\a}$. Thus the
 extremal projector $\P_{-\a}$, related to $\n_{-\a}=\CC e_{-\a}$, acts in
$M_{\m^{s_\a}}(\A')$. Denote, following the notation of Section
\ref{section1.2}, its image in $\End M_{\m^{s_\a}}(\A')$ by $\PP_{-\a}$.
By \rf{hom12} and the properties of the extremal projector, we have
\begin{equation}\label{hom10}
\bar{y}=\PP_{-\a}(\bar{y}) \qquad \text{in}\quad \A'/\A'\m^{s_\a}\,.
\end{equation}
The equality \rf{hom10} can be read as
$$y\in \N(\A'\n_{-\a})\mod \A'\m^{s_\a}\,,$$
that is $y=y'+z$, where  $y'\in \N(\A'\m_{-\a})$ and $z\in \A'\m^{s_\a}$.

Indeed, $\PP_{-\a}(\bar{y})\in \N(\A'\m_{-\a})\mod \A'\m_{-\a}$ and
$\A'\m_{-\a}\subset \A'\m^{s_\a}$. Due to  Theorem \ref{theorem1},
see \rf{hom18}, and the
properties of the maps $\q_{\a,\m}$, we have
$$\q_{\a,\m}(xy)=\q_{\a,\m}(xy')=\q_{\a,\m}(x)\q_{\a,\m}(y')=
\q_{\a,\m}(x)\q_{\a,\m}(y)\,.$$
 \hfill{$\square$}

\medskip
Combining  Theorem \ref{theorem1} with the statements of Proposition
\ref{prophom1} and Proposition \ref{prophom2}, we conclude that the
restriction of the map $\q_{\a,\m}$ to  the normalizer
$\N(\A'\m^{s_\a})$ defines an isomorphism of the algebras
$Z^{\m^{s_\a}}(\A)$ and $Z^{\m^{}}(\A)$.

Iterations of these conclusions give the following statement.

\begin{proposition}\label{theorem3}
For any $w',w\in W$, such that $l(w'w)=l(w')+l(w)$ and $\m=\n^{w'}$,
the restriction of the map $\q_{w,\m}$ to the normalizer $\N(\A'\m^w)$
defines an isomorphism of the algebras $Z^{\m^w}(\A)$
 and $Z^{\m}(\A)$
such that for any $x\in Z^{\m^w}(\A)$  and $d\in\Dh$ we have
$$[d,\q_{{w},\m}(x)]=\q_{{w},\m}([d,x]),\qquad
\q_{{w},\m}(xd)=\q_{{w},\m}(x)\cdot
\shift_{w'(\rho)-w'w(\rho)}(d)\,
.$$
\end{proposition}
\medskip

The case of a finite-dimensional reductive Lie algebra $\g$ and
the element $w_0$ of the maximal length is special.
In this case the following statements hold:

\begin{proposition}\label{propw0}${}$
\begin{itemize}
\item[(i)]
The map $\q_{w_0}$ defines an isomorphism of the double coset algebra
$\Cnn$ and the Mickelsson algebra $Z^\n(\A)$ such that for any
$d\in \Dh$, $\w\in \W$ and $x\in Z_{\P_{\n_-}}(\A)$ we have
$$\q_{w_0}(\w)=z'_\w,
\qquad[d,\q_{{w_0}}(x)]=\q_{{w_0}}([d,x]),\qquad
 \q_{w_0}(xd)=
\q_{w_0}(x)\shift_{\rho-w_0(\rho)}(d)\, ;
$$
\item[(ii)]
the restriction of the map $\q_{w_0}$ to the normalizer $\N(\A'\n_-)$
defines an isomorphism of the Mickelsson algebras
$Z^{\n_-}(\A)$ and  $Z^\n(\A)$ such that for any
$d\in \Dh$, $\w\in \W$ and $x\in Z^{\n_-}(\A) $ we have
$$\q_{w_0}(z_{\n_-,\w})=z'_\w,
\qquad[d,\q_{{w_0}}(x)]=\q_{{w_0}}([d,x]),
\qquad \q_{w_0}(xd)=
\q_{w_0}(x)\shift_{\rho-w_0(\rho)}(d)\, .
$$
\end{itemize}
\end{proposition}
\smallskip

\setcounter{equation}{0}
\section{Braid group action}\label{section6}

\subsection{Operators $\qq_i$ }\label{section5.1}

Suppose that the automorphisms $\T_w:\, U(\g)\to U(\g)$, $w\in W$ admit extensions
 to automorphisms $\T^\AAA_w:\A\to\A$ of a $\g$-admissible algebra $\A$.
Such an extension is uniquely determined by automorphisms $\T_i^\AAA:\A\to\A$,
defined for all simple positive roots $\a_i$, which extend the automorphisms
$\T_i:U(\g)\to U(\g)$, see  Section \ref{section1.1} and satisfy braid group relations,
related to $\g$:
\begin{align}\label{b2}
%\T_i^\AAA(g)&=\T_i(g), && g\in U(\g)\,,\\ \label{b2}
\underbrace{\T_i^\AAA\T_j^\AAA\T_i^\AAA\cdots }_{m_{i,j}}&=
\underbrace{\T_j^\AAA\T_i^\AAA\T_j^\AAA\cdots }_{m_{i,j}}\,,&& i\not=j\,,
\end{align}
where
$m_{i,j}=2$, if $a_{i,j}=0$; $\,m_{i,j}=3$, if $a_{i,j}a_{j,i}=1$;
$\,m_{i,j}=4$, if $a_{i,j}a_{j,i}=2$; $\,m_{i,j}=6$, if $a_{i,j}a_{j,i}=3$;
$\,$
and $m_{i,j}=\infty$, if $a_{i,j}a_{j,i}>3$.
\ Having \rf{b2}, for any $w\in W$ we define the automorphism
$\T_w^\AAA:\A\to\A$ by the relation $\T_w^\AAA=\T_{i_1}^\AAA\cdots \T_{i_k}^\AAA$, where
$w=s_{\a_{i_1}}\cdots s_{\a_{i_k}}$ is a reduced decomposition of $w$.
Then the elements $\T_w^\AAA$ do not depend on a choice of the reduced decomposition
and satisfy the relations
\begin{equation}\label{b5}
\T^\AAA_{w'w}=\T_{w'}^\AAA\cdot\T_w^\AAA,\qquad\text{if}\quad l(w'w)=l(w')+l(w)\,.
\end{equation}
The automorphisms $\T_i^\AAA$ (and $\T_w^\AAA$) admit unique extensions
to automorphisms of
the algebra $\A'$, satisfying  \rf{b2}. We denote them by the
same symbols.

{}For any maximal nilpotent subalgebra $\m=\n^{w}$, where $w\in W$, and
a root $\alpha$, such that the root vector $e_\alpha$ is a simple positive
root vector of $\m$, we have the following relations:
\begin{equation}\label{b4}
\q_{\a,\m}=\T_w^\AAA \q_{w^{-1}(\a),\n}\T_w^{-1}.
\end{equation}

{}For each $\a_i\in\Pi$ define operators $\qq_i:\A'/\A'\n\to \A'/\A'\n$
 by the
relations
\begin{equation}\label{b3}
\qq_i=\q_{s_{\a_i}}\cdot \T_i^\AAA.
\end{equation}
In \rf{b3}, we understand  $\q_{s_{\a_i}}\equiv\q_{{\a_i},\n}$
as maps $\q_{{\a_i},\n}:\A'/\A'\n^{s_{\a_i}}\to
\A'/\A'\n^{}$, given by the relations \rf{q1a}--\rf{q8}.
 With the same agreement for any $w\in W$ we define operators
 $\qq_w:\A'/\A'\n\to \A'/\A'\n$
 by the relations
\begin{equation}\label{b3a}
\qq_w=\q_{w}\cdot \T_w^\AAA\ .\qquad
\end{equation}
The relation  \rf{b4}, Proposition \ref{prop3.5} and its analog for the maps $\nq_{\ow,\m}$ imply

\begin{proposition}\label{prop5.1}
Operators $\qq_i$  satisfy the braid group relations
\begin{equation}\label{b6}
\underbrace{\qq_i\qq_j\cdots }_{m_{i,j}}=
\underbrace{\qq_j\qq_i\cdots }_{m_{i,j}},
\qquad
 i\not=j\,.
\end{equation}
\end{proposition}

In other words, for any reduced decomposition
$w=s_{\a_{i_1}}\cdots s_{\a_{i_m}}$ we have  equalities
$\qq_w=\qq_{i_1}\cdots \qq_{i_m}$,
so
\begin{equation}\label{b6a}
\qq_{w'w}=\qq_{w'}\qq_w
\qquad\text{if}\quad l(w'w)=l(w')+l(w)\,.\end{equation}
{\it Proof}. Let $w',w\in W$ and $l(w'w)=l(w')+l(w)$.
We have by \rf{coc} and \rf{b4}
\begin{equation*}
\begin{split}
\qq_{w'w}= &\q_{w'w,\n}\T_{w'w}^\AAA=\q_{w',\n}\q_{w,\n^{w'}}\T_{w'w}^\AAA\\
=&\q_{w',\n}\left(\T_{w'}^\AAA\q_{w,\n}\T_{w'}^{-1}\right)
\T_{w'w}^\AAA=
\q_{w',\n}\T_{w'}^\AAA\q_{w,\n}\T_{w}^\AAA=\qq_{w'}\qq_{w}\,.
\end{split}
\end{equation*}
Thus we have the relations \rf{b6a}, which are equivalent to the braid group relations.
\hfill{$\square$}

{}For any $w\in W$ denote by $w\circ\ $ the shifted action of $w$ in $\H^*$:
$$w\circ\mu=w(\mu+\rho)-\rho\,.$$
It induces  the shifted action by automorphisms of $W$ in $\Dh$,
 characterized by the relations
\begin{equation}
\label{circ}w\circ h_\a=h_{w^{}(\a)}+\langle h_\a, w(\rho)-\rho\rangle\,.
\end{equation}
 Theorem \ref{theorem2} and Proposition \ref{theorem3} imply
\begin{proposition}\label{theorem5}
 For any $x\in \A'/\A'\n$, $y\in Z^\n(\A)$
and $d\in \Dh$ we have
\begin{equation}
\begin{split}\label{b7}
\qq_{i}([d,x])&=[s_{\a_i}(d),\qq_{i}(x)],\qquad
\qq_i(xd)=\qq_i(x)\cdot\left( s_{\a_i}\circ d\right)\,,\\
\qq_i(xy)&=\qq_i(x)\cdot\qq_i(y)\ .
\end{split}\end{equation}
\end{proposition}

\begin{corollary}${}$\label{cor5.3}
\begin{itemize}
\item[(i)]
The restriction of operators $\qq_i$ to $Z^\n(\A)$ defines an
automorphism of $Z^\n(\A)$, satisfying \rf{b7}.
\item[(ii)] For any $x\in Z^\n(\A)$
we have
\begin{equation*}
\qq_i^2(x)=(h_{\a_i}+1)^{-1}\T_i^2(x)(h_{\a_i}+1).
\end{equation*}
\end{itemize}
\end{corollary}

{\it Proof}. The statement (i) of the Corollary is a direct consequence
 of Proposition \ref{theorem5}. The statement (ii) follows from 
Proposition \ref{propinverse} and Proposition \ref{propqn}. Namely,
for any $x\in Z^\n(\A)$ we have by Proposition \ref{propinverse}
$$ \qq_i^2(x)=\q_{\a_i,\n}T_i^\AAA\q_{\a_i,\n}T_i^\AAA(x)=\q_{\a_i,\n}
\q_{-\a_i,\n^{s_\a}}(T_i^2(x))=(h_{\a_i}+1)^{-1}
\T_i^2(x)(h_{\a_i}+1).$$
\hfill{$\square$}

Clearly, all the statements of  Proposition \ref{theorem5} remain valid for
all operators $\qq_w$, $w\in W$. The properties \rf{b7},
\rf{b7a} look as
\begin{equation}\label{b9}
\begin{split}
\qq_{w}([h,x])&=[w(h),\qq_{w}(x)],\qquad
\qq_w(xd)=\qq_w(x)\cdot \left( w\circ d\right),\\
\qq_w(xy)&=\qq_w(x)\cdot\qq_w(y)
\end{split}
\end{equation}
\medskip

%%%%%%%%%%%%%%%%%%%%%%%%%%%%%%%%%%%%%%%%%%%%
\subsection{Calculation of $\qq_i(z_{\w})$}\label{section5.2}

Denote by $I_{\n_-}$ the image of the right ideal $\n_-\A'$ in  $\A'/\A'\n$:
$$I_{\n_-}= \left(\n_-\A'+\A'\n\right)/\A'\n\,.$$

\begin{lemma}\label{lemma3} For any $\a_i\in\Pi$ we have an inclusion
\begin{equation}\label{b15}
\qq_i(I_{\n_-})\subset I_{\n_-}.
\end{equation}
\end{lemma}

{\it Proof}.
We have to prove that for any $\gamma\in\Delta_+$ and $x\in\A'$
\begin{equation}\label{b16}
\q_{\a_i}(\T_i^\AAA(e_{-\gamma}x))=e_{-\mu}y
\end{equation}
for some $\mu\in\Delta_+$ and  $y\in\A'/\A'\n$. If $\gamma\not=\a_i$ then
$$\q_{\a_i}(\T_i^\AAA(e_{-\gamma}x))=\q_{\a_i}
(e_{-\gamma'}\cdot\T_i^\AAA( x))\,,$$
where $\gamma'=s_{\a_i}(\gamma)\in\Delta_+\setminus\a_i$, and the statement
of the lemma follows from the invariance of the subalgebra $\n_-(\a_i)$ with
respect to the action of $\e_{\a_i}$. Here $\n_-(\a_i)$ is generated by
root vectors $e_{-\gamma}$, where $\gamma\in\Delta_+\setminus\a_i$.

If $\gamma=\a_i$ then the right hand side of \rf{b16} vanishes due to
\rf{q4}, (i).\hfill{$\square$}
\medskip

Suppose that an $\ad$-invariant generating subspace $\W$ of a $\g$-admissible
algebra $\A$ is invariant
with respect to the action of the automorphisms $\T_i^\AAA$, $\T_i^\AAA(\W)=\W$.
Suppose that $\A$ satisfies the highest weight (HW) condition.
 With these assumptions we calculate the elements $\qq_i(z_\w)$, where
$\w\in\W$ and $z_\w$ are the generators of $Z^n(\A)$, defined in
Section \ref{section2.6}.

Keep the notation of  Section \ref{section2.6} and Section \ref{section2.8}.

\begin{proposition}${}$\label{prop5.5}
Assume that $\A$ satisfies the HW condition. Then
\begin{equation*}%\label{b20}
\qq_i(z_{\w})= z_{\BB_{\a_i}^{(2)}[-\rho](1\ot \T_i^\AAA(\w))}\,.
\end{equation*}
\end{proposition}

The operators  $\BB_{\pm\a}^{(2)}[\lambda]$ are defined in  Section \ref{section2.8}.

{\it Proof}.  Assume that $\A$ satisfies the HW condition.
By definition of the elements
$z_\w$, we have $z_\w=\w \mod I_{\n_-}$. Lemma \ref{lemma3} then
implies that
$$\qq_i(z_\w)=\qq_i(\w) \mod I_{\n_-}\,.$$
We have
\begin{equation*}
\begin{split}
\qq_i(z_\w)&=\qq_i(\w) \mod I_{\n_-}\\
&=\sum_{n\geq 0}\frac{(-1)^n}{n!}
 \e_{\a}^n(\w)e_{-\a}^ng_{\a,n}  \mod I_{\n_-} \\
&=
\sum_{n\geq 0}\frac{1}{n!}\e_{-\a}^n \e_{\a}^n(\w)g_{\a,n}\ \mod
I_{\n_-}\\
&=\sum_{n\geq 0}\frac{(-1)^n}{n!}\frac{(-1)^n}{(\Hh_\a-h_\a)(\Hh_\a-h_\a+1)
\cdots (\Hh_\a-h_\a+n-1)} \e_{-\a}^n \e_{\a}^n(\w) \mod I_{\n_-}.
\end{split}
\end{equation*}
This is precisely  the statement of the proposition. \hfill{$\square$}

{\bf Remark}. Let $w=s_{\a_1}\cdots s_{\a_n}$ be a reduced
decomposition of an element $w\in W$. Let $\gamma_1,\ldots ,\gamma_n$
be the corresponding sequence of positive roots: $\gamma_1=\a_1$,
$\gamma_2=s_{\a_1}(\a_2),\ldots$.
Then the properties of the maps $\qq_i$, see  Proposition
\ref{theorem5}, imply the following relation
$$\qq_w(z_{\w})=z_{\BB_{\gamma_1}^{(2)}[-\rho]\cdots
\BB_{\gamma_n}^{(2)}[-\rho](1\ot \T_w^\AAA(\w))}
%\qquad
%\qqq_w(\tilde{z}_\w)=
%\tilde{z}_{\BB_{-\gamma_1}^{(1)}[-\rho]\cdots
%\BB_{-\gamma_n}^{(1)}[-\rho](\T_w^\AAA(\w)\ot 1)}
\,.
$$
\medskip

\subsection{Calculation of $\qq_i(z'_{\w})$}\label{section5.3}

In this Section we assume that $\g$ is an arbitrary contragredient
Lie algebra of finite growth and the generating subspace $\W$ of a $\g$-admissible
algebra $\A$ is invariant with respect to the action of the automorphisms $\T_i^\AAA$.
 With this assumption we calculate the elements $\qq_i(z'_\w)$, where
$\w\in\W$ and $z'_\w$ are the generators of $Z^n(\A)$, defined in 
Section \ref{section2.6}.

\begin{proposition} %Let $\g$ be a finite-dimensional reductive Lie algebra.
Assume that the element $z'_\w\in \Zp(\A)$
 is defined. Then the element $\qq_i(z'_\w)$
 is also defined and is given by the relation
\label{theorem9}
\begin{equation}\label{b10}
\begin{split}
\qq_i(z'_\w)&=
z'_{\B_{-\a_i}^{(2)}[\rho]\left(1\ot\T_i^\AAA(\w)\right)}
%=z'_{\T_i\ot\T_i^\AAA\left(\B_{\a_i}[\rho](1\ot\w)\right)}
\,.
\end{split}\end{equation}
\end{proposition}

The proof of Proposition \ref{theorem9} is based on the following Lemma.

Let $\a$ be a simple positive root, $\n_\pm(\a)$ subalgebras of $\n_\pm$,
generated by all root vectors except $e_{\pm\a}$.

\begin{lemma}\label{lemmab1}
Suppose that the element $z'_\w$ exists. Then it admits a presentation
\begin{equation}\label{b11}
z'_\w=\q_{\a}(\w)+
\sum_{n\geq0,j}e_{-\a}^n{\w}_{j,n} {d}_{j,n}{f}_{j,n} \,,
\end{equation}
where ${\w}_{j,n}\in\W$, ${d}_{j,n}\in\Dh$, ${f}_{j,n}\in
\n_-(\a)U(\n_-(\a))$.
\end{lemma}

{\it Proof} of Lemma. Using the PBW theorem, present $z'_\w$ in a form
$$z'_\w= x+y,$$
where
\begin{equation*}\begin{split}
x&=\w+\sum\nolimits_k \w_k d_k e_{-\a}^k, \qquad \w_k\in \W, d_k\in\Dh,\\
y&=
\sum_{n\geq0,j}\tilde{{\w}}_{j,n}e_{-\a}^n {d}_{j,n}{f}_{j,n},\qquad
\tilde{{\w}}_{j,n}\in\W,  {d}_{j,n}\in\Dh, {f}_{j,n}\in
\n_-(\a)U(\n_-(\a)).
\end{split}
\end{equation*}
By definition, $z'_w$ satisfies the equation $[e_\a,z'_\w]=0 \mod\A'\n$.
Since the algebra $\n_-(\a)$ is invariant with respect to the adjoint action
of $e_{\a}$, the commutator $[e_\a,y]$ is an element $z$ of the same kind,
so we have two equations
$$[e_\a,x]=0 \mod\A'\n \qquad\text{and}\qquad
[e_\a,y]=0 \mod\A'\n\,.$$
The first equation has a unique solution $x=\q_\a(\w)$. To finish the proof of
the Lemma, we move all factors $e_{-\a}^n$ in the presentation of $y$ to the
left, using commutation relations in $\A$. \hfill{$\square$}

{\it Proof} of  Proposition \ref{theorem9}.\
The application of the automorphism $\T_i^\AAA$ to the presentation
\rf{b11} gives
\begin{equation}\label{b12}
\T_i^\AAA\q_{\a_i,\n}(u)=\T_i^\AAA\q_{\a_i,\n}(\w)+
\sum_{n\geq0,j} e_{\a_i}^n{\bar{\w}}_{j,n}
\bar{d}_{j,n}\bar{f}_{j,n}\,,
\end{equation}
where $\bar{\w}_{j,n}=\T_i^\AAA(\tilde{\w}_{j,n})\in\W$,
$\,\bar{f}_{j,n}=\T_i^\AAA(\tilde{f}_{j,n})\in \n_-(\a_i)U(\n_-(\a_i))$ and
$\,\bar{d}_{j,n}=\T_i^\AAA(\tilde{d}_{j,n})\in\Dh$.
Now we apply the map $\q_{\a_i,\n}$ to both sides of \rf{b12}.
The images of the terms in the last sum with $n>0$ vanish, since
$\q_{\a_i,\n}(e_{\a_i}x)=0$ for any $x\in\A'$ by \rf{q4}, (iv).
The images of the terms in the last sum with $n=0$ do not contribute
to the 'leading term', since the algebra $\n_-(\a_i)$ is
 invariant with respect to the adjoint action
of $\e_{\a_i}$. We obtained the statement:
\begin{lemma}\label{lemma5.3}
The leading term of  $\qq_i(z'_\w))$ is equal to the leading term of
$\q_{\a_i,\n}\T_i^\AAA\q_{\a_i}(\w)$.
\end{lemma}
Let us compute the leading term of $\q_{\a_i,\n}\T_i^\AAA\q_{\a_i}(\w)$.
We have
\begin{align*}&\q_{\a_i,\n}\T_i^\AAA\q_{\a_i}(\w)=
\q_{\a_i,\n}\T_i^\AAA\left(\sum_{n=0}^\infty\frac{(-1)^n}{n!}\e_{\a_i}^n(\w)
e_{-\a_i}^ng_{n,\a_i}
\right)\\
=&\q_{\a_i,\n}\left(\sum_{n=0}^\infty\frac{1}{n!}
\e_{-\a_i}^n(\T_i^\AAA(\w))e_{\a_i}^n\left((h_{\a_i})(h_{\a_i}+1)\cdots
(h_{\a_i}+n-1)\right)^{-1}\right)\\
=&\q_{\a_i,\n}\left(\sum_{n=0}^\infty\frac{1}{n!}
\e_{-\a_i}^n(\T_i^\AAA(\w))e_{\a_i}^n\right)\left((h_{\a_i}+2)(h_{\a_i}+3)\cdots
(h_{\a_i}+n+1)\right)^{-1}\,,
\end{align*}
by the property \rf{q4}, (iii). Since $\q_{\a_i,\n}(e_{\a_i}x)=0$ for any
$x\in\A'$, we get further
$$\q_{\a_i,\n}\T_i^\AAA\q_{\a_i}(\w)=
\sum_{n=0}^\infty\frac{(-1)^n}{n!}\q_{\a_i,\n}\left(\e_{\a_i}^n
\e_{-\a_i}^n(\T_i^\AAA(\w))\right)\left((h_{\a_i}+2)\cdots
(h_{\a_i}+n+1)\right)^{-1}\,.$$
Since $\e_{\a_i}^n\e_{-\a_i}^n(\T_i^\AAA(\w))$ belongs to $\W$ , the leading
term of $\q_{\a_i,\n}\left(\e_{\a_i}^n\e_{-\a_i}^n(\T_i^\AAA(\w))\right)$
is equal to $\e_{\a_i}^n\e_{-\a_i}^n(\T_i^\AAA(\w))$ and the leading term
of $\q_{\a_i,\n}\T_i^\AAA\q_{\a_i,\n}(\w)$ is equal to the sum
$$
\sum_{n=0}^\infty\frac{(-1)^n}{n!}\e_{\a_i}^n
\e_{-\a_i}^n(\T_i^\AAA(\w))\left((h_{\a_i}+2)\cdots
(h_{\a_i}+n+1)\right)^{-1}\,,$$
which can be written as
$\B^{(2)}_{-\a_i}[\rho]\left(1\ot\T_i^\AAA(\w)\right)$.
\hfill{$\square$}

As well as in the previous section, for a reduced decomposition
 $w=s_{\a_1}\cdots s_{\a_n}$
of an element $w\in W$ and a corresponding  sequence of positive roots:
$\gamma_1=\a_1$,
$\gamma_2=s_{\a_1}(\a_2),\ldots$ we have, assuming the existence of
$z'_{\w}$,
$$\qq_w(z'_{\w})=z'_{\B_{-\gamma_1}^{(2)}[\rho]\cdots
\B_{\gamma_n}^{(2)}[\rho](1\ot \T_w^\AAA(\w))}\,.
$$

\setcounter{equation}{0}
\section{Mickelsson algebra $\Zm(\A)$}\label{sectionzm}

In this section we collect the results of the previous section for
a Mickelsson algebra  $\Zm(\A)$, see the definition below.
This algebra deserves a special attention, since it acts on $\n_-$-coinvariants, which are
sometimes more convenient than $\n$-invariants.
If a $\g$-admissible algebra $\A$ admits an antiinvolution, which extends
the Cartan antiinvolution ${}^t: U(\g)\to U(\g)$, then all the results of this
section can be obtained by an application of this antiinvolution to
corresponding results of the previous sections.

\subsection{Algebra $\Zm(\A)$ and related structures}\label{sectionzn1}

{}For any associative algebra  $\A$, which contatins $U(\g)$, we
define Mickelsson algebras $S_{\n_-}(\A)$ and $\Zm(\A)$ as the quotients
$$S_{\n_{-}}(\A)=\n_-\A\backslash \N(\n_-\A),
\qquad
Z_{\n_{-}}(\A)=\n_-\A'\backslash \N(\n_-\A')\ ,
$$
where $\N(\n_-\A)$ ($\N(\n_-\A')$) is the normalizer of the right ideal
$\n_-\A$ ($\n_-\A'$). The algebra $\Zm(\A)$ is a localization of the
algebra $\Sm(\A)$: $\Zm(\A)=\Dh\ot_{U(\H)}\Sm(\A)$.

Alternatively, the Mickelsson algebra $Z_{\n_-}(\A)$ can be described
as a subspace of $\n_-$-invariants in a right $U(\g)$-module
$\tilde{M}_{\n_-}(\A')=\n_-\A'\backslash\A'$:
$$\Zm(\A)=\left(\tilde{M}_{\n_-}(\A')\right)^{\n_-}=
\{m\in \tilde{M}_{\n}(\A')\,|\,m\n_-=0\}\, .$$

As well as  $S^{\n{}}(\A)$, the space $S_{\n_{-}}(\A)$ is an associative
algebra, containing $U(\H)$, and for any left $\A$-module $M$, the space
$M_{\n_-}=M/\n_-M$ of $\n_-$-coinvariants is a  $S_{\n_-}(\A)$-module,
see Proposition \ref{prop5}.  The algebra $\Zm$ acts in the space $M_{\n_-}$ of
 $\n_-$-coinvariants of any left $\A'$-module $M$.

Suppose that a $g$-admissible algebra $\A$  satisfies the additional
 {\it local lowest weight} condition
\begin{itemize}
\item[(LW)]\ For any $\w\in\W$  the adjoint action of elements
$x\in U(\n_-)_\mu$ on $\w$ is nontrivial, $\Hat{x}(\w)\not=0$,
 only for a finite number of
$\mu\in\H^*$.
\end{itemize}

Then the quotient $\tilde{M}_{\n_-}(\A')=\n_-\A'\backslash\A'n$
has a structure of a right  $F_{\g,\n}$-module, extending the
action of $\A'$ by the right multiplication.
In particular, the extremal
projector $\P$ acts in the right $F_{\g,\n}$-module
$\tilde{M}_{\n_-}(\A')$. Denote the corresponding operator by
$\PPP\in\End \tilde{M}_{\n_-}(\A')$, see  Section \ref{section1.2}.

The properties of the extremal projectors imply the relation
\begin{equation}\label{?6}
\Zm(\A)= \text{Im}\PPP\subset \tilde{M}_{\n_-}(\A')\, .
\end{equation}

Equip the double coset space $\Ann$, see \rf{Znn},
with a multiplication $\circ:\Znn\ot\Znn\to\Znn$:
\begin{equation}
\label{?7}
a\circ b=a\P b\ \stackrel{def}{=}\PPP(a)\ \cdot b\ .
\end{equation}

We also call the double coset space $\Ann$, equipped with the operation \rf{?7},
the {\it double coset algebra} $\Znn$. In a case, when both conditions
(HW) and (LW) are satisfied, the multiplication rules \rf{?3} and
\rf{?7} coincide.

Define  linear maps $\phi_-:\Zm(\A)\to\Znn$ and $\psi_-:\Znn\to\Zm(\A)$
by the rules
\begin{equation}\label{?8}
\phi_-(x)=x\mod \A'\n,\qquad \psi^+(y)=\PPP(y),\qquad x\in\Zm(\A),y\in\Znn\ .
\end{equation}

\begin{proposition}\label{propgener2}
Assume that a $\g$-admissible algebra $\A$ satisfies the condition (LW). Then
\begin{itemize}
\item[(i)]
The operation \rf{?3} equips $\Znn$ with a structure of an associative
algebra.
\item[(ii)]
The linear maps $\phi_-$ and $\psi_-$ are inverse to each other
and establish an isomorphism of the algebras $\Zm(\A)$ and $\Znn$.
\end{itemize}
\end{proposition}

We have the 'lowest weight counterparts' of Propositions
\ref{propr3} and \ref{theorem00}.

\begin{proposition}
\label{propr3a}
Let $\A$ be a $\g$-admissible algebra satisfying the condition (LW). Then
\begin{itemize}
\item[(i)] Each element of the double coset algebra $\Znn$ can be uniquely
presented in a form $x=\sum_id_i\w_i$, where $d_i\in\Dh$, $w_i\in\W$, so
that $\Znn$ is a free left (and right) $\Dh$-module, isomorphic to
$\Dh\ot\W$.
\item[(ii)]
{}For each $\w\in \W$ there exists a unique element
$\tilde{z}_\w\in \Zm(\A)$
 of the form
\begin{equation}
\label{r7a}\tilde{z}_\w= \w+\sum_{i=1,...,k}\w_i  e_i d_i\ ,\qquad
e_{i}\in \n_{}U(\n_{}),\ d_{i}\in \Dh, \ \w_i\in \W \ ,
\end{equation}
so that the algebra $\Zm(\A)$ is a free left (and right) $\Dh$-module,
generated by the elements $\tilde{z}_\w$. The element $\tilde{z}_\w$ is equal to
$\PPP(\w)$.
\item[(iii)]
{}For each $\w\in \W$ there exists
 at most one element $\tilde{z}'_\w\in\Zm(\A)$ of the form
\begin{align}
\label{r8a}
\tzprim_\w&= \w+\sum\nolimits_{j} e_j  \w_jd_j\ ,
&&e_{j}\in \n_{}U(\n_{}),\ d_{j}\in \Dh, \ \w_j\in \W\,.
\end{align}
\end{itemize}
\end{proposition}

Next, we have analogs of Theorems \ref{proposition3a} and \ref{theorem0}:

\begin{theorem} \label{theorem0a} Let $\g$ be reductive and
finite-dimensional,
$\A$  a $\g$-admissible algebra with a generating subspace $\W$.
Then for any $\w\in \W$
\begin{itemize}
\item[(i)] there exists a unique element $\tilde{z}'_\w\in \Zm(\A)$
of the form \rf{r8a}. The algebra $\Zm(\A)$ is generated by the elements
$\tilde{z}'_\w$ as a free left (and right) $\Dh$-module;
\item[(ii)]
 we have the following equality in $Z_{\n_-}(\A)$
\begin{equation}\label{r14a}
\tilde{z}_\w=\tilde{z}'_{\B_-^{(2)}[\rho](1\ot\w)}\,.
\end{equation}
\end{itemize}
\end{theorem}
Here
$
\tilde{z}_{d\ot\w}=d\cdot \tilde{z}_\w,\quad \tilde{z}'_{d\ot\w}=d\cdot
\tilde{z}'_\w,\quad
\tilde{z}_{\w\ot d}=\tilde{z}_\w\cdot d ,\quad \tilde{z}'_{\w\ot d}=
\tilde{z}'_\w\cdot d$.

\medskip
\subsection{Zhelobenko maps }%Maps, related to $Z_{\n_-}(\A)$}

For any real root $\a$, the relations
\begin{align*}
%\notag
\nq_{\a}^{}(x)&=\sum\limits_{n=0}^\infty
\frac{1}{(n)!}{g_{n,\a}}\cdot e_{\a}^n\cdot\e_{-\a}^{n}(x) \
 &&\!\! \mod\ \n_{-\a}\A',&x&\in\A\ ,\\ %\label{n2}
\nq_{\a}^{}(dx)&=\shift_{\a}(d)\nq_{\a}^{}(x),&&&& d\in\Dh
\end{align*}
 define a map
$\nq_{\a}^{}:\A'\to \n_{-\a}\A'\backslash\A'$, such that for any
$x\in\A'$
\begin{equation}\label{n3}
\nq_{\a}^{}(e_{\a} x)=0,\qquad
 \nq_\a^{}(xe_{-\a})= \nq_\a^{}(x)e_{-\a}= 0\ .%\notag
\end{equation}
Here   $g_{n,\a}$ is given by \rf{q2},  $\n_{\pm\a}=\CC e_{\pm\a}$.
 We have analogs of Propositions \ref{propqa},
\ref{propinverse} and
Theorem \ref{theorem1}.

\begin{theorem}${}$
\label{propqn}
\begin{itemize}
\item[(i)]
The map $\nq_\a$ defines an isomorphism of algebras
$\nq_\a:\,{}_{\n_{\a}}\A_{\n_{-\a}}$ $\to$ $ Z_{\n_{-\a}}(\A)$,
such that for any $d\in\Dh$,
\begin{align*}%\label{q30}
[d,\nq_\a(x)]&=\nq_\a([d,x]),\qquad\nq_\a(dx)=\shift_\a(d)\nq_\a(x)\,.
%\\ \label{q27a}
%\q_\a(w)&=z'_w,& \w&\in \W\,.
\end{align*}
\item[(ii)] For any $ x\in Z_{\n_{\a}}(\A)$ and $y\in Z_{\n_{-\a}}(\A)$ we
have
\begin{equation}\label{inverse1}
%\begin{split}
\nq_{-\a}\nq_\a(x)=(h_\a+1)^{}x(h_\a+1)^{-1},\qquad
\nq_{\a}\nq_{-\a}(y)=(h_\a+1)^{-1}y(h_\a+1)^{}\,.
%\end{split}
\end{equation}
\end{itemize}
\end{theorem}

Theorem \ref{propqn} implies the equality
$$
\nq_\a(xy)=\nq(x)\nq(y) \qquad\text{for any}\qquad x\in\N(\n_\a\A'),\,
y\in\A'\,.
$$

{}For a maximal nilpotent
subalgebra $\m$ of $\g$, such that $e_\a$ is a simple positive root vector
of $\m$, define the linear map
$\nq_{\a,\m}^{(k)}: \A'\to\m_-\A'\backslash\A'$
by the prescriptions
\begin{align}\label{n4a}
\nq_{\a,\m}^{}(x)&=\sum\limits_{n=0}^\infty
\frac{1}{(n)!}{g_{n,\a}}\cdot e_{\a}^n\cdot\e_{-\a}^{n}(x) \
 & &\mod\ \m_{-}\A',& x&\in\A,\\ \label{n4b}
\nq_{\a,\m}^{}(dx)&=\shift_{\a}(d)\nq_{\a}^{}(x),&&& &d\in\Dh\ .
\end{align}
The assignment $\nq_{\a,\m}^{}$ satisfies the relation
$\ \nq_{\a,\m}^{}\left(\m_-^{s_\a}\A'\right)=0\,$
%(see Lemma \ref{lemma1})
and determines the map
$$
\nq_{\a,\m}^{}:\ \m_-^{s_\a}\A'\backslash\A'\to \m_-^{}\A'\backslash\A'\,
.$$
We have, see Proposition \ref{prophom1} and Theorem \ref{theorem2}:
\begin{proposition} ${}$
\label{prophom2}
\begin{itemize}
\item[(i)]
The map $\nq_{\a,\m}$ transforms  $\N(\m_-^{s_\a}\A')$ to
the Mickelsson
algebra $Z_{m_-}(\A)$ and % $=$ $\m_-\A'\backslash\N(\m_-\A')/$.
$$\nq_{\a,\m}(zu)=\nq_{\a,\m}(z)\nq_{\a,\m}(u)\qquad
 \text{for any} \quad z\in\N(\m^{s_\a}_-\A'),\ u\in \A'.$$

\item[(ii)] The restriction of the map $\nq_{\a,\m}$ to the normalizer
$\N(\m_-^{s_\a}\A')$ defines an isomorphism of the algebras
$Z_{m_-^{s_\a}}(\A)$ and $Z_{m_-}(\A)$, satisfying the relations
\begin{equation}\label{q8a}
[d,\nq_{\a,\m}^{}(x)]=\nq_{\a,\m}^{}([,x]),\qquad
\nq_{\a,\m}^{}(xd)=\nq_{\a,\m}^{}(x)\shift_\a(d),\qquad d\in\Dh\, .
\end{equation}
\end{itemize}
\end{proposition}

{}For any $w\in W$, satisfying the condition \rf{q27}, and its reduced
decomposition \rf{q19} we define a map
$\nq_{w,\m}:\ \m_-^{w}\A'\backslash\A'\to\m_-^{}\A'\backslash\A'\ $
by the relation
$$\nq_{w,\m}=\nq_{\gamma_1,\m_1}\cdot \nq_{\gamma_2,\m_2}\cdot\ldots\cdot
 \nq_{\gamma_n,\m_n}\, ,$$
where positive roots $\gamma_k$ and maximal nilpotent subalgebras $\m_k$
are defined by the prescriptions \rf{q19a}--\rf{q19b}. This map does not
depend on a choice of the reduced decomposition of $w$ and satisfies the
relations
\begin{equation}\label{coc1}
\nq_{w'w,\m}=\nq_{w',\m}\nq_{w,\m^{w'}} \qquad\text{if}\quad
l(w'w)=l(w')+l(w)\,,
\end{equation}
\begin{equation}\label{coc11}
[h,\nq_{{w},\m}(x)]=\nq_{{w},\m}([h,x]),\qquad
\nq_{{w},\m}(dx)=\shift_{w'(\rho)-w'w(\rho)}(d)\cdot\nq_{{w},\m}(x)
\,\end{equation}
for any\ $x\in\A'$, $h\in\H$ and $ d\in\Dh$.

The restriction of  $\nq_{w,\m}$ to the normalizer $\N(\m_-^w\A')$
defines an isomorphism of the algebras $Z_{\m^w_-}(\A)$ and $Z_{\m_-}(\A)$,
satisfying \rf{coc11}. We denote $\nq_w\equiv\nq_{w,\n}$.

The following counterpart of Proposition \ref{propqw} is valid.

\begin{proposition}
\label{prop3.6} Let $\g$ be a finite-dimensional reductive Lie algebra.
Let $w_0$ be the longest element of the Weyl group $W$. Then
\begin{itemize}
\item[(i)]
The map $\nq_{w_0}$ defines an isomorphism of the algebras
$ \Ann$ and  $ Z_{\n_-}(\A)$,
such that for any $x\in \Ann$,  $d\in\Dh$
and $\w\in\W$
\begin{align}
\label{q30n}
[d,\nq_{w_0}(x)]&=\nq_{w_0}([d,x]),\qquad
\nq_{w_0}(dx)=\shift_{\rho-w_0(\rho)}(d)\nq_{w_0}(x)\,,\\
\notag%\label{q27n}
\nq_{w_0}(\w)&=\tilde{z}'_\w\,.
\end{align}
\item[(ii)] The restriction of $\nq_{w_0}$ to the normalizer
$\N(\n\A')$ defines an isomorphism of the algebras $Z_\n(\A)$ and
$Z_{\n_-}(\A)$, satisfying \rf{q30n}, such that
\begin{equation}
\notag%\label{q27nn}
\nq_{w_0}(\tilde{z}_{\n,\w})=\tilde{z}'_{\n_-,\w}, \qquad\w\in \W\,.
\end{equation}
\end{itemize}
\end{proposition}

Here $\tilde{z}_{\n,\w}$ are the generators of the Mickelsson algebra
$Z_\n(\A)$ of the 'first type', see Proposition \ref{propr3a} (i),
 $\tilde{z}'_{\n_-,\w}$ are the generators of the Mickelsson algebra
$Z_{\n-}(\A)$ of the 'second type', see Proposition \ref{propr3a} (ii).

\subsection{Braid group action}

Keep the notation of  Section \ref{section6}.
We suppose again that the automorphisms $\T_i:\, U(\g)\to U(\g)$, $i=1,...,r$,
see  Section \ref{section1.1} admit extensions
 to automorphisms $\T^\AAA_i:\A\to\A$ of a $\g$-admissible algebra $\A$,
satisfying the braid group relations \rf{b2}.

{}Then for each maximal nilpotent subalgebra $\m=\n^{w}$, where $w\in W$, and
a root $\alpha$, such that the root vector $e_\alpha$ is a simple positive
root vector of $\m$, we have the relations:
\begin{equation}\label{b4n}
\nq_{\a,\m}=\T_w ^\AAA\nq_{w^{-1}(\a),\n}\T_w^{-1}.
\end{equation}
{}For each $\a_i\in\Pi$ define operators
$\qqq_i:\n_-\A'\backslash\A'\to \n_-\A'\backslash\A' $
and  $\qqq_w:\n_-\A'\backslash\A'\to \n_-\A'\backslash\A' $ as
$$\qqq_i=\nq_{s_{\a_i}}\cdot \T_i^\AAA,
\qquad
\qqq_w=\nq_{w}\cdot \T_w^\AAA\ .$$
The operators $\qqq_i$ satisfy the braid group relations,
\begin{equation*}
\underbrace{\qqq_i\qqq_j\cdots }_{m_{i,j}}=
\underbrace{\qqq_j\qqq_i\cdots }_{m_{i,j}},
\qquad
 i\not=j\,,
\end{equation*}
that is,
\begin{equation*}
\qqq_{w'w}=\qqq_{w'}\qqq_w,
\qquad\text{if}\quad l(w'w)=l(w')+l(w)\,.\end{equation*}

{}For any  $w\in W$, $x\in Z_{\n_-}(\A)$, $y\in \n_-\A'\backslash\A'$
and $d\in \Dh$ we have by Proposition \ref{prophom2}:
\begin{equation}
\begin{split}\label{b7a}
\qqq_{w}([h,x])&=[w(h),\qqq_{w}(x)],\qquad
\qqq_w(dx)=\left( w\circ d\right)\cdot\qqq_w(x)\,,\\
\qqq_w(xy)&=\qqq_w(x)\cdot\qqq_w(y)\ .
\end{split}\end{equation}

\begin{corollary}${}$\label{cor5.3n}
\begin{itemize}
\item[(i)] The restriction of $\qqq_i$ to $Z_{\n_-}(\A)$ defines an
automorphism of ${}$ $Z_{\n_-}(\A)$, satisfying \rf{b7a}.
\item[(ii)]  For any $y\in Z_{\n_-}(\A)$
we have
\begin{equation*}
\qqq_i^2(y)=(h_{\a_i}+1)^{-1}\T_i^2(y)(h_{\a_i}+1).
\end{equation*}
\end{itemize}
\end{corollary}

The following proposition describes the action of the automorphisms $\qqq_i$
on the canonical generators of the Mickelsson algebra $\Zm(\A)$.
\begin{proposition}${}$
\begin{itemize}
\item[(i)]  Assume that $\A$ satisfies the LW condition. Then
\begin{equation*}%\label{b20a}
\qqq_i(\tilde{z}_\w)= \tilde{z}_{\BB_{-\a_i}^{(1)}[-\rho](\T_i^\AAA(\w)\ot 1)}\,.
\end{equation*}
\item[(ii)]
Assume that the element
$\tilde{z}'_\w\in \Zm(\A)$ is defined. Then
$\qqq_i(\tilde{z}'_\w)$ is also defined and is given by the relation
\begin{equation*}%\label{b10}
\begin{split}
\qqq_i(\tilde{z}'_\w)&=
\tilde{z}'_{\B_{\a_i}^{(1)}[\rho]\left(\T_i^\AAA(\w)\ot 1)\right)}
\,.
\end{split}\end{equation*}
\end{itemize}
\end{proposition}
The operators $\BB^{(1)}_{\pm\a}[\lambda]$ and $\B^{(1)}_{\pm\a}[\lambda]$
are defined in Section \ref{section2.8}.

{\bf Remark}. Let $w=s_{\a_1}\cdots s_{\a_n}$ be a reduced
decomposition of an element $w\in W$. Let $\gamma_1,\ldots ,\gamma_n$
be the corresponding sequence of positive roots: $\gamma_1=\a_1$,
$\gamma_2=s_{\a_1}(\a_2),\ldots$.
Then the properties of the maps $\qq_i$ and $\qqq_i$, see Proposition
\ref{theorem5} imply the following relations
$$
\qqq_w(\tilde{z}_\w)=
\tilde{z}_{\BB_{-\gamma_1}^{(1)}[-\rho]\cdots
\BB_{-\gamma_n}^{(1)}[-\rho](\T_w^\AAA(\w)\ot 1)}\,,
\qquad
\qqq_w(\tilde{z}'_\w)=
\tilde{z}'_{\B_{\gamma_1}^{(1)}[\rho]\cdots
\B_{-\gamma_n}^{(1)}[\rho](\T_w^\AAA(\w)\ot 1)}\,.
$$

\setcounter{equation}{0}
\section{Standard modules and dynamical Weyl group}
\label{section7}

\subsection{Double coset space}
\label{section6.1}

Recall the notation $\ \Ann\ = \n_-\A'\backslash \A'/\A'\n\ $ from Section
\ref{sectiondcs}.
\begin{lemma}
The multiplication in $\A'$ equips the double coset space $\Ann$ with
the structure of the left $\Zm(\A)$-module and the right $\Zp(\A)$-module.
\end{lemma}

{\it Proof\,} follows from the definition of normalizers. \hfill{$\square$}

If $\A$ satisfies the HW condition, Proposition \ref{propgener1} says
that the double coset space $\Ann$ is a free right $\Zp(\A)$-module
of rank one, generated by the class of $1$.
If $\A$ satisfies the LW condition, Proposition \ref{propgener2} says
that the double coset space $\Ann$ is a free left $\Zm(\A)$-module
of rank one, generated by the class of $1$.

Lemma \ref{lemma3} says that the operators $\qq_i:\A'/\A'\n\to \A'/\A'\n$
and $\qqq_i:\n_-\A'\backslash \A'\to \n_-\A'\backslash \A'$
correctly define operators in the double coset space $\Ann$. We denote them
by the same symbol. According to definitions,
 they are given by the formulas
\begin{align*}
\qq_{i}(x)&=\sum\limits_{n=0}^\infty
\frac{(-1)^n}{n!}\e_{\a_i}^{n}(\T_i^\AAA(x))
\cdot e_{-\a_i}^n\cdot{g_{n,\a_i}}&
\mod\ &\n_-\A'+\A'\n\,,\\
\qqq_{i}(x)&=\sum\limits_{n=0}^\infty
\frac{1}{n!}{g_{n,\a_i}}\cdot e_{\a_i}^n\cdot\e_{-\a_i}^{n}(\T_i^\AAA(x))&
\mod\ &\n_-\A'+\A'\n\,,
\end{align*}
where $g_{n,\a_i}=\left(
h_{\a_i}(h_{\a_i}-1)\cdots(h_{\a_i}-n+1)\right)^{-1}$. Equivalently,
\begin{align}
\label{b17}
\qq_{i}(x)&=\sum\limits_{n=0}^\infty
\frac{1}{n!} \e_{-\a_i}^n\e_{\a_i}^{n}(\T_i^\AAA(x))
\cdot{g_{n,\a_i}}&
\mod\ &\n_-\A'+\A'\n\,,\\
\label{b18}
\qqq_{i}(x)&=\sum\limits_{n=0}^\infty
\frac{1}{n!}{g_{n,\a_i}}\cdot \e_{\a_i}^n\e_{-\a_i}^{n}(\T_i^\AAA(x))&
\mod\ &\n_-\A'+\A'\n\,.
\end{align}
They satisfy all the properties, mentioned in  Proposition
\ref{prop5.1}, Proposition \ref{theorem5} and Corollary \ref{cor5.3}.

In the notation of  Section \ref{section2.8}, the formulas \rf{b17} and
\rf{b18} mean that for any $\w\in\W$ we have the following equalities in
$\Znn$:
\begin{equation}\label{b19}
\qq_i(\w)= m\left(\BB^{(2)}_{\a_i}[-\rho](1\ot\T_i^\AAA(\w))\right),\qquad
\qqq_i(\w)= m\left(\BB^{(1)}_{-\a_i}[-\rho](\T_i^\AAA(\w)\ot 1)\right)\,,
\end{equation}
where $m:\Dh\ot\W\to\A'$ and $m:\W\ot\Dh\to\A'$ are the multiplication maps.

Let $M$ be a module over an associative algebra $U$. Denote by
$\xi_M$ the corresponding homomorphism $\xi_M:U\to \End(M)$.
Let $T:U\to U$ be an automorphism of the algebra $U$. Denote by $M^T$
the $U$-module $M$, conjugated by the automorphism $T$. It can be
 described as follows. $M^T$
 coincides with $M$ as a vector space, while the map
 $\xi_{M^T}:U\to \End(M)\equiv\End(M^T)$ is
$$\xi_{M^T}=\xi_M\cdot T\,.$$
In this notation, Proposition \ref{theorem5} states the equivariance of
the maps $\qq_i$ and $\qqq_i$:

\begin{proposition} ${}$
\begin{itemize}
\item[(i)] The map $\qq_i$, given by \rf{b17}, is a morphism
of the right $\Zp(\A)$-modules:
$$\qq_i: \Ann\to \left(\Ann\right)^{\qq_i}\, .$$
\item[(ii)] The map $\qqq_i$, given by \rf{b18}, is a morphism
of the left $\Zm(\A)$-modules:
$$\qqq_i: \Ann\to \left(\Ann\right)^{\qqq_i}\,.$$
\end{itemize}
\end{proposition}

The same  statement holds for the operators $\qq_w$, defined as products of \rf{b17},
and the operators $\qqq_w$, defined as products of \rf{b18}, where $w\in W$.

Let $\lambda\in\H^*$ be a generic weight, that is,
$\langle h_\a,\lambda\rangle\not\in\ZZ$ for any $\a\in\Delta$.
Then the following quotients of the double coset space are well
defined:
\begin{equation}\begin{split}\label{b20}
\Annl&=\n_-\A'\backslash \A'/\A'\cdot
\left(\n,h-\langle h,\lambda\rangle)|_{h\in\H}\right)\,,\\
\Alnn&=
\left((h-\langle h,\lambda\rangle)|_{h\in\H},\n_-\right)
\A'\backslash \A'/\A'\n
\,.
\end{split}
\end{equation}
The space $\Annl$ is a left $\Zm(\A)$-module,
the space $\Alnn$ is a right $\Zp(\A)$-module.

\begin{corollary} For a generic $\lambda\in\H^*$,
\begin{itemize}
\item[(i)] the map  \rf{b17} defines a morphism
of the right $\Zp(\A)$-modules:
$$\qq_{i,\lambda}: \Alnn\to
\left({}_{s_{\a_i}\circ\lambda,\n_-}{\A}_{\n}\right)^{\qq_i}\, ;$$
\item[(ii)] the map  \rf{b18} defines a morphism
of the left $\Zm(\A)$-modules:
$$\qqq_{i,\lambda}: \Annl\to
\left({}_{\n_-}{\A}_{\n,s_{\a_i}\circ\lambda}\right)^{\qqq_i}\,.$$
\end{itemize}
\end{corollary}

Due to \rf{b19}, we have
$$\qq_{i,\lambda}(v)=\Hat{\PP}_{\a_i}[\lambda-\rho](\T^\AAA_i(v)),
\qquad
\qqq_{i,\lambda}(v)=\Hat{\PP}_{-\a_i}[\lambda-\rho](\T^\AAA_i(v)),$$
that is,
\begin{equation*}
\begin{split}
\qq_{i,\lambda}(v)&=\sum_{n\geq 0}\frac{(-1)^n}{n!}
\prod_{k=1}^n
\left(\Hh_{\a_i}+\langle h_{\a_i},\lambda-\rho\rangle+k\right)^{-1}
\e_{-\a_i}^n\e_{\a_i}^n(T_i^\AAA(v)),\\
\qqq_{i,\lambda}(v)&=\sum_{n\geq 0}\frac{(-1)^n}{n!}
\prod_{k=1}^n
\left(-\Hh_{\a_i}+\langle h_{\a_i},\lambda-\rho\rangle+k\right)^{-1}
\e_{\a_i}^n\e_{-\a_i}^n(T_i^\AAA(v)).
\end{split}
\end{equation*}

More generally, for any element $w\in W$, the maps
 $\qq_w:\A'/\A'\n\to \A'/\A'\n$ and $\qqq_w:\n_-\A'\backslash\A'\to
\n_-\A'\backslash\A'$ define morphisms of modules over Mickelsson algebras
$$\qq_{w,\lambda}:\Alnn\to
\left({}_{w\circ\lambda,\n_-}{\A}_{\n}\right)^{\qq_w}\ \ \ \ {\mathrm{and}}\ \ \ \ 
\qqq_{w,\lambda}: \Annl\to
\left({}_{\n_-}{\A}_{w\circ\lambda,\n}\right)^{\qqq_w}\ .$$
 For a reduced decomposition
 $w=s_{\a_1}\cdots s_{\a_n}$
of an element $w\in W$ and a corresponding  sequence of positive roots
$\gamma_1,\ldots ,\gamma_n$,
 we have
$$\qq_{w,\lambda}(v)=\Hat{\PP}_{\gamma_1}[\lambda-\rho]\cdots
\Hat{\PP}_{\gamma_n}[\lambda-\rho](\T^\AAA_w(v)),
\quad
\qqq_{i,\lambda}(v)=\Hat{\PP}_{-\gamma_1}[\lambda-\rho]\cdots
\Hat{\PP}_{-\gamma_n}[\lambda-\rho](\T^\AAA_w(v)).$$

{\bf Remark}.
%We now pass to Mickelsson algebras $S^{\n}(\A)$ and $S_{\n_-}(\A)$.
By the definition, the double coset space
$\,{}_{\n_-}{\widetilde{\A}}_{\n,\lambda} $
 $=$ $
\n_-\A\backslash \A/\A\n$ is a left $S_{\n_-}(\A)$- and a right
$S^{\n}(\A)$-module. Its quotients
\begin{equation*}\begin{split}
\Bnnl&=\n_-\A\backslash \A/\A\cdot
\left(\n,(h-\langle h,\lambda\rangle)|_{h\in\H}\right)\,,\\
\Blnn&=
\left(\n_-,(h-\langle h,\lambda\rangle)|_{h\in\H}\right)
\A\backslash \A/\A\n
\,
\end{split}
\end{equation*}
coincide with the spaces \rf{b20} for generic $\lambda$: $\Bnnl=\Annl$ and
$\Blnn=\Alnn$. They have correspondingly the
structure of  a left $S_{\n_-}(\A)$- and a right
$S^{\n}(\A)$-module. The module
$\Bnnl$ can be interpreted as a space of
$\n_-$-coinvariants in the left $\A$-module $M_{\n,\lambda}(\A)=
\A/\A\cdot\left(\n,h-\langle h,\lambda\rangle|h\in\H\right)$.
The module
$\Blnn$ can be interpreted as a space of
$\n$-coinvariants in the right $\A$-module
 $\tilde{M}_{\nml}(\A)=
\left(\n_-,h-\langle h,\lambda\rangle|h\in\H\right)
\A\backslash \A$.

{}For each $i$, the operators $\qq_i$ and $\qqq_i$
define homomorphisms of Mickelsson algebras $S^{\n}(\A)$ and
$S_{\n_-}(\A)$ to their localizations with respect
 to denominators, generated by monomials
$(h_{\a_i}+k)$, $k\in \ZZ$. If $\lambda\in\H^*$ satisfies
the condition $\langle h_{\a_i},\lambda\rangle\not\in\ZZ$, these
localizations act in
$\left({}_{s_{\a_i}\circ\lambda,\n_-}{\widetilde{\A}}_{\n}\right)^{\qq_i}$ and
$\left({}_{\n_-}{\widetilde{\A}}_{\n,s_{\a_i}\circ\lambda}\right)^{\qqq_i}$
correspondingly. In this sense the operators $\qq_i$ and $\qqq_i$
define morphisms of the right
$S^{\n}(\A)$-modules and of the left $S_{\n_-}(\A)$-modules:
\begin{equation}\label{b21}
\qq_{i,\lambda}: \Blnn\to
\left({}_{s_{\a_i}\circ\lambda,\n_-}{\widetilde{\A}}_{\n}\right)^{\qq_i}
\quad\text{and}\quad
\qqq_{i,\lambda}: \Bnnl\to
\left({}_{\n_-}{\widetilde{\A}}_{\n,s_{\a_i}\circ\lambda}\right)^{\qqq_i}\,.
\end{equation}
One can regard  \rf{b21} as a family of operators
with a meromorphic dependence on a parameter $\lambda$, study their
singularities, residues etc.

\subsection{Quotients of free modules}

As any associative algebra with unit, the Mickelsson algebra
is a free left and a free right module over itself of rank one.
Let us restrict ourselves to the Mickelsson algebra $\Zp(A)$ and its
free right module of rank one.

Let $\lambda\in\H^*$ be a generic weight. Consider the following quotient of the
 free right $\Zp(A)$-module
\begin{equation}\label{b23a}
\Phi_\lambda(\A)=(h-\langle h,\lambda\rangle)|_{h\in\H}\Zp(\A)\backslash
\Zp(\A).
\end{equation}
It can be realized as follows.
The multiplication $m$ in $\A$ induces an isomorphism of the two left
$U(\g)$-modules:
\begin{equation}\label{b23}
\W\ot M_\n(\g) \rightsquigarrow     M_\n(\A')\,,
\end{equation}
where $M_\n(\g)$ is the 'universal Verma module' $U'(\g)/U'(\g)\n$,
$M_\n(\A')=\A'/A'\n$,
 and
$\W$ is taken with a structure of the adjoint representation of $\g$.
The map \rf{b23} is also an isomorphism of the right $\Dh$-modules, where
the structure of $\Dh$-modules in the left hand side of \rf{b23} is given by
the prescription $(\w\ot m)\cdot d=\w\ot (m\cdot d)$ for any $\w\in\W$,
$m\in M_\n(\g)$, $d\in\Dh$.
 The Mickelsson algebra
is the space of highest weight vectors in
$M_\n(\A')$, so with the identification \rf{b23} we have the following
isomorphism of $\Dh$-bimodules:
\begin{equation}\label{b24}
\Zp(\A)\approx \left(\W\ot M_\n(\g)\right)^\n\,.
\end{equation}
Recall that $\Zp(\A)$ is a $U(\H)$-bimodule and admits the weight decomposition
with respect to the adjoint action of $\H$.
This implies that the right $\Zp(\A)$-module $\Phi_\lambda(\A)$ is a
semisimple right $U(\H)$-module and admits a decomposition
 $$\Phi_\lambda(\A)= \oplus_{\nu}\Phi_{\lambda,\lambda-\nu},$$
where the sum is taken over the weights $\nu$ of $\A$ such that for any
$\varphi\in \Phi_{\lambda,\mu}(\A)$ we have
\begin{equation}\label{b25}
\varphi\cdot h=\langle h, \mu\rangle \varphi\
,\end{equation}
 and  $\Phi_{\lambda,\mu}$ coincides with the double coset
\begin{equation*}%\label{b25}
\Phi_{\lambda,\mu}=
(h-\langle h,\lambda\rangle)|_{h\in\H}\cdot\Zp(\A)\backslash
\Zp(\A)/\Zp(\A)\cdot(h-\langle h,\mu\rangle)|_{h\in\H}\ .
\end{equation*}
Analogously, the left
$\Zp(A)$-module $\Phi_{.\,,\mu}(\A)$
\begin{equation*}
\Phi_{.\,,\mu}(\A)=\Zp(\A)/\Zp(\A)\cdot (h-\langle h,\mu\rangle)|_{h\in\H}
\end{equation*}
admits a presentation
 $$\Phi_{.\,\mu}(\A)= \oplus_{\nu}\Phi_{\mu+\nu,\mu}.$$

With the identification \rf{b24}, the relation \rf{b25} can be interpreted as
follows: any element $\phi\in \Phi_{\lambda,\mu}(\A)$ is a highest weight
vector in the tensor product $\W\ot M_{\mu}$ of the weight $\lambda$.
 Here $M_\mu= U(\g)/U(\g)\cdot(\n,(h-\langle h,\mu\rangle)h\in\H)$
is the Verma module of $g$ with the highest weight $\mu$.
We have proved

\begin{lemma} Let $\lambda\in\H^*$ be generic and $\nu\in\H^*$ a weight
of $\Zp(\A)$. Then
the weight space $\Phi_{\lambda,\lambda-\nu}$ of the right
 $\Zp(\A)$-module
$\Phi_\lambda(\A)$ is isomorphic to the space of intertwinig operators
\begin{equation}\label{b26}
\Phi_{\lambda,\lambda-\nu}\approx {\rm Hom\,}_{U(\g)}\left(M_\lambda,
\W\ot M_{\lambda-\nu}\right)\,.
\end{equation}
\end{lemma}
\medskip

{}Denote by $\mathbb{I}_\lambda$ the class of unit $1\in \Zp(\A)$ in
$\Phi_\lambda(\A)$.
The vector $\mathbb{I}_\lambda$ generates $\Phi_\lambda(\A)$ as a $\Zp(\A)$-module.
{}For any $\w\in\W$ denote by $\Phi_\lambda^v$ the vector of the right
$\Zp(\A)$-module $\Phi_\lambda(\a)$, obtained by an application of the element
$z'_\w$ to  $\mathbb{I}_\lambda$. It is equal to the class of $z'_\w$ in
$\Phi_\lambda(\A)$. Let $\nu$ be the weight of $\w$. In the description \rf{b26},
$\Phi_\lambda^v$ presents a map  $\Phi_\lambda^v$ $\in$
${\rm Hom\,}_{U(\g)}\left(M_\lambda,\W\ot M_{\lambda-\nu}\right)$,
such that
$$\Phi_\lambda^v(1_\lambda)=\w\ot 1_{\lambda-\nu}+ \text{l.o.t.}$$
where $1_\lambda$ is the highest weight vector of Verma module $M_\lambda$
and $\text{l.o.t.}$ mean terms which have lower weight on the second tensor
component.

The properties of the operators $\qq_w$ imply that
$\qq_w(\Phi_{\lambda,\mu})=\Phi_{w\circ\lambda,w\circ\mu}$,
so that each operator $\qq_w$ defines morphisms of the right and left
$\Zp(\A)$-modules:
\begin{equation}\label{b28a}
\qq_{w,\lambda}:\Phi_\lambda(\A)\to
\left(\Phi_{w\circ\lambda}\right)^{\qq_w}(\A),\qquad
\qq_{w,\,.\,\mu}:\Phi_{.\,,\mu}(\A)\to
\left(\Phi_{.\,,w\circ\mu}\right)^{\qq_w}(\A)
\,.
\end{equation}
Proposition \ref{theorem9} gives a formula for transformations of
vectors $\Phi_\lambda^v$:
\begin{equation}\label{b28}
\qq_{w,\lambda}(\Phi_\lambda^v)=\Phi_{w\circ\lambda}^{\Hat{\PP}_{-\gamma_1}[\lambda+\rho]\cdots
\Hat{\PP}_{-\gamma_n}[\lambda+\rho](\T_w^\AAA(\w))}\,,
\end{equation}
where $\gamma_1,...,\gamma_n$ is the sequence of positive roots,
 attached to a reduced decomposition $w=s_{\a_{i_1}}\cdots s_{\a_{i_1}}$
by the standard rule $\gamma_1=\a_{i_1}$, $\gamma_2=s_{\a_{i_1}}(\a_{i_2})$,
$\ldots$ .
\smallskip

\subsection{ Relations to dynamical Weyl group}

Let $\W$ be a  $U(\g)$-module algebra with a locally nilpotent action
of real root vectors and  $\A_\W=U(\g)\ltimes \W$
be a smash product of $U(\g)$ and $\W$, see the example 2 in Section
\ref{section2.1}.

The typical examples are: the tensor algebra of an integrable highest weight representation
of $\g$ and the symmetric algebra of an integrable highest weight representation of $\g$.

Due to assumptions above, the $\g$-module structure in $\W$ lifts
to an action of the Weyl group of $\g$ in $\W$ by standard formulas
\begin{equation}\label{b32}
\hat{\T}_i=\exp \e_{\a_i}\cdot \exp -\e_{-\a_i}\cdot\exp \e_{\a_i}\ .
\end{equation}
Due to \rf{b32}, operators $\hat{T}_i$ are automorphisms of the algebra
$\W$.

The operators $\hat{T}_i$ admit a lift to automorphisms of $\A$ by the relation
$\T_i(g\w)=\T_i(g)\hat{T}_i(\w)$, where $g\in U(\g)$, $\w\in\W$, and
$\T_i(g)$ is the automorphism of $U(\g)$, as in Section
\ref{section1.1}. Actually, they are given by the same relation
\rf{b32} with respect to the adjoint action of $\g$ in $\A$.

Elements of  right $\Zp(\A)$-modules $\Phi_\lambda(\A)$
are known in this case under the name {\it intertwining operators}.
The morphisms $\qq_{w,\lambda}$ generate a so called {\it dynamical Weyl group}
action, see \cite{EV,TV}.

The intertwining operators  $\Phi_\lambda(\A)$ form
an algebroid with respect to the composition operation. The composition of
intertwining operators can be described as follows.

{}For a generic $\lambda\in\H^*$, any morphism
 $\varphi_\lambda: M_\lambda\to \W\ot M_{\lambda-\nu}$ of
$\g$-modules admits a lift to a morphism $\bar{\varphi}_\lambda:\W\ot M_{\lambda}\to
\W\ot M_{\lambda-\nu}$ of $\A'$-modules by the rule
$\bar{\varphi}_\lambda(\w\ot m)=\w\cdot \varphi_\lambda(m)$
for any $m\in M_\lambda$. Then the composition
${\varphi}'_{\lambda-\nu}\circ
 {\varphi}_\lambda$ of intertwinig operators
$\varphi_\lambda\in\Phi_{\lambda,\lambda-\nu}$ and
$\varphi'_{\lambda-\nu}\in\Phi_{\lambda-\nu,\lambda-\nu-\nu'}$ is an
element
$\varphi''_{\lambda}\in\Phi_{\lambda,\lambda-\nu-\nu'}$, such that
$$\bar{\varphi}''_{\lambda}\ =\ \bar{\varphi}'_{\lambda-\nu}\circ
 \bar{\varphi}_\lambda\,.$$
The composition of intertwining operators coincides with the structure
of the right $\Zp(\A)$-module in $\Phi_\lambda$. Namely, in the notation of
the previous section, for any $x\in \Zp(\A)$,
denote by $\Phi_\lambda^x$ its class in $\Phi_\lambda$, considered as intertwining operator.
 Then we have

\begin{proposition}\label{b33}
\cite{Kh} Let $z',z''\in\Zp(\A)$. Assume that the weight of $z'$ with respect to the
adjoint action of $\H$ is $\nu$. Then
$$\Phi_\lambda^{z'z''}=\Phi_{\lambda-\nu}^{z''}\circ\Phi^{ z'}_\lambda\ .$$
\end{proposition}

The statement that the maps $\qq_{w,\lambda}$ are morphisms of
$\Zp(\A)$-modules, see \rf{b28a}, is equivalent in this context to the
statement that the dynamical Weyl group action respects the composition of
intertwining operators.

%%%%%%%%%%%%% New Section %%%%%%%%%%%%%%
\setcounter{equation}{0}
\section{Quantum group settings}\label{section8}

\subsection{Notation and assumptions}\label{section7.1}

In this section we announce basic statements of the paper for Mickelsson
algebras, related to reductions over quantum groups. We restrict our
attention to the Mickelsson algebras $\Zp(\A)$.

Keep the notation of  Section \ref{section1.1}.
Let $\nu$ be an indeterminate; $d_i\in\NN$ are defined by the condition
that the matrix $(\a_i,\a_j)=d_ia_{i,j}=d_ja_{j,i}$ is symmetric.
Here $a_{i,j}$
is the Cartan matrix of $\g$.
For any root $\gamma\in\Delta$ put $\nu_\gamma=\nu^{(\gamma,\gamma)/2}$,
 $[a]_{p}=\frac{p^a-p^{-a}}{p-p^{-1}}$,
and $(a)_{p}=\frac{p^a-1}{p-1}$
for any symbols $a$ and $p$.
We also use the notation $\nu_i=\nu_{\a_i}=\nu^{d_i}$ for  simple roots
$\a_i$.

Denote by $U_\nu(\g)$ the Hopf algebra, generated by
Chevalley generators
$e_{\a_i}\in U_\nu(\n_{})$,
$e_{-\a_i}=f_{\a_i}\in U_\nu(\n_-)$, $k_{\alpha_i}^{\pm 1}=
\nu_i^{\pm h_{\alpha_i}}\in U_\nu(\H)$,
%$h_{\a_i}\in\H$,
 where $\a_i\in\Pi$, so that
\begin{align*}
[e_{\a_i},e_{-\a_j}]&=\delta_{i,j}\,[h_{\a_i}]_{q_i}=
\delta_{i,j}\,\frac{k_{\a_i}^{}-k_{\a_i}^{-1}}{\nu_i-\nu_i^{-1}}\,,&&\\
k_{\a_i}e_{\pm\a_j}k_{\a_i}^{-1}&=\nu_i^{\pm a_{i,j}}e_{\pm\a_j}=
\nu^{\pm(\a_i,\a_j)}e_{\pm\a_j}\,,&&\\
\sum_{r+s=1-a_{i,j}}(-1)^r& e_{\pm\a_i}^{(r)}e_{\pm\a_j} e_{\pm\a_i}^{(s)}
=0,
\quad i\not= j\,, &%\qquad
\text{where}&\quad e_{\pm\a_i}^{(k)}=
\dfrac{e_{\pm\a_i}^{k}}{[k]_{\nu_i}!}\,,\\ %&&\\
\Delta(e_{\a_i})&=e_{\a_i}\ot 1+ k_{\a_i}\ot e_{\a_i},&
\Delta(e_{-\a_i})&=1\ot e_{-\a_i}+  e_{-\a_i}\ot k_{\a_i}^{-1}\,,\\
\Delta(k_{\a_i})&=k_{\a_i}\ot k_{\a_i},& S(k_{\a_i})&=k_{\a_i}^{-1}\,,\\
S(e_{\a_i})&=-k_{\a_i}^{-1}e_{\a_i},&
S(e_{-\a_i})&=-e_{-\a_i}k_{\a_i}\,.
\end{align*}
The adjoint action \rf{ad} of the Chevalley generators has the form
\begin{equation}
\begin{split}
\label{nu1}
\e_{\a_i}(x)&\equiv\ad_{e_{\a_i}}(x)=
e_{\a_i}x-k_{\a_i}xk_{\a_i}^{-1}e_{\a_i},\\
\e_{-\a_i}(x)&\equiv\ad_{e_{-\a_i}}(x)= [e_{-\a_i},x]\cdot k_{\a_i}\, .
\end{split}
\end{equation}
Let $T_i\equiv T_{s_i}: U_\nu(\g)\to U_\nu(\g)$ be automorphisms of
$U_\nu(\g)$, defined by the relations
\begin{align}%\label{nu2}
%\begin{split}
T_i(e_{\a_i})&=-k_{\a_i}e_{-\a_i},\qquad\notag
T_i(e_{-\a_i})=-e_{\a_i}k_{\a_i}^{-1},&&T_i(k_{\a_j})=k_{s_{\a_i}(\a_j)}\,,\\
T_i(e_{\a_j})&=\sum_{r+s=1-a_{i,j}}(-1)^r\nu_i^r e_{\a_i}^{(r)} e_{\a_j}
e_{\a_i}^{(s)},& i&\not= j, \label{nu0}\\
T_i(e_{-\a_j})&=\sum_{r+s=1-a_{i,j}}(-1)^r\nu_i^{-r} e_{-\a_i}^{(s)}
e_{-\a_j}
e_{-\a_i}^{(r)},& i\not&= j\,.\notag
%T_i(k_{\a_j})&=k_{s_{\a_i}(\a_j)}\,.
%\end{split}
\end{align}
Here we use the Cartan generators $k_\gamma$,
$\gamma\in\Delta$. They are defined by the rules $k_{-\a}=k_\a^{-1}$, and
$k_{\a+\b}=k_\a k_\b$.

In Lusztig's notation \cite{L}, $T_i\equiv T'_{i,+}$. For any $w\in W$ they
 define automorphisms $T_w: U_\nu(\g)\to U_\nu(\g)$, as in Section
\ref{section5.1}.

Denote by $\Dh_\nu$ the localization of the  commutative algebra $U_\nu(\H)$
with respect to the multiplicative set of denominators, generated by
$$ \{[h_\a+k]_{\nu_\a}|\a\in\Delta\,,k\in\ZZ\}\,.$$
Denote by $U'_\nu(\g)$ the extension of $U_\nu(\g)$ by means of $\Dh_\nu$:
$$U'_\nu(\g)=U_\nu(\g)\ot_{U_\nu(\H)}\Dh_\nu\approx
\Dh_\nu\ot_{U_\nu(\H)}U_\nu(\g)\, .$$
As well as in the nondeformed case $(\nu=1)$, there exists an extension $F_{\g,\n}^\nu$
of the algebra $U'_\nu(\g)$ and an element $\P=\P_\n\in F_{\g,\n}^\nu$
 (the extremal projector),
satisfying the conditions $e_{\a_i}\P=\P e_{-\a_i}=0$, $\P^2=\P$, see \cite{KhT}.

In particular, for the algebra $U_q({\mathfrak sl}_2)$, generated by
$e_{\pm\a}$ and $k_\a^{\pm 1}$, we have two projection operators,
$\P=\P_\a[\rho]$, and $\P_-=\P_{-\a}[\rho]$, where
\begin{equation}
\label{nuP6}
\begin{split}
\P_\a[\lambda]&=\sum\limits_{n=0}^\infty
\frac{(-1)^n}{[n]_{\nu_\a}!} \nu_\a^{-\langle h_\a,
\la-\rho\rangle}f_{\a,n}[\lambda]e_{-\a}^ne_\a^n,\\
\P_{-\a}[\lambda]&=\sum\limits_{n=0}^\infty
\frac{(-1)^n}{[n]_{\nu_\a}!} \nu_\a^{\langle h_\a, \la-\rho\rangle}
g_{\a,n}[\lambda]e_{\a}^ne_{-\a}^n
\end{split}
\end{equation}
and
\begin{equation*}
%\label{nuP5}
f_{\a,n}[\lambda]=\prod\limits_{j=1}^n[h_\a+\langle
h_\a,\lambda\rangle+j]_{\nu_\a}^{-1},\qquad
g_{\a,n}[\lambda]=\prod\limits_{j=1}^n[-h_\a+\langle
h_\a,\lambda\rangle+j]_{\nu_\a}^{-1}.
\end{equation*}

Let $\A$ be an associative algebra, which contains $U_\nu(\g)$.
We call $\A$ a $U_\nu(\g)$-admissible algebra if
\begin{itemize}
\item[(a)] there is a subspace $\W\subset \A$, invariant with respect
to the adjoint action of $U_\nu(\g)$, such that the multiplication $m$ in $\A$
induces an isomorphisms of vector spaces
 $${\rm (a1)}\qquad m: U_\nu(\g)\ot \W\to\A\,,\qquad{\rm (a2)}\qquad
 m: \W \ot U_\nu(\g)\to\A\, ;$$
\item[(b)] the adjont action in $\W$ of all real root vectors
 $e_\gamma\in U_\nu(\g)$,
 related to any normal ordering of the root system,
  is locally nilpotent.
The adjoint action of the Cartan subalgebra $U_\nu(\H)$ in $\W$ is semisimple.
% Here $\gamma\in\Delta^{re}$.
\end{itemize}
In particular, $\A$ is isomorphic to $U_\nu(\g)\ot \W$ and to
$\W\ot U_\nu(\g)$ as a
$U_\nu(\g)$-module with respect to the adjoint action.

Denote by $\nnq$ the linear subspace of $U_\nu(\g)$ generated by the elements
$e_{\a_i}$, $\a_i\in\Pi$. Denote by $\nnq_{-}$ the linear subspace
of $U_\nu(\g)$ generated by the elements $e_{-\a_i}$, $i\in\Pi$. Let $U_\nu(\n)$
be the subalgebra of $U_\nu(\g)$, generated by $\nnq$ and
$U_\nu(\n_-)$ the subalgebra of $U_\nu(\g)$, generated by $\nnq_{-}$.

For any $U_\nu(\g)$-admissible algebra $\A$ put
$\A'=\A\ot_{U_\nu(\H)}\Dh_\nu$ and define Mickelsson algebras
$S^{\nnq}(\A)$ and
$Z^{\nnq}(\A)$,
as in  Section \ref{section2.2}:
\begin{align*}
S^{\nnq}(\A)&=\N(\A \nnq)/\A \nnq_{},&
Z^{{\nnq}}(\A)&=\N(\A'\nnq_{})/\A'\nnq_{}
\end{align*}
and the double coset algebra $\Ann$, see  Section \ref{sectiondcs},
$$
\Ann=\n_-\A'\backslash \A'/\A'\n\equiv\A'/\left(\n_-\A'+\A'\n\right)\,,
$$
equipped with the multiplication structure \rf{?3}.

With the conditions of  Section \ref{section2.6}, the Mickelsson algebra
$Z^{\nnq}(\A)$  has distinguished
generators $z_\w$, $z'_\w$,  $\w\in\W$, determined by the relations \rf{r7}
and \rf{r8}.

\subsection{Basic constructions}

Let $\a\in\Pi$ be a simple root, $\n_\a=\CC e_\a$, $\n_{-\a}=\CC e_{-\a}$.
 Let $x\in \A$ be an element of $\A$,  finite with respect to the
adjoint action of $e_{\a}$. Denote by $\q_\a(x)$ the following element of
$\A'/\A'\n_\a$:
\begin{equation}
\label{nu3}
\begin{split}
\q_{\a}(x)&=\sum_{n\geq 0}\frac{(-1)^n}{[n]_{\nu_\a}!}
\left(\hat{k}_\a^{-1}\e_\a\right)^n(x)e_{-\a}^n g_{n,\a}\\
&=\sum_{n\geq 0}\frac{(-1)^n}{[n]_{\nu_\a}!}
\nu_\a^{-n(\hat{h}_\a-n+1)}\e_\a^n(x)e_{-\a}^n g_{n,\a}
\qquad
\text{mod}\quad \A'\n_\a\,,
\end{split}
\end{equation}
where $g_{n,\a}=\left([h_\a]_{\nu_\a}[h_\a-1]_{\nu_\a}\cdots
[h_\a-n+1]_{\nu_\a}\right)^{-1}$.
The assignment \rf{nu3} has the properties
\begin{align}
%\begin{align}
\label{q4}
&(i)& \notag
\q_{\a}^{}(xe_{-\a})&=0\ ,%\notag
\\
&(ii)&\notag
\q_{\a}^{}(x k_\a^{-1})&= \nu_\a^{-2}\q_{\a}^{}(x)k_\a^{-1}
,&
\q_{\a}^{}(k_\a^{-1}x )&= \nu_\a^{-2}k_\a^{-1}\q_{\a}^{}(x)\ ,
\\
&(iii)&\notag
\q_{\a}^{}(x k_\gamma)&= \q_{\a}^{}(x)k_\gamma
,&
\q_{\a}^{}(k_\gamma x )&= k_\gamma\q_{\a}^{}(x),
\quad \text{if}\quad \langle h_\a,\gamma\rangle =0,
\\
&(iv)&\notag
  k_\a^{-1} e_\a\q_{\a}^{}(x)&=
\q_{\a}^{}( k_\a^{-1} e_\a x)=0\ .
%-k\oq_{\a}^{(k-1)}(x)\ .%\notag
\end{align}%\end{align}
We extend the assignment \rf{nu3} to the map
$\q_{\a}: \A'\to \A'/\A'\n_\a$ with the help of the properties (ii) and
(iii). It satisfies the property
\begin{equation}
\label{nu4a}
\q_\a^{}(xd)=\oq_\a^{}(x)\shift_\a(d),\qquad
\q_\a^{}(dx)=\shift_\a(d)\q_\a^{}(x),\qquad
d\in \Dh_\nu\, ,
\end{equation}
where  $\shift_\mu:\Dh_\nu\to\Dh_\nu$, $\mu\in\H^*$ is uniquely characterized
 by the conditions
$$\tau_\a(k_{\gamma})=\nu^{(\mu,\gamma)}k_{\gamma}\,.$$

Due to the property (iv) the map $\q_\a$ defines a map
$\q_\a:{}_{\n_{\a}}\A_{\n_{-\a}}\to Z^{\n_\a}(\A)$.
We have an analog of  Theorem \ref{theorem1}:

\begin{proposition}\label{propnu3}
The map
$$\q_\a:{}_{\n_{\a}}\A_{\n_{-\a}}\to Z^{\n_\a}(\A)$$
is an isomorphism of algebras.
\end{proposition}

In the following we assume that the automorphisms \rf{nu0} admit
extensions $T_i:\A'\to\A'$, which satisfy the braid group relations
\rf{b2}, though, as well as in the classical case,
part of the results below do not depend on this assumption.

Let $w\in W$ be an element of the Weyl group of $\g$, $\a\in\Pi$
 a simple root, such that $l( ws_\a)=l(w)+1$. Set $\gamma=w(\a)$,
$e_{\pm\gamma}=T_w(e_{\pm\a})$ and $T_\gamma= T_w T_\a T_w^{-1}$.
 Denote by $\m=\nnq^w$ the
 linear span of
the vectors $T_w(e_{\a_i})$, $\a_i\in\Pi$ and by $\m^{s_\gamma}$ the space
$T_\gamma(\m)$. Let
 $\q_{\gamma,\m}$ be the linear map $\q_{\gamma,\m}:\A'\to A'/\A'\m$, defined by
the prescription:
\begin{equation*}
%\label{nu4}
\begin{split}
\q_{\gamma,\m}(x)&=\sum_{n\geq 0}\frac{(-1)^n}{[n]_{\nu_\gamma}!}
\left(\hat{k}_\gamma^{-1}\e_\gamma\right)^n(x)e_{-\gamma}^n g_{n,\gamma}\\
&=\sum_{n\geq 0}\frac{(-1)^n}{[n]_{\nu_\gamma}!}
\nu_\gamma^{-n(\hat{h}_\gamma-n+1)}\e_\gamma^n(x)e_{-\gamma}^n g_{n,\gamma}
\qquad
\text{mod}\quad \A'\m\,,
\end{split}
\end{equation*}
for $x$, which are adjoint finite with respect to $e_\gamma$, and then extended
 to $\A'$ by means of the properties, analogous to (i) and (ii) for the map
$\q_\a$.

Here $g_{n,\gamma}$ $=$
 $\left([h_\gamma]_{\nu_\gamma}[h_\gamma-1]_{\nu_\gamma}
\cdots[h_\gamma-n+1]_{\nu_\gamma}\right)^{-1}$.

\begin{proposition}\label{propnu1}${}$
\begin{itemize}
\item[(i)]
We have
$\q_{\gamma,\m}(\A'\m^{s_\gamma})=0$, so $\q_{\gamma,\m}$ defines a map
$\q_{\gamma,\m}:\A'/\A'\m^{s_\gamma}\to \A'/\A'\m^{}$.
\item[(ii)] We have an equality ($\a=w^{-1}(\gamma)$)
\begin{equation}\label{nu5}
\q_{\gamma,\m}=T_w \q_{\a,\n} T_w^{-1}\,.
\end{equation}
\end{itemize}
\end{proposition}

\noindent
Note that  the statement (ii) is  nontrivial for $\nu\not=1$,
since
$$T_w \ad_x T_w^{-1}(y)\not=
\ad _{T_w(x)}(y)\,.$$

Let $\overline{w}=\{w=s_{\a_{i_1}}s_{\a_{i_2}}\cdots s_{\a_{i_n}}\}$ be a reduced
decomposition of the element $w\in W$. Let $\gamma_1,...,\gamma_n$ be a
related sequence of positive roots: $\gamma_1=\a_{i_1},...,\gamma_k=s_{\a_{i_1}}$ $\cdots$
$s_{\a_{i_{k-1}}}(\a_{i_k}),...$. Proposition \ref{propnu1} (i) implies
that there is a well defined map
$$\q_{\overline{w}}:\A'/\A'\n^w\to \A'/\A'\n:\qquad
 \q_{\overline{w}}=\q_{\gamma_1,\n}\q_{\gamma_2,\n^{s_{\gamma_1}}} \cdots
\q_{\gamma_n,\n^{s_{\gamma_{n-1}}\cdots s_{\gamma_1}}}\,.$$
\begin{proposition}\label{propnu2}
Let $\g$ be of finite dimension.
Then for any reduced decomposition $\overline{w}_0$ of   the longest element
 $w_0$ of $W$ the map $\q_{\overline{w}_0}$ sends a vector $\w\in\W$ to the generator
$z'_\w$ of the Mickelsson algebra $\Zp(\A)$.
\end{proposition}

Proposition \ref{propnu2} implies that the maps $q_{\gamma,m}$ satisfy the cocycle
conditions, that is, the maps $\q_{\overline{w}}$ do not depend on a reduced decomposition  $\overline{w}$ of $w\in W$; they can thus be denoted as $q_w$.

Set $\qq_i=\q_{s_{\a_i}}\cdot \T_i:\, \A'/\A'\nnq\to \A'/\A'\nnq$.  They satisfy the braid group relations:
 \begin{equation*}
\underbrace{\qq_i\qq_j\cdots }_{m_{i,j}}=
\underbrace{\qq_j\qq_i\cdots }_{m_{i,j}},\qquad\qquad
 i\not=j\,,
\end{equation*}
and, due to Proposition \ref{propnu3},
 their restriction to  $\Zp(\A)$ are automorphisms of the Mickelsson
algebra, such that
$$\qq_i(dx)=\left(s_{\a_i}\circ d\right)\cdot x, \qquad
\qq_i(xd)=x\cdot\left(s_{\a_i}\circ d\right), \qquad d\in\Dh_\nu,x\in\A'.$$
Here $w\circ d$ is the natural extension of the shifted action of $w\in W$
in $\H^*$, see \rf{circ}, to the automorphism of $\Dh_\nu$, defined by the
conditions
$$ w\circ k_\gamma=k_{w(\gamma)}\cdot q_\gamma^{\langle \gamma, w(\rho-\rho)
\rangle} .$$

\subsection{Some calculations}\label{section7.3}

There is a standard construction of the extension of the Hopf algebra $U_\nu(\g)$
by means of automorphisms $T_i$. Namely, let
$U^W_\nu(\g)$ be the smash product
of $U_\nu(\g)$ and of the algebra, generated by elements $\TT_i^{\pm^1}$, satisfying
the braid group relations
$$ \underbrace{\TT_i\TT_j\cdots }_{m_{i,j}}=
\underbrace{\TT_j\TT_i\cdots }_{m_{i,j}},\qquad\qquad
 i\not=j\,.$$
The cross-product relations are
\begin{equation}\label{nu8}\TT_i g \TT_i^{-1}=T_i(g),\qquad g\in U_\nu(\g).
\end{equation}
Due to coalgebraic properties of Lusztig automorphism, the smash product
$U^W_\nu(\g)$ can be
equipped with a structure of a Hopf algebra, if we put
$$\Delta(\TT_i)=\TT_i\ot \TT_i\cdot \tilde{R}_i^{},$$
where
$$\tilde{R}_i=\exp_{\nu_i^{-2}}\left((\nu_i-\nu_i^{-1})e_{-\a_i}\ot e_{\a_i}\right)=
\sum_{n\geq0} \frac{(\nu_i-\nu_i^{-1})^n}{(n)_{\nu_i^{-2}}!}e_{-\a_i}^n\ot e_{\a_i}^n.$$
In the same way we  extend  the algebra $\A'$ to the cross-product $\A^W$,
using  the  relations \rf{nu8}. Since $U_\nu^W$ is a Hopf algebra, the adjoint action
$\hat{\TT}_i$ of $\TT_i$ in $\A^W$ is well defined. It respects the subalgebra
$\A'\subset \A^W$: $\hat{\TT}_i(\A')\subset \A'$.
The following statement is nontrivial for $\nu\not=1$ and important for calculations of the maps
$\qq_i$.
\begin{proposition}\label{propnu4}
For any $x\in\A'$ we have
\begin{equation*}
\qq_i(x)=\q_i({\rm T}_ix{\rm T}_i^{-1})=\q_i( \hat{{\rm T}}_i(x)).
\end{equation*}
\end{proposition}

Now we describe the squares of the automorphisms $\qq_i: \Zp(\A)\to\Zp(\A)$.
Assume that elements $\w_{m,j}\in\W$, $m\in\ZZ_{\geq0}$, $j=0,1,...,m$
 form a finite-dimensional representation of the algebra
$U_\nu({\mathfrak sl}_2)$, generated by $e_{\pm\a_i}$ and $k_{\a_i}^{\pm 1}$,
with respect to the adjoint action, so that we have:
 \begin{equation}\label{vmj}
\e_{\a_i}^{j+1}\left(\w_{m,j}\right)=\e_{-\a_i}^{m-j+1}\left(\w_{m,j}\right)
=0, \qquad \hat{h}_{\a_i}\left(\w_{m,j}\right)=
(m-2j)\w_{m,j}\, .
\end{equation}
In particular, $\w_m=\w_{m,0}$ is the highest weight vector of
 this  representation, and $\w_{m,j}=\e_{-\a_i}^{(j)}\left(\w_m\right)$.

\begin{proposition}
 Assume that $\w_{m,j}\in\W$ satisfy \rf{vmj}. Then
\begin{equation}\label{qsquare}
\qq_i^2(\w_{m,j})=\nu^{-j(m-j+1)-(j+1)(m-j)}\cdot [h_{\a_i}+1]_{\nu_i}^{-1}\cdot \hat{{\rm T}}_i^2(\w_{m,j})\cdot [h_{\a_i}+1]_{\nu_i}\,.
\end{equation}
\end{proposition}

The property \rf{qsquare} simplifies under natural assumptions on operators
$\hat{{\rm T}}_i:\W\to\W$.  Namely,
suppose that the operators $\hat{{\rm T}}_i:\W\to\W$ satisfy properties of Lusztig symmetries $T'_{i,+}$, that is,  (see \cite{L}, 5.2.2)
\begin{equation}\label{Tmj}
\hat{{\rm T}}_i(\w_{m,j})=\nu^{j(m+1-j)}\w_{m,m-j}\,.
\end{equation}
\begin{corollary}  With the  conditions \rf{Tmj}
 for any $x\in\Zp(A)$ we have
\begin{equation}\label{qsquare2}
\qq_i^2(x)= [h_{\a_i}+1]_{\nu_i}^{-1}\cdot x\cdot [h_{\a_i}+1]_{\nu_i}\,.
\end{equation}
\end{corollary}

Keep the notation \rf{r10} of  Section \ref{section2.6}. For a real root
$\a\in\Delta^{re}$, set
\begin{equation*}
\begin{split}%\label{r11a}
\bar{f}_{n,\a}^{(2)}[\mu]&=
\nu_\a^{-2n}\left(k_\a^{(1)}\right)^{-n}
\prod_{k=1}^{n}
\left([\Hat{h}_\a^{(2)}-h^{(1)}_\a+\langle h_\a,\mu\rangle+k
]_{\nu_\a}\right)^{-1},\\
\bar{g}_{n,\a}^{(2)}[\mu]&=
\left(k_\a^{(1)}\right)^{n}
\prod_{k=1}^{n}
\left([-\Hat{h}_\a^{(2)}+h^{(1)}_\a+\langle h_\a,\mu\rangle+k]_{\nu_\a}
\right)^{-1}.
\end{split}
\end{equation*}
Here $\Hat{h}_\a^{(2)}=\ad_{h_\a}^{(2)}$ is the adjoint action of $h_\a$ in
$\W$, $h^{(1)}_\a$ is the operator of multiplication by $h_\a$ in $\Dh$.
{}For $\mu\in\H^*$
define  operators
 $\BB_{\a}^{(2)}[\mu]:\Dh\ot\W\to \Dh\ot\W$ and
$\B_{-\a}^{(2)}[\mu]:\Dh\ot\W\to \Dh\ot\W$
 by the relations \rf{r12}:
\begin{equation*}%\label{r12}
\begin{split}
\BB_{\a}^{(2)}[\mu]&=\sum_{n=0}^\infty
\frac{(-1)^n}{[n]_{\nu_\a}!}\bar{f}_{n,\a}^{(2)}[\mu]\big(\e_{-\a}^{(2)}\big)^n
\big(\e_\a^{(2)}\big)^n \!, \\
\B_{-\a}^{(2)}[\mu]&=\sum_{n=0}^\infty
\frac{(-1)^n}{[n]_{\nu_\a}!}\bar{g}_{n,\a}^{(2)}[\mu]\big(\e_{\a}^{(2)}\big)^n
\big(\e_{-\a}^{(2)}\big)^n \!.
\end{split}
\end{equation*}

\begin{proposition}
For any $\w\in\W$ we have
$$
\qq_i(z'_\w)=
z'_{\B_{-\a_i}^{(2)}[\rho]\left(1\ot\T_i(\w)\right)},\qquad
\qq_i(z_{\w})= z_{\BB_{\a_i}^{(2)}[-\rho](1\ot \T_i(\w))}\,.
$$
\end{proposition}

\subsection{Another adjoint action}

In this Section we sketch the modification
of the above constructions for the second adjoint action.

Let $U^{op}_\nu(\g)$ be the Hopf algebra $U_\nu(\g)$, described in Section
\ref{section7.1}, with the same multiplication and opposite comultiplication,
\begin{align*}
\Delta^{op}(e_{\a_i})&= k_{\a_i}\ot e_{\a_i}+e_{\a_i}\ot 1,&
\Delta^{op}(e_{-\a_i})&= e_{-\a_i}\ot k_{\a_i}^{-1}+1\ot e_{-\a_i}\, .
\end{align*}
The adjoint action \rf{ad} for the algebra $U^{op}_\nu(\g)$
looks slightly different:
\begin{equation}
\begin{split}
\label{nu24}
\e_{\a_i}(x)&\equiv\ad_{e_{\a_i}}(x)= [e_{\a_i},x]\cdot k_{\a_i}^{-1}\,,\\
\e_{-\a_i}(x)&\equiv\ad_{e_{-\a_i}}(x)=
e_{-\a_i}x-k_{\a_i}^{-1}xk_{\a_i}e_{-\a_i}\ .
\end{split}
\end{equation}
Let $\A$ be a $U_\nu^{op}(\g)$-admissible algebra.
 We define Zhelobenko
operators, starting from the assignment
\begin{equation*}
\begin{split}
\q_{\a}(x)&=\sum_{n\geq 0}\frac{(-1)^n}{[n]_{\nu_\a}!}
\left(\e_\a\right)^n(x)\left(e_{-\a}^n{k}_\a^{} \right)g_{n,\a}\\
&=\sum_{n\geq 0}\frac{(-1)^n}{[n]_{\nu_\a}!}
\e_\a^n(x)e_{-\a}^n \nu_\a^{n({h}_\a-n+1)} g_{n,\a}
\qquad
\text{mod}\quad \A'\n_\a\, .
\end{split}
\end{equation*}
Corresponding maps $\q_{\a,\m}$ satisfy the cocycle conditions, such that
the operators
$$\qq_i=\q_{\a_i,\n}\cdot T_i^{-1}$$ are automorphisms of the Mickelsson algebra
$\Zp(\A)$, satisfying the braid group relations.

\subsection{Relations to dynamical Weyl group}

Let $\W$ be a $U_\nu(\g)$-module algebra with a locally nilpotent
action of real root vectors of $U_\nu(\g)$. It means that $\W$ is
an associative algebra and a $U_\nu(\g)$-module, such that for any
$g\in U_\nu(\g)$, $\w_1$, $\w_2\in\W$, we have the equality
$$\hat{g}(\w_1\cdot \w_2)=
\sum\nolimits_i \hat{g}'_i(\w_1)\cdot\hat{g}''_i(\w_2)\,.$$
Here $\Delta(g)=\sum\nolimits_i g'_i\ot g''_i$, and $\hat{g}(\w)$
is an action of $g\in U_\nu(\g)$ on $\w\in\W$. Suppose that $\W$ is
equipped with an action of operators $\hat{T}_i$, which are automorphisms
of the algebra, satisfy the braid group
relations, and form an equivariant structure with respect to the operators
\rf{nu0}, that is, for any $\w\in\W$, $g\in U_\nu(\g)$,
$$\hat{T}_i\left(\hat{g}(\w)\right)=T_i(g)\hat{T}_i(\w)\,,$$
where $T_i(g)$ means the application of operators \rf{nu0}.

These conditions imply that $\W$ is a module algebra over $U_\nu^W(\g)$.
So we have a well defined smash product $U_\nu^W(\g)\ltimes
\W$, which contains $U_\nu(\g)\ltimes\W$ and the elements $\TT_i$.
The automorphisms $T_i$: $U_\nu^W(\g)\ltimes\W\to U_\nu^W(\g)\ltimes\W$,
given as $T_i(x)=\TT_i x\TT_i^{-1}$, preserve the subalgebra
$U_\nu(\g)\ltimes\W$. Moreover, the restriction of the adjoint action of
$\TT_i$ on $\W$ coincides with $\hat{T}_i$:
$\hat{\TT}_i|_{\W}=\hat{T}_i$. Thus the smash product
$\A=U_\nu(\g)\ltimes\W$ is a $U_\nu(\g)$-admissible algebra.

Elements of the right $\Zp(\A)$-module $\Phi_\lambda(\A)$, defined in
\rf{b23a}, are intertwining operators $\Phi_\lambda^\w: M_\lambda\to
\W\ot M_\lambda$. Operators $\qq_i$ give rise to the operators $\qq_{i,\la}$
of the dynamical Weyl group
$$\qq_{i,\la}(\Phi_\lambda^\w)=
\Phi_{s_{\a_i}}^{\check{\PP}_{-\a_i}[\la+\rho](\hat{T}_i(\w))}\,$$
and
\begin{equation*}
\qq_{w,\lambda}(\Phi_\lambda^v)=\Phi_{w\circ\lambda}^{\Hat{\PP}_{-\gamma_1}[\lambda+\rho]\cdots
\Hat{\PP}_{-\gamma_n}[\lambda+\rho](\hat{T}_w(\w))}\,,
\end{equation*}
where $\gamma_1,...,\gamma_n$ is the sequence of positive roots,
 attached to a reduced decomposition $w=s_{\a_{i_1}}\cdots s_{\a_{i_1}}$
by the standard rule $\gamma_1=\a_{i_1}$, $\gamma_2=s_{\a_{i_1}}(\a_{i_2})$,
$\ldots$; here $\check{\PP}_{-\gamma_k}[\la+\rho]$ is the
adjoint action of the operators \rf{nuP6}, $e_{\pm\gamma_k}=
T_{\a_{i_1}}\cdots T_{\a_{i_{k-1}}}(e_{\pm\a_{i_k}})$ are the Cartan-Weyl
generators.

\section{Concluding remarks}

We conclude with  remarks on the assumptions on a $\g$-admissible algebra
$\A$, used in the paper. They are listed in  Section \ref{section2.1}.

The assumption (a) requires an existence of an $\ad$-invariant subspace
$\W\subset\A$, such that the multiplication $m$ in $\A$ induces
isomorphisms of vector spaces
$${\rm (a1)}\qquad m: U_\nu(\g)\ot \W\to\A\,,\qquad
 {\rm (a2)}\qquad m: \W \ot U_\nu(\g)\to\A\,.$$
Assumptions (a1) and (a2) are not of equal use. For the Mickelsson algebra
$\Zp(\A)$ we need the condition (a1) only when we use generators $z_\w$, that
is, in Proposition \ref{propr3}, Theorem \ref{theorem0},
 Corollary \ref{cor1}, in Section
\ref{section5.2} and in corresponding statements of  Section \ref{section8}.
On the contrary, the construction of the Zhelobenko operators for the
algebra  $\Zp(\A)$ requires the condition (a1) from the very beginning.
The condition (a2) is also necessary for the existence of the generators
$z'_\w$.

{}For the algebra $\Zm(\A)$ the situation is opposite. We need the condition
(a1) for the construction of the Zhelobenko maps and generators $\tilde{z}'_\w$,
while the condition (a2) is related only to the generators $\tilde{z}_\w$.
Both conditions (a1) and (a2) are satisfied for basic examples, listed in 
Section \ref{section2.1}.

%%%%%%%%%%%%%%%%%newsection%%%%%%%%%%%%%%%%%%%%%%%%%%
The condition (b) requires a local nilpotency of the adjoint action of real
root vectors in $\W$. It always takes place if the space
$\W$ is a sum of integrable representations or an affinization $V(z)$ of
a locally finite representation of an affine algebra $U_\nu'(\g)$ with the
gradation element dropped.

In the latter case the generators $z_\w$ or $z'_\w$  do not formally
exist, since neither the highest weight (HW) from
 Section \ref{sectiondcs} nor the lowest weight (LW) condition from
 Section \ref{sectionzn1} is satisfied.
 Nevertheless, in this case the generators of the Mickelsson algebra
exist as formal series and could be used with a proper attention to
convergences.

\medskip
%\section*{Acknowledgement}
\begin{center}
\textsc{Acknowledgement}
\end{center}
\medskip
A big part of the work was done during the visit of the first author
to CRM and Marseille University.
Authors thank CRM and  Marseille University for the hospitality
and stimulating scientific atmosphere. S.Kh was supported
by the RFBR grant 05-01-01086 and RFBR grant for the support of
scientific schools NSh-8065.2006.2
. S.Kh. and O.O.
were supported by the ANR project GIMP No.
ANR-05-BLAN-0029-01.

%%%%%%%%%%%%% END OF PROOFS %%%%%%%%%%%%%%%%

\end{document}